\documentclass[dvips,
a4paper,twoside,
11pt
]{report} 

\usepackage{tocbibind} 
\usepackage{verbatim} 
\usepackage{comment} 
\usepackage{amsmath, amsthm, amssymb}

\usepackage{helvet}

\usepackage{setspace} 
\usepackage{color}  
\usepackage{epsfig}  
\usepackage{floatflt} 
\usepackage{nomencl} 
\usepackage{tabularx} 
\usepackage{afterpage} 
\usepackage{layout} 

\usepackage[a4paper,twoside,includeheadfoot,inner=4.1cm,outer=2.4cm,top=2.4cm,bottom=2.4cm,marginparwidth=90pt]{geometry} 

\newcommand{\newchapter}[2]                   
{                                                       
\chapter {#1}                                   
\markright{\color{draftgrey}{\chaptername  \,\thechapter :\;{#2}} 
} 
}                                                       

\newenvironment{draftcomment}{\color{draftgrey}\small}

\newtheorem{theorem}{Theorem}[chapter]

\newtheorem{corollary}[theorem]{Corollary}

\newtheorem{lemma}[theorem]{Lemma}

\DeclareMathOperator{\trace}{trace}
\DeclareMathOperator*{\osc}{osc}
\DeclareMathOperator{\dist}{dist}
\DeclareMathOperator{\diam}{diam}
\DeclareMathOperator{\diag}{diag}
\DeclareMathOperator{\graph}{graph}

\newcommand{\pd}[2]{\dfrac{\partial#1}{\partial#2}}
\newcommand{\pdd}[3]{\dfrac{\partial^2#1}{\partial#2 \partial#3}}
\newcommand{\bigR}{{\mathbb R}}
\newcommand{\bigS}{{\mathbb S}}
\newcommand{\pD}[2]{\dfrac{\partial^2#1}{\partial#2^2}}
\newcommand{\pDDD}[2]{\dfrac{\partial^3#1}{\partial#2^3}}
\newcommand{\smallplus}{\mspace{-1mu}+\mspace{-1mu}}
\newcommand{\smallminus}{\mspace{-1mu}-\mspace{-1mu}}
\newcommand{\smallequals}{\mspace{-1mu}=\mspace{-1mu}}
\newcommand{\nabsq}[1]{\left|\nabla #1 \right|^2}

\newcommand{\erf}[1]{{\mathcal{E}}{\rm{rf}}\left({#1}\right)}
\newcommand{\inverf}[1]{{\mathcal{E}}{\rm{rf}}^{-1}\left({#1}\right)}

\DeclareMathOperator{\Lip}{Lip}

\DeclareMathOperator{\diverg}{div}

\renewcommand{\proof}{\par\smallskip {\noindent\textbf{Proof:}}\quad} 
\newcommand{\remark}{\par\smallskip{\noindent\textbf{Remark:}}\quad}
\newcommand{\margincomment}[1]

\newcommand{\halmos}{\qedsymbol}

\newcommand{\limbda}{d}

\newcommand{\longpage}{\enlargethispage{\baselineskip}}

\makeatletter
\renewcommand{\@pnumwidth}{1.68em}\renewcommand{\@tocrmarg}{2.7em}
\makeatother

\makeglossary

\begin{document}

\definecolor{draftgrey}{gray}{0.6}

\begin{titlepage}

\thispagestyle{empty}
\hbox{ } {\Huge\bf{
\vskip5cm
\centerline{
Parabolic Equations}
\centerline{
with }
\centerline{
Continuous Initial Data}  }}

\vskip2cm

\centerline{\Large{Julie Clutterbuck}}
\vskip3cm

\smallskip
\bigskip
\vskip2cm

\centerline{April 2004}

\vskip2cm
\bigskip

{\sc{ \centerline{A thesis submitted for the degree of } 
\centerline{Doctor of Philosophy at the Australian National University}}}

\end{titlepage}

\chapter*{Declaration}
The work in this thesis is my own except where otherwise stated.
\vskip5cm
\hskip 5cm 
Julie Clutterbuck

\chapter*{Thanks}

My fellow students have been an unfailing source of good conversation and interesting mathematics:  thankyou Nick Dungey, Andreas Axelsson, Denis Labutin, Jason Sharples and others at the Australian National University.

My supervisor Ben Andrews has been surprisingly patient over the years, and an excellent and insightful guide.

At the Mathematical Sciences Institute,  the administrative and academic staff have been friendly and helpful.    Alan McIntosh has been generous with his good advice (which I ignored) and his friendship (which I appreciate).

I thank my dear friend Maria Athanassenas for her support as my mentor.

My family have been wonderful:  thankyou to my lovely parents, Sue and Neil, my grandmother Marjorie Anderson, and the rest of them.  My friends have remained my friends which is more than I could ask for:    thankyou Annabel Barbara, Nathalie Jitnah, Rohan Baxter, Joseph Baxter, Dan Rosauer, Allison Foster,  and Nick Lake.  

And above all, thanks to my beloved Brendan.


\chapter*{Abstract}

The aim of this thesis is to derive new gradient estimates for parabolic equations.   The gradient estimates found are independent of the regularity of the initial data.  
This allows us to prove the existence of solutions to problems that have non-smooth, continuous initial data.  We include existence proofs for problems with both Neumann and Dirichlet boundary data.

The class of equations studied is modelled on 
mean curvature flow for graphs.  It includes anisotropic mean curvature flow, and other operators that  have no uniform non-degeneracy bound.

We arrive at similar estimates by three different paths:  a 'double coordinate' approach, an approach examining the intersections of a solution and a given barrier, and a classical geometrical approach.

\tableofcontents


\newchapter{Preface}{Preface}
\subsection*{Mean curvature flow}
Let $\lbrace M_t \rbrace$ be a family of hypersurfaces, each smoothly embedded  in $\bigR^{n+1}$ and indexed by $t\in[0,T]$.   We say   $M_t$ is moving by mean curvature flow when
\begin{equation*}
\pd x t = \mathbf{H},
\end{equation*}
where $\mathbf{H}$ is the mean curvature vector at $x\in M_t$.   

In the last twenty-five years, this flow  has been the subject of concerted study, as have other geometric flows such as the Ricci flow and Gauss curvature flow.  Notable results have included Grayson's proof that the curve shortening flow 
(which is mean curvature flow, reduced to one space dimension)
  shrinks embedded curves to a spherical point \cite{grayson:plane-curves} and Huisken's proof that convex surfaces become spherical under mean curvature flow \cite{huisken:spheres}.     %

If we observe that ${\mathbf{H}}=\Delta_{M_t} x$, where $\Delta_{M_t}$ is the Laplace-Beltrami operator on the manifold $M_t$, then it seems natural 
to consider mean curvature flow as the heat flow for manifolds.   

Unlike classical heat flow, this is a nonlinear operator, but it still has some of the same attributes; in particular, this flow exhibits a \emph{smoothing property}.     In \cite{eh:interior}, Ecker and Huisken showed that if the initial surface is given by a locally Lipschitz graph, then there exists a smooth solution for positive times.     In this thesis, this result is extended to non-Lipschitz initial conditions. 

\hfill %
 \medskip  %

\medskip %

\subsection*{Mean curvature flow and parabolic differential equations}
Chapter \ref{curvature flow chapter} is a short introduction to parabolic differential operators and some key results in mean curvature flow.

If $M_t$ is locally represented as $\graph u $ for some $u:\Omega\times[0,T]\rightarrow\bigR$, then $u$ will satisfy the parabolic equation
\begin{equation*}
\pd u t =\sqrt{1+|Du|^2}\text{\rm{div}}\left(\frac{Du}{\sqrt{1+|Du|^2}}\right)= \left(\delta^{ij}-\frac{D_i uD_ju}{1+|Du|^2}\right)D_{ij}u.
\end{equation*}
This places mean curvature flow in the setting of classical parabolic partial differential equations, the framework for most of this thesis. 

Motivated by this setting, we examine other parabolic operators that have similar diffusion properties to   mean curvature flow.  This includes anisotropic mean curvature flow,   a generalization of mean curvature flow arising in many physical applications.

\subsection*{Gradient estimates using a `double coordinates' approach}
In this thesis three distinct methods are used for deriving gradient estimates.  The first of these is introduced in Chapter \ref{onedimensionalchapter}.  It  originated with Kru{\v{z}}kov in \cite{kruzhkov:nonlinear}.   Given a parabolic equation in one space dimension 
$$u_t=F\left(u_{xx},u_x,u,x,t\right)$$
one can form an evolution equation for the difference $w(x,y,t)=u(x,t)-u(y,t)$,
$$w_t=F\left(u_{xx},u_x,u(x),x,t\right)-F\left(u_{yy},u_y,u(y),y,t\right).$$
Under favourable conditions on $F$, one can use the maximum principle and an appropriate barrier to find estimates for $w$ that depend on $|x-y|$.   Letting $|x-y|\rightarrow 0$ will then give a gradient estimate for $u$.  

In this thesis,   Kru{\v{z}}kov's method is extended by making use of the full Hessian 
\begin{equation*}
[D^2w]
=\begin{bmatrix} w_{xx} & w_{xy} \\ w_{yx} & w_{yy}\end{bmatrix} \end{equation*}
rather that just the diagonal elements $w_{xx}$ and $w_{yy}$.  In the one-dimensional case this is of little importance, but in the higher-dimensional case it gives us greater scope to choose barriers. 

In Chapter \ref{onedimensionalchapter},  we also describe a barrier which begins with an unbounded gradient, but instantly becomes smooth.  

Such barriers allow estimates that are independent of initial gradient bounds, and that therefore may be used to prove the existence of solutions to parabolic equations with continuous initial data.    This is one of the main features of the gradient estimates in this thesis.

The `double coordinate' method is extended to higher dimensions in Chapter \ref{higher dimensional chapter} for a class of operators that have similar diffusion to  the mean curvature flow.  

As this class includes anisotropic curvature flow %
,  this is a significant improvement to the existing regularity theory.  

Gradient estimates are found for both   entire periodic solutions and boundary value problems.    

\subsection*{Existence results}

The gradient estimates of Chapters  \ref{onedimensionalchapter} and \ref{higher dimensional chapter}  may be used to extend standard existence results.

We do this for a class of parabolic equations in the one dimensional case in \mbox{Chapter \nolinebreak[1] \ref{short-time chapter};} for mean curvature flow with Neumann boundary conditions in Chapter \ref{neumann chapter}; and for mean curvature flow with Dirichlet boundary conditions in Chapter \ref{dirichlet chapter}.
 
\subsection*{Gradient estimates by counting intersections}

The second technique for finding gradient estimates is found in Chapter \ref{zero counting chapter}, and it involves  
examination of the intersections between a given solution $u$ and a barrier $\varphi$.   

In \cite{ang:zeroset}, Angenent proved that the number of points in the \emph{zero set} --- the set where $u(x,t)=0$ ---  of a solution to a parabolic equation in one space dimension is non-increasing.  

This is applied to the difference $u-\varphi$.
  The intersections of $u$ and $\varphi$ are the zeroes of  $u-\varphi$,  and so Angenent's results allow us to show that $u$ and $\varphi$ do not develop new intersections as they evolve.

When two functions intersect only once,  the gradient of one of them dominates the gradient of the other at that point.    Tailoring the barriers gives us estimates for $u_x(x,t)$ in terms of the height of $u(x,t)$, the time $t$, and (in the case of bounded domains) the distance of $x$ from the boundary.  Again, there is no dependence on an initial gradient estimate.

As Angenent's results are limited to equations in one space dimension ---  it is difficult to imagine what a generalization of these results to higher dimensions would look like --- this technique applies only to parabolic operators in one space dimension.

These methods will apply to a wide range of parabolic operators, provided that suitable barriers exist. In particular, we can apply this method to the class of operators studied in   Chapter \ref{short-time chapter}.

\subsection*{Gradient estimates using a geometric approach}

The third method for finding gradient estimates, in  Chapter \ref{newwork chapter}, is a rather geometrical approach found in the classic Ecker--Huisken curvature flow papers \cite{eh:mean},  \cite{eh:interior},  \cite{huisken:spheres}, \cite{hu:boundary} and \cite{hu:local}.
    
 A maximum principle is applied to the difference $Z=v-\varphi$, where $v$ is a ``gradient function'' (for example, $v=\sqrt{1+|Du|^2})$, and $\varphi$ is a barrier.    In contrast to the earlier two approaches, this creates a direct estimate for the gradient $Du$ itself, rather than for    the difference $u(x)-u(y)$.

We apply this to the mean curvature flow (re-creating some of the results obtained in earlier chapters) and also to the anisotropic mean curvature flow, under some restrictions on the degree of anisotropy allowed.    Results for entire periodic solutions and strictly interior results are found in both cases.  The estimates found are again independent of initial gradient bounds, but dependent on the height.  

\subsection*{Appendices}

An inventory of standard results for parabolic partial differential equations, relevant to the existence results in Chapters \ref{short-time chapter}, \ref{neumann chapter} and \ref{dirichlet chapter}, is included for the convenience of the reader.  There is also a nomenclature listing.

\newchapter{Mean curvature flow and parabolic equations}{Mean curvature flow and parabolic equations}
\label{curvature flow chapter}

\section{Parabolic equations}
An operator $\boldsymbol{P}:\bigR\times\ {\mathbb S}^{n}\times\bigR^n\times\bigR\times\Omega\rightarrow \bigR\times[0,T]$ \label{page where parabolic operator is defined}
is considered parabolic on a domain when 
$$\boldsymbol{P}(z+\sigma,r+\eta,p,q,x,t)>\boldsymbol{P}(z,r,p,q,x,t)$$
 for any positive definite $\eta\in{\mathbb S}^n$, positive number $\sigma$, and any $(z,r,p,q,x,t)$ in the domain.  Here, ${\mathbb S}^n$ is the set of $n\times n$ symmetric matrices. 

In this thesis we look only at operators of the form \begin{equation*} \boldsymbol{P}u=\boldsymbol{P}(-u_t,D^2u,Du,u,x,t)=-u_t+F(D^2u,Du,u,x,t).\end{equation*}   
If $F$ is differentiable with respect to the first variable,  then $\boldsymbol{P}$ will be parabolic if the %
matrix of derivatives $F_r=\left[\pd F {r_{ij}}\right]$ is positive definite.

If we can write $F(r,p,q,x,t)=a^{ij}(p,q,x,t)r_{ij}+b(p,q,x,t)$ for some symmetric $a:\bigR^n\times\bigR\times\Omega\times[0,T]\rightarrow \bigS^{n}$, then we call the operator  {\emph{quasilinear}}.    It is to quasilinear operators that we will pay most attention   in the following pages.

In this case, $F_r=[a^{ij}]$ and so the operator is  parabolic on ${\cal {S}} $  if $[a^{ij}(p,q,x,t)]$ is positive definite for all $(p,q,x,t)\in {{\cal {S}}}$,  where 
 ${\cal{S}}$ is a subset of $\bigR^n\times\bigR\times\Omega\times[0,T]$.   
 
We denote the maximum and minimum eigenvalues of $[a^{ij}(p,q,x,t)]$ (or $F_r$) by $\lambda(p,q,x,t)$ and $\Lambda(p,q,x,t)$.  

If the ratio $\Lambda/\lambda$ is bounded on ${\cal{S}}$, then $\boldsymbol{P}$ is called {\emph{uniformly parabolic}} on ${\cal{S}}$.

An operator is parabolic with respect to a function $u$ when $\boldsymbol{P}(D^2u,Du,u,x,t)$ is parabolic.

When $\lambda\not>0$ --- for example, when $a^{ij}(Du,u,x,t)= |Du|^{p-1}\delta^{ij}$, as in the parabolic $p$-Laplacian equation --- such an operator is called {\emph{degenerate}}.

The key to much of the theory used here is the \emph{Comparison Principle}.   As presented in \cite{li:parabolic}:
\begin{theorem}[Quasilinear comparison principle 1] \label{comparison principle}
Suppose that $\boldsymbol{P}$ is a quasilinear parabolic operator 
$$\boldsymbol{P}u=-u_t+a^{ij}(Du,u,x,t)D_{ij}u+b(Du,u,x,t).$$ 
Let $u$ and $v$ be in $C^{2,1}(\overline\Omega\setminus{\cal P}\Omega)\cap C(\overline\Omega)$ and let $\boldsymbol{P}$ be parabolic with respect to either $u$ or $v$.  Then if $\boldsymbol{P}u>\boldsymbol{P}v$ in $\overline\Omega\setminus{\cal P }\Omega$ and if $u<v$ on ${\cal P}$, then $u<v$ in $\overline\Omega$.
\end{theorem}
Here ${\cal P}$ denotes the \emph{parabolic boundary} 
$$ {\cal P}\left(\Omega\times[0,T]\right):=\Omega\times\lbrace 0\rbrace \cup \partial\Omega\times[0,T].$$ 
The proof of this theorem is simple, and is an excellent illustration of later arguments.
\proof   Suppose that there is an  interior point $x_0$ at time $t_0>0$ where $u=v$ for the first time.  Since this is an internal maximum of $w=u-v$,   $Dw=Du(x_0,t_0)-Dv(x_0,t_0)=0$ and $D_{ij}w=D_{ij}u-D_{ij}v$ must be negative semi-definite.    Now,
\begin{align*}
w_t&=u_t-v_t \\
&=-\boldsymbol{P}u+ a^{ij}(Du,u,x,t)D_{ij}u +b(Du,u,x,t) 
\\&\phantom{spacespacespace} + \boldsymbol{P}v-a^{ij}(Dv,v,x,t)D_{ij}v-b(Dv,v,x,t) \\
&< a^{ij}(Du,u,x,t)D_{ij} w.
\end{align*}  
However, as this is the first such maximum we must have $w_t\ge 0$, which give us a contradiction. 
It follows that $u<v$.  \halmos

We also use the following form, again as in \cite{li:parabolic} :
\begin{theorem}[Quasilinear comparison principle 2] \label{comparison 1} Suppose that $\boldsymbol{P}$ is a quasilinear parabolic operator 
$$\boldsymbol{P}u=-u_t+a^{ij}(Du,x,t)D_{ij}u+b(Du,u,x,t),$$ 
and  that there is an increasing function $k(M)$ such that $b(p,q,x,t)+k(M)q$ is a decreasing function of $q$ on 
$\Omega\times[-M,M]\times{\mathbb R}^n$ for any $M>0$.  If $u$ and $v$ are functions in $C^{2,1}(\Omega)\cap C(\overline\Omega)$ with $\boldsymbol{P}u\ge \boldsymbol{P}v$ in $\overline\Omega\setminus{\cal P}\Omega$ and $u\le v$ on ${\cal P}\Omega$, and if $\boldsymbol{P}$ is parabolic with respect to  either $u$ or $v$, then $u\le v$ in $\overline\Omega$. 
\end{theorem}

\section{Mean curvature flow}
If our family of hypersurfaces $M_t$ is also a family of embeddings ${\mathbf F}_t:M^n\rightarrow \bigR^{n+1}$, then we can   write mean curvature flow as 
\begin{equation}
\pd { } t {{\mathbf{F}}_t} =\mathbf{H},  \label{mcf in geometric setting}
\end{equation}
where $\mathbf{H}$ is the mean curvature vector at ${\mathbf F}_t(x)\in M_t$.  In the case that $M_t$ can be written locally as a graph over a set $\Omega\in\bigR^n$, we write ${\mathbf F}_t(x)=(x,u(x,t))$, and can calculate geometric quantities such as the upwards unit normal  
\begin{equation*}  \label{nu} 
\nu= \frac{\left(-Du,1\right)}{\sqrt{1+|Du|^2}},
\end{equation*}
the metric on the surface  
\begin{equation*} \label{g_ij}
g_{ij}=  \delta_{ij}+ {D_i u D_j u},
\end{equation*}\nomenclature[gij]{$g_{ij}$}{metric on a Riemannian manifold\refpage   }
the second fundamental form   \label{h_ij}
\begin{equation*} 
h_{ij}=-\frac{D_{ij}u}{\sqrt{1+|Du|^2}},
\end{equation*} \nomenclature[h1ij]{$h_{ij}$ }{ second fundamental form\refpage}
and the mean curvature  \label{H}
\begin{equation*}
H=g^{ij}h_{ij}= -\diverg \left(\frac{ Du}{\sqrt{1+|Du|^2}}\right).
\end{equation*}  \nomenclature[H2]{$H$ }{mean curvature\refpage}

The mean curvature vector is $\mathbf{H}=H\nu$, and so (if we remove movement tangential to the surface) mean curvature flow for graphs is given by 
\begin{equation*}  \label{mcf for graphs}
\pd u t = {\sqrt{1+|Du|^2}}\diverg \left(\frac{ Du}{\sqrt{1+|Du|^2}}\right)= \left(\delta_{ij}-\frac{D_iu D_ju}{ 1+|Du|^2}\right)D_{ij}u.
\end{equation*}
In the case when $n=1$, this reduces to curve-shortening flow 
\begin{equation*}
\pd u t = \frac{u_{xx}}{1+u_x^2}. 
\end{equation*}

With reference to the previous section,  note that the largest and smallest   eigenvalues for mean curvature flow for graphs are $\Lambda=1$ and $\lambda=(1+|p|^2)^{-1}$,  so it will be uniformly parabolic only when the gradient is bounded.  

Whether studied in a geometric setting as \eqref{mcf in geometric setting}, or as  a special case of a quasilinear parabolic differential equation as \eqref{mcf for graphs}, the comparison principle has been crucial.  Applied to mean curvature flow, it gives:
\begin{theorem} Let $M_t$ and $M'_t$ be two smooth compact surfaces moving under mean curvature flow.  If they are disjoint at the initial time, they are disjoint at later times.
\end{theorem}
We can also make similar comparisons between surfaces with boundaries, and between other quantities (such as the \emph{gradient function} $v=\sqrt{1+|Du|^2}$ ).    The following theorem from \cite{eh:interior} is one such result.

\begin{theorem}[Interior gradient estimate] 
Suppose that $u$ satisfies the mean curvature flow equation \eqref{mcf for graphs} on a cylinder $B_R(y_0)\times[0,T]$.  Then we have an estimate for the gradient at the center of the ball at later times:
\begin{align*}
&\sqrt{1+|Du(y_0,t)|^2}\\
&\phantom{space}\le C_1 \sup_{y\in B_R(y_0)} \sqrt{1+|Du(y,0)|^2} \exp\left[\frac {C_2}{R^2}\left(\sup_{B_R(y_0)\times[0,T]} u(y,t)-u(y_0,t)\right)^2\right], 
\end{align*}
where $C_1$ and $C_2$ depend only on $n$.  
\end{theorem}

We will also use the following {\it{a priori}} estimates for higher derivatives, from the same paper:
\begin{theorem}[$C^2$ interior estimate for mean curvature flow] 
 \label{interior-c2-estimate}
Suppose that $u$ satisfies  \eqref{mcf for graphs} on $B_R(y_0)\times[0, T]$.  Then for arbitrary $0\le \theta<1$ the estimate
\begin{equation*}
\sup_{B_{\theta R}(y_0)} |A|^2(t)\le c(n)(1-\theta^2)^{-2}\left(\frac1 {R^2} +\frac1 t\right) \sup_{ B_R(y_0)\times [0,t]}(1+|Du|^2)^2
\end{equation*}
holds for all $0\le t\le T$.  Here $|A|^2=h_{ij}h_{kl}g^{ik}g^{kl}$ .
\end{theorem}

\begin{theorem}[$C^k$ interior estimate for mean curvature flow]
\label{interior-higher-estimate}
Suppose that $u$ satisfies  \eqref{mcf for graphs} on $B_R(y_0)\times[0, T]$. 
Then for $m\ge 0$ and
arbitrary $0\le \theta<1$ we have the estimate
\begin{equation*}
\sup_{B_{\theta R}(y_0)} |\nabla^mA|^2(t)\le c_m(1-\theta^2)^{-2}\left(\frac1 {R^2} +\frac1 t\right)^{m+1},
\end{equation*}
where $c_m$ is a constant depending on $n,m$ and $\sup_{ B_R(y_0)\times [0,t]}(1+|Du|^2)^{1/2}$.  
\end{theorem} 

A bound on $|A|$ gives a bound on $|u|_{C^2}$, since (using coordinates in which   $h$ is diagonal) 
\begin{align*}
|A|^2&=h_{ij}h_{kl}g^{ik}g^{kl} \\
&= h_{ii} g^{ik}g^{kl}h_{ll} \\
&\ge \frac1{(1+|Du|^2)^2}\sum (h_{ii})^2 \\
&\ge \frac1{(1+|Du|^2)^3}\sup_{ij} |D_{ij}u|^2,
\end{align*}
where we have used that the smallest eigenvalue of $g^{ij}$ is $(1+|Du|^2)^{-1}$.  So,
\begin{equation*}
|D_{ij}u|\le (1+|Du|^2)^{3/2} |A|,
\end{equation*}
and in a similar manner, bounds on derivatives of $|A|$ give bounds on higher derivatives of $u$.

These estimates may be used to show long-time existence results, such as the following from \cite{eh:mean}:
\begin{theorem}
If $u_0$ is a locally Lipschitz,  entire graph over $\bigR^n$, then there is smooth solution to \eqref{mcf for graphs} for all $t>0$. 
\end{theorem}

In the following pages, we will derive new existence results of this sort.
\newchapter{Gradient estimates for parabolic equations of curve shortening flow type
 in one space dimension}{Gradient estimates in one space dimension}
\label{onedimensionalchapter}

In this chapter we outline gradient estimates for a class of parabolic equations in one space dimension.  %

  This chapter takes inspiration from the work of Huisken in  \cite{gh:dist}, where he investigated embedded plane curves evolving by curve shortening flow by looking at the evolution equation for the quotient of  $d(p,q)$, the distance between two points $p$ and $q$ in the metric of the plane, and $l(p,q)$, the length of curve between $p$ and $q$.  This introduced a double set of space coordinates (those around the point $p$ and those around the point $q$). At a maximum point, the first and second derivative conditions give strong conditions at both $p$ and $q$, allowing close examination of all possible situations.  An application of the maximum principle
resulted in a new proof of Grayson's theorem regarding the evolution of embedded curves.

In this chapter, we follow the approach of     Kru{\v{z}}kov in \cite{kruzhkov:nonlinear} (and well described in Lieberman's book \cite{li:parabolic}, chapter XI, section 6). 

If  $u(x,t)$ solves a parabolic partial differential equation in one space variable, then $v(x,y,t)=u(x,t)-u(y,t)$ solves 
a parabolic equation in two space variables, for which we can seek a barrier.  

In the paper cited, Kru{\v{z}}kov was interested in fully nonlinear equations 
\begin{equation*}
u_t=F(u_{xx},u_x,u,x,t)
\end{equation*}
with uniform parabolicity condition
\begin{equation*}
\frac \partial {\partial r}F(r,p,q,x,t)\ge A>0.
\end{equation*}
In this section, we do not require uniform parabolicity, but in order to show existence of the barriers, we will require a 
scaling similar to that of the curve shortening flow equation 
\begin{equation*}
u_t=\frac{u_{xx}}{1+{u_x}^2}, \label{Curve Shortening Flow}
\end{equation*}
in that $\pd F r \sim |p|^{-2}$ for large $|p|$.

We begin with a description of the ideas motivating the method.
The notation $'$ will indicate derivatives with respect to the space variable, which I hope I will use only  where this is unambiguous.

\section{Outline of the `double coordinate' method}
\label{general discussion section of 1d chapter}

Consider a smooth  $u:\mathbb R \times [0,T)\rightarrow \mathbb R$ satisfying 
\begin{gather*} u_t=a(u_x,u,x,t)u_{xx} +b(u_x) \\
u(x,0)=u_0. 
\end{gather*}

Let $Z:\mathbb R \times \mathbb R\times (0,T)\rightarrow \mathbb R $ be given by  $$Z(x,y,t)={u(y,t)-u(x,t)-{\phi(|y-x|,t)}},$$ where $\phi$ is some smooth function.

Suppose now that $Z$ attains a maximum at some point $(x,y,t)$, with $y> x$.  At the maximum point, the first derivatives  are zero, and so
\begin{equation} \label{one dimensional chapter, first derivatives}
\begin{split}
0&=%
Z_x=-u'(x,t)+\phi'(y-x,t) \\
0&=%
Z_y=u'(y,t)-\phi'(y-x,t). 
\end{split}
\end{equation}
Similarly, the matrix of second order partial derivatives is non-positive, by which we mean that for all $v\in \bigR^2$, $v^T [D^2 Z] v\le 0$, where $[D^2 Z]$ is the Hessian matrix
\begin{equation} \label{one dimensional chapter, second derivatives}
\begin{bmatrix} Z_{xx} & Z_{xy} \\
                        Z_{yx} & Z_{yy} \end{bmatrix} 
 = \begin{bmatrix} -u''(x,t)-\phi''(y-x,t) & \phi''(y-x,t) \\
        \phi''(y-x,t) &  u''(y,t)-\phi''(y-x,t) \end{bmatrix}.
\end{equation}

 If we now consider the evolution equation satisfied by $Z$,
\begin{align*}
\frac{\partial Z}{\partial t}&= u_t (y,t)- u_t (x,t)- \phi_t (|y-x|,t) \\
&= a\left(u_y,u(y,t),y,t\right)u_{yy}(y,t)-a\left(u_x,u(x,t),x,t\right)u_{xx}(x,t) \\
& \phantom{spacespacespacespace} 
{+b(u_y) -b(u_x) - \phi_t (|y-x|,t)};
\end{align*}
{and if we take this at the local maximum we have} \begin{align*}
\frac{\partial Z}{\partial t}&= a\left(\phi',u(y,t),y,t\right)\left(Z_{yy}+\phi''\right)-a\left(\phi',u(x,t),x,t\right)\left(-Z_{xx}-\phi''\right)   \\
& \phantom{space} +b(\phi')-b(\phi') -  \phi_t (|y-x|,t)\\
&= {\trace} \left( 
\begin{bmatrix} a\left(\phi',u(x,t),x,t\right) & c_1 \\ c_2 & a\left(\phi',u(y,t),y,t\right) \end{bmatrix}
        \begin{bmatrix} Z_{xx} & Z_{xy} \\ Z_{yx} & Z_{yy} \end{bmatrix}\right) \\
&\phantom{space} + a\left(\phi',u(y,t),y,t\right)\phi'' 
 + a\left(\phi',u(x,t),x,t\right)\phi'' 
-(c_1+c_2)\phi'' - \phi_t,
\end{align*}
for some $c_1$ and $c_2$.  If the first matrix above is positive semi-definite, then as $[D^2Z]$ is negative semi-definite, the trace above is non-positive and
\begin{equation*}
\frac{\partial Z}{\partial t}\le \phi''\left[a\left(\phi',u(x,t),x,t\right) +a\left(\phi',u(y,t),y,t\right)-c_1-c_2 \right]- \phi_t.
\end{equation*}
A useful choice for $c_1$ and $c_2$ that makes the first matrix positive semi-definite is $c_1=-a\left(\phi',u(x,t),x,t\right)$, $c_2=-a\left(\phi',u(y,t),y,t\right)$; then 
\begin{equation}
\frac{\partial Z}{\partial t}\le 2\big(a\left(\phi',u(x,t),x,t\right)+a\left(\phi',u(y,t),y,t\right)\big)\phi''- \phi_t.  \label{equation from discussion}
\end{equation}

The idea now is to choose $\phi$ in a way so that at the local maximum, $Z_t\le0$.   We begin by observing that for simple equations, a solution to a simplified version of the equation itself is acceptable for use as the barrier $\phi$.  

\remark We could simplify the method by choosing $c_1=c_2=0$, in which case the factor of $2$ is absent from \eqref{equation from discussion}.  The use of cross-derivatives will be important when we extend this method to higher dimensions.

\section{An estimate for periodic solutions}

\begin{theorem}
Suppose $u:\mathbb R\times [0,T)\rightarrow \mathbb R$ is a $H_2$, periodic and bounded solution of 
\begin{align} &u_t=a(u_x)u_{xx} +b(u_x) \label{simpleequation}\\
\intertext{ with initial condition}
&u(\cdot,0)=u_0, \notag
\end{align}
with  $u(x,t)=u(x+L,t)$ and $\osc u(\cdot,t) \le M$.

If $\varphi:\mathbb R\times (0,T)\rightarrow \mathbb R$ is a solution of 
 \begin{equation*} \varphi_t \ge a(\varphi_x)\varphi_{xx},  \notag \\
\end{equation*}
 with the initial and boundary conditions
\begin{align*}
&\varphi(x,t)\rightarrow 1 \text{  as } t\rightarrow 0\text{ for } x>0,  \\
&\varphi(0,t)=0 \text{ for }t>0, \\
&\varphi(x,t)\rightarrow 1\text{ as }x\rightarrow \infty\text{ for  }t>0, 
\end{align*}
then
\begin{equation*}|u(x,t)-u(y,t)|\le M \varphi\left(\frac{|y-x|}M,\frac{4t}{M^2}\right).
\label{first1dresult}
\end{equation*} 
\end{theorem}

\proof  Following on from the previous remarks, we set $$Z(x,y,t):=u(y,t)-u(x,t)-\phi(|y-x|,t)$$ and choose 
 $\phi(z,t)=M \varphi(z/ M,  {4t}/{M^2})$. 

As $t\rightarrow0$, $\phi\rightarrow M$ for all $z\not=0$, and so $Z(x,y,0)\le0$.  

As $u$ is periodic, $Z(x,y,t)=Z(x+L,y+L,t)$ and so $Z$ is periodic over strips 
$$\left\lbrace \,(x,y): 2nL\le y+x\le 2(n+1)L \,\right\rbrace. $$ 
\begin{floatingfigure}[v]{5.5cm}
\begin{center}
\epsfig{file=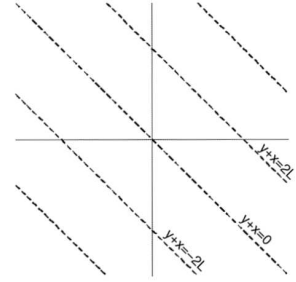, width=4cm}
\caption{$Z$ is periodic over strips} 
\end{center}
\end{floatingfigure}
Within each strip, $Z(y,y,t)=0$ and $$Z(x,y,t)=u(y,t)-u(x,t)-\phi\le M-\phi\rightarrow 0$$ as $|y-x|\rightarrow \infty$,  so $Z$ attains its spatial maximum in each strip, and hence in the entire plane, for each $t>0$.  
We can calculate
{\allowdisplaybreaks{\begin{align*}
\phi_t& =\frac 4 M \varphi_t\left(\frac z M, \frac {4t}{ M^2}\right) \notag \\ 
&\ge \frac 4 M a(\varphi')\varphi''  \notag \\
&= \frac 4 M a(\phi')M\phi'' \notag \\
&= 4a(\phi')\phi'' \label{ineqforphi},
\end{align*}}}
and in particular at a maximum point $(x,y,t)$ with $x\not=y$ and $Z$ non-negative,  equation \eqref{equation from discussion} becomes
\begin{equation*}
\pd Z t\le  4a(\phi')\phi''- \phi_t\le0.
\end{equation*}

Therefore, at such a maximum point $Z$ is non-increasing in time and so $Z\le0$.

The reason for the restriction $x\not=y$ is that $\phi$ is not differentiable here.  When the maximum {\it{is}} attained at such a point, then $Z(x,x,t)=0$, so in either case $Z\le0$ for all $t$.  The result %
follows.
\halmos

We can find explicit estimates for more general equations by choosing an explicit barrier.

\section{{Description of a barrier}}
\label{section where psi is defined}
This barrier, $\psi$, will be used often. 
\label{pagereference for definition of psi}  

Let $\Phi$ be the fundamental solution to the heat equation, $$\Phi(y,t)=\frac1{\sqrt t}\exp\left(\frac{-cy^2}t\right),$$ so that $\Phi_{yy}=4c\Phi_t$, where $c$ is a positive constant and $t>0$.    
Implicitly define $\psi(z,t)$ by 
\begin{equation*}z=\Phi(\psi-1,t)-\Phi(\psi+1,t). \label{defn for psi}\end{equation*}

This function has the property that as $t\rightarrow 0$, 
\begin{equation*}
\psi(z,t)\rightarrow \begin{cases} \phantom{-}1 & z>0, \\
                                -1 & z<0. 
                        \end{cases}  
\end{equation*} and that $\psi(0,t)=0$ for all $t>0$.    

We can calculate
{\allowdisplaybreaks[1]{\begin{gather*}   
\psi'
=\frac 1 {\Phi_y(\psi-1,t)-\Phi_y(\psi+1,t)}, \notag \\
\psi'' 
=- \left( \psi' \right)^3 \notag
\left[\Phi_{yy}(\psi-1,t)-\Phi_{yy}(\psi+1,t)\right],  
\intertext{and (third derivatives are included for completeness but not used until a later chapter)} 
\psi'''
=3\frac{(\psi'')^2}{\psi'}-(\psi')^4\left[ {\Phi_{yyy}(\psi-1,t)-\Phi_{yyy}(\psi+1,t)}\right], \label{third derivatives of psi} 
\intertext{while} 
\frac{\partial \psi}{\partial t}=-\frac{\Phi_t(\psi-1,t)-\Phi_t(\psi+1,t)}
{\Phi_y(\psi-1,t)-\Phi_y(\psi+1,t)}.  \notag
\end{gather*}}}

Routine calculations yield
{\allowdisplaybreaks \begin{gather}
\Phi_y=-\frac{2cy}{t^{3/2}}\exp\left(-\frac{cy^2}t\right), \notag \\
\Phi_{yy}=\frac{2c}{t^{3/2}}\left[\frac{2cy^2}t-1\right]\exp\left(-\frac{cy^2}t\right), \notag \\
\Phi_{yyy}=\frac{4c^2y}{t^{5/2}}\left[3-\frac{2cy^2}t\right]\exp\left(-\frac{cy^2}t\right),\label{third derivative of heat kernel p} \\
\Phi_t=\frac1{2t^{3/2}}\left[-1+\frac{2cy^2}{t}\right]\exp\left(-\frac{cy^2}t\right). \notag
\end{gather}
The partial differential equation satisfied by $\psi$ is 
\begin{equation*}
\frac{\partial \psi}{\partial t}= \frac1{4c}\frac{\psi''}{{\psi'} ^2}.
\end{equation*} }

\section{An explicit estimate for periodic equations}
\begin{theorem} \label{theorem for periodic gradient bounds} Let $u:\mathbb R \times [0,T)\rightarrow \mathbb R$ be a $H_2$ solution of
\begin{gather*} u_t=a(u_x,u,x,t)u_{xx} +b(u_x) \\
u(\cdot,0)=u_0 
\end{gather*}
where
$u_0$ is continuous; 
both $u_0$ and $a$ are periodic, \mbox{$u_0(x+L)=u_0(x)$},  $a(p,q,x,t)=a(p,q,x+L,t)$ (and therefore $u$ is also periodic); $\osc u(\cdot,t) \le M$; and where we can find positive constants $A$ and $P$ such that 
\begin{equation}
a(p,q,x,t)p^2\ge  A \text{ for all $|p|\ge P$.} \label{condition for a}
\end{equation}

Then there is a $T'>0$ such that for $t\in (0,T']$, $$|u_x|\le C_1 \sqrt t\left(1+t\right) \exp(C_2/t),$$
where $T'$, $C_1$ and $C_2$ are dependent on $M$, $A$ and $P$.

\end{theorem}

\proof 
 Let $Z$ be as before, with $\phi(z,t):=2M\psi(z/2M,t/4M^2)$ for the $\psi$ defined in Section \ref{pagereference for definition of psi}, with the constant $c$ given by $c=\max\left(\frac1{16A},CP^2\right)$, where $C$  will be chosen later.

Consider the region 
\begin{equation*} \label{definition of G}
G:=\left\lbrace\,(x,y,t): 0\le y-x \le z_M(t), 0\le t\le 2cM^2/3 \,\right\rbrace,
\end{equation*} where  $z_M(t)$ satisfies $\phi(z_M(t),t)=M$.  Explicitly,
\begin{equation}z_M(t)=\frac{4M^2}{\sqrt t}\left[\exp(-{cM^2}/t)-\exp(-{9cM^2}/t)\right].\label{defn of zm}\end{equation}   
\begin{center}
\begin{figure}[h]
\centering
\epsfig{file=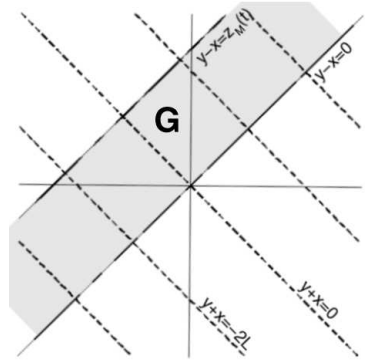, width=5cm}
\caption{The (periodic) region $G$ at some time $t$}
\end{figure}
\end{center}

As before, $Z$ is periodic over strips parallel to $y+x=0$ and so it attains its maximum on $G$.   We first show that $Z\le0$ on the boundary of $G$.

 For $y-x\not= 0$, as $t\rightarrow0$, $\phi\rightarrow 2M$ and so $Z<0$.

At $y-x=0$, $\phi(0,t)=0$ for all $t>0$, and so $Z(x,x,t)=0$.  At $y-x=z_M(t)$, $\phi(z_M,t)=M$ and so $Z\le 0$. 

Now suppose that $Z$ attains a maximum on the interior of $G$.  

It follows from \eqref{equation from discussion} that at the maximum,
\begin{align*}
\pd Z t &\le 
2\left[a\left(\phi',u(x,t),x,t\right)+a\left(\phi',u(y,t),y,t\right)\right]\phi'' 
-\phi_t \\
&=2\left[a\left(\psi',u(x,t),x,t\right)+a\left(\psi',u(y,t),y,t\right)\right]\frac{\psi''}{2M} -\frac{\psi_t}{2M} \\
&= \frac{\psi''}{2M} \left[2a\left(\psi',u(x,t),x,t\right)+2a\left(\psi',u(y,t),y,t\right)-\frac{1}{4c{\psi'}^2} \right] .
\end{align*}
The second derivative $$\psi''=- \left(\psi'\right)^3\left[\phi_{yy}(\psi-1,t/4M^2)-\phi_{yy}(\psi+1,t/4M^2)\right]  $$
is negative, as $\psi'$ is positive, and the part inside the square brackets $[\,\cdot\,]$ is positive in $G$ if $t\le2 cM^2/3=T'$.  

{\allowdisplaybreaks{ %
We can estimate  
\begin{align*}{\psi'}^2&\ge \inf_G {\psi '}^2 \\*
&\ge \left(\sup_{\substack{ 0\le\psi\le 1/2 \\  0\le t\le T'}}
                \Phi_y(\psi-1,t/4M^2)-\Phi_y(\psi+1,t/4M^2) \right)^{-2} \\
&\ge \left(\sup_{0\le\psi\le 1/2} \, \sup_{0\le t\le T'}
                \left|\Phi_y(\psi-1,t/4M^2)\right|+\left|\Phi_y(\psi+1,t/4M^2)\right| \right)^{-2} \\  
        &\ge  \left(\sup_{0\le\psi\le 1/2}      \left|\Phi_y\left(\psi-1,\frac23c(\psi-1)^2\right)\right|+\left|\Phi_y\left(\psi+1,\frac23c(\psi+1)^2\right)\right| \right)^{-2} \\ 
&= \left(\sup_{0\le\psi\le 1/2} \frac {2e^{-3/2} }{(2/3)^{3/2}\sqrt c |\psi-1|^2}
    + \frac{2e^{-3/2}} {(2/3)^{3/2}\sqrt c |\psi+1|^2} \right)^{-2} \\
&\ge c \frac{2 e^3}{3^3 5^2}  \\
&= {P}^2,
\end{align*} }}
where the last line follows by choosing  $C=3^35^2/(2e^3)$  and recalling that  $c\ge CP^2$.  

Now we can use the condition \eqref{condition for a}, controlling the degeneracy of $a$,  to estimate
\begin{align*}
\pd Z t &\le  \frac{\psi''}{2M} \left[2a\left(\psi',u(x,t),x,t\right)+2a\left(\psi',u(y,t),y,t\right)-\frac{1}{4c{\psi'}^2} \right] \\
&\le  \frac{\psi''}{2M} \left[ 4\frac{A}{{\psi'}^2} -\frac{1}{4c{\psi'}^2}  \right] \\
&\le 0
\end{align*}
as $c\ge (16A)^{-1}$.  
So,
$Z_t\le0$
at an internal maximum, $Z\le0$ on the boundary, and so %
$Z\le0$ on $G$.

Explicitly,  for $(x,y,t)\in G$, $Z\le 0$ means
\begin{align}
u(x,t)-u(y,t)&\le 2M\psi\left(\frac{|y-x|}{2M},\frac t{4M^2}\right) \notag \\
&\le  |y-x|\sup_z
\psi'\left(\frac z{2M},\frac t{4M^2}\right) \notag \\
&= |y-x| \psi'\left(0,\frac t{4M^2}\right) \notag  \\
&= |y-x| \frac{t^{3/2}}{4c(2M)^3}\exp\left(\frac {4cM^2}t\right),   \label{grad estimate one}
\end{align}
which is an estimate for the difference quotient $|u(x,t)-u(y,t)|/|x-y|$.

We can obtain an estimate for $(x,y,t)\not\in G$, $y>x$  by observing that for $z>z_M(t)$,
\begin{equation*}
u(y,t)-u(x,t)-\frac{Mz}{z_M(t)}\,\le\, u(y,t)-u(x,t)-M\,\le\, 0
\end{equation*}
so that as $|y-x|>z_M(t)$, 
\begin{equation}
u(y,t)-u(x,t)\,\le \,\frac {M|y-x|}{z_M(t)}\,\le\,  |y-x| \frac{\sqrt t}{2M} \exp\left(\frac{cM^2}t\right). 
\label{grad estimate two}
\end{equation}

So far, we have  estimates for when $y\ge x$.    We can find identical estimates for the region where $y<x$ by reflecting in the line $x-y=0$.  

Letting $y\rightarrow x$ in the estimates for the difference quotients \eqref{grad estimate one} and \eqref{grad estimate two} gives the result.
\halmos

Comparing $\varphi$ with  the special barrier $\psi$ gives us a gradient estimate for solutions of quasilinear equations    with this scaling.  We will use this estimate in 
 Chapter \ref{higher dimensional chapter}.

\begin{corollary}[Gradient estimates for the barrier $\varphi$] \label{particular varphi'(0) estimate}
Let $\varphi$ be a smooth solution of  
\begin{equation*}\varphi_t\le a(\varphi',\varphi,z,t)\varphi'', \end{equation*}
on $(0,\infty)\times(0,T)$, with the initial and boundary conditions 
\begin{align*}
&\varphi(z,0)= 1 \text{ for } t\ge0 \text{ and }z\not=0,  \\
&\varphi(0,t)=0 \text{ for }t>0, \\
&\varphi(x,t)\rightarrow 1\text{ as }x\rightarrow \infty\text{ for  }t>0.
\end{align*}
If
 \begin{equation*}
a(p,r,x,t)p^2\ge  A >0 \text{ for all $|p|\ge P>0$,}
\end{equation*}
then there is a $T'>0$ such that for $t\in(0,T']$
\begin{equation*}
\varphi'(0,t)\le C_1 \sqrt{t}(1+t)\exp{(C_2/t)}
\end{equation*}
where $T', C_1$ and $C_2$ depend on $A$ and $P$.
\end{corollary}

\proof 
We apply the comparison principle to $\varphi$ and $\psi^\epsilon(x,t):=\epsilon(1+t)+2\psi(x/2,t/4)$, where  $\psi$ is as in Section \ref{section where psi is defined} with the constant $c=(4A)^{-1}$.   

The barrier $\psi^\epsilon$ dominates $\varphi$ on $(0,\infty)\times\lbrace 0\rbrace$ and $\lbrace 0 \rbrace \times [0,T]$.  

Choose $T'>0$ and $z_1$ so that on $(0,z_1)\times(0,T')$, ${\psi^\epsilon}''\le 0$, ${\psi^\epsilon}'\ge P$ and at $z_1$, $\psi^\epsilon(z_1,t)> 1\ge \varphi(z_1,t)$.
 
Then $\boldsymbol{P}\varphi=-\varphi_t+a(\varphi',\varphi,z,t)\varphi''\ge 0$   and \begin{align*}
\boldsymbol{P}\psi^\epsilon&= -\epsilon -\psi_t + a(\psi',\psi^\epsilon,x,t)\psi'' \\
&=-\epsilon-\frac{\psi''}{(\psi')^2}\left[\frac1{4c}-a(\psi',\psi^\epsilon,x,t)\right] \\
&< 0,
\end{align*}  
so Theorem \ref{comparison principle} implies that $\psi^\epsilon>\varphi$ on $(0,z_1)\times(0,T')$ for all $\epsilon>0$, and as $\epsilon\rightarrow 0$, $\psi^0\ge\varphi$.    Since $\psi^0(0,t)=\varphi(0,t)$ 
this gives us the boundary gradient estimate %
\begin{equation*}
\varphi'(0,t)\le\psi'(0,t/4)\le  C_1  \sqrt{t}(1+t)\exp{(C_2/t)}.\end{equation*}
\halmos

\section{Interior estimates for non-periodic equations with $b=0$}

\begin{theorem} \label{interior 1d estimate}
 Let $u:\Omega\times[0,T]\rightarrow \bigR$ be a $H_2$ solution to 
\begin{equation*}
u_t=a(u_x,u,x,t)u_{xx},
\end{equation*}
where $\Omega\subset \bigR$ is an open interval 
and where there are positive constants $A$ and $P$ so that
\begin{equation}
a(p,q,x,t)p^2\ge A \quad {\text{ for all $|p|\ge P$. } } \label{csf-type condition}
\end{equation} 
If $\osc_{\Omega\times[0,T]} u \le M$, then we can find an estimate for $0<t<T'$ 
\begin{equation*}
|u'(x,t)|\le C_1 \sqrt{t}(1+t)  \exp{\left(\frac{C_2}t\right)},
\end{equation*}
where $T'$, $C_1$ and $C_2$ are dependent on $A$, $P$, $M$ and $\dist(x,\partial\Omega)$.
 
\end{theorem}

We  modify the previous proof, introducing a new boundary in the $x,y$ coordinates, since we no longer have compactness of the domain through periodicity.  We will seek to avoid a maximum of $Z$ occurring on the new boundary.

\proof  Firstly, suppose that $\Omega=(-1,1)$.  

We consider the sub-region 
$$ G:=\{\,(x,y,t)\in (-1,1)^2\times [0,T']: 0\le t\le T'\le  T, |y+x|<1, 0\le y-x< z_M(t)\,\}, $$
where  $z_M(t)$ is as before in \eqref{defn of zm} and $T'$ is chosen so that $z_M(t)< 1$ for $t\le T'$.

Define $Z$ on $G$ by 
\begin{equation*}
Z(x,y,t):=u(y,t)-u(x,t)-\phi(x,y,t).
\end{equation*}

In order to avoid a positive maximum of $Z$ on the boundary, we will ensure that $\phi$ satisfies 
\begin{itemize}
\item $\phi(x,y,0)\ge M $ for $y\not=x$ %
\item $\phi(x,y,t) \ge M$ for $|y-x|= z_M(t)$  %
\item $\phi(x,x,t)\ge 0$ for $t>0$  %
\item $\phi(x,y,t)\ge M$ for $|y+x|= 1$.  %
\end{itemize}

\begin{center}
\begin{figure}[h]
\centering
\epsfig{file=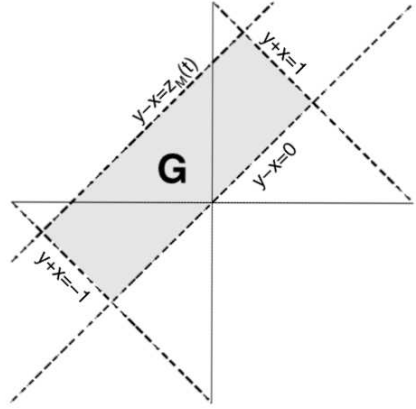, width=6cm}
\caption{The region $G$ at some time $t$}
\end{figure}
\end{center}

We choose 
\begin{equation*}\phi=2M\psi\left(\frac{y-x}{2M},\frac t{4M^2}\right) + \gamma(x+y-2\beta)^2, 
\end{equation*}
where $\psi$ is the explicit barrier defined in Section \ref{pagereference for definition of psi}, for positive constants $c$ and $\gamma$.    We will also choose $\beta$ with $|\beta|<1/2$ later.  

This $\phi$ satisfies the first three conditions.  In order to fulfill the final condition, choose $\gamma=M/(1-2|\beta|)^2$.  Then at $|x+y|=1$, 
$$\gamma(x+y-2\beta)^2\ge(\gamma\left(|x+y|-2|\beta|\right)^2=\gamma\left(1-2|\beta|\right)^2=M. $$ 

Now suppose that $Z$ first reaches a positive   maximum at an internal point $(x,y,t)\in G$.  As usual, we calculate that at this point first derivatives are zero
\begin{equation*}
\begin{split}
0&=Z_x=-u'(x,t)+\psi'-2\gamma(x+y-2\beta), \\
0&=Z_y=u'(y,t)-\psi'-2\gamma(x+y-2\beta), 
\end{split} \end{equation*}
 and the matrix of second derivatives is negative semi-definite \begin{equation*}
 0  \ge \begin{bmatrix} Z_{xx} & Z_{xy} \\
                        Z_{yx} & Z_{yy} \end{bmatrix} 
 = \begin{bmatrix} -u''(x,t)-\frac1{2M}\psi'' - 2\gamma & \frac1{2M}\psi'' -2\gamma\\
        \frac1{2M}\psi''-2\gamma &  u''(y,t)-\frac1{2M}\psi'' -2\gamma \end{bmatrix}.
\end{equation*}
 
We use these in the evolution equation for $Z$
{\allowdisplaybreaks\begin{align*}
\pd Z t&=a(u_y,u(y),y,t)u_{yy}- a(u_x,u(x),x,t)u_{xx}- \phi_t \\
&= a(u_y,u(y),y,t)\left[Z_{yy}+\frac{\psi''}{2M} + 2\gamma\right]-  a\left(u_x,u(x),x,t\right)\left[-Z_{xx}-\frac{\psi''}{2M} - 2\gamma\right] 
\\*&\phantom{spacespacespace}
-\frac{\psi_t}{2M} 
-\left(a(u_y,u(y),y,t)+ a(u_x,u(x),x,t)\right)Z_{xy} \\*
&\phantom{spacespacespace} + (a(u_y,u(y),y,t)+ a(u_x,u(x),x,t))\left(\frac{\psi''}{2M}-2\gamma\right)\\
&=\trace\left(\begin{bmatrix} a(u_x,u(x),x,t)  & -a(u_x,u(x),x,t) \\
                -a(u_y,u(y),y,t)         & a(u_y,u(y),y,t) \end{bmatrix}  
\begin{bmatrix} Z_{xx} & Z_{xy} \\
                        Z_{yx} & Z_{yy} \end{bmatrix}\right) \\*
& \phantom{spacespacespace} +2a(u_y,u(y),y,t)\frac{\psi''}{2M} +2a(u_x,u(x),x,t)\frac{\psi''}{2M} -\frac{\psi_t}{2M} \\
&\le \left[2a(u_y,u(y),y,t)+2a(u_x,u(x),x,t)-\frac1{4c{\psi'}^2}\right]\frac{\psi''}{2M},
\end{align*}}
the last line applying at the internal maximum .  As before, $\psi''\le0$ for $t\le T'\le 2cM^2/3$.  The first derivative condition for $Z$ implies that $u_x=\psi'-2\gamma(x+y-2\beta)$ and $u_y=\psi'+2\gamma(x+y-2\beta)$ and so 
\begin{equation}
\max(|u_x|,|u_y|)\ge|\psi'|\ge P,\label{lower bound on psi}\end{equation}
 where the final inequality comes from choosing $c\ge CP^2$ as in the previous section.   

We can then exploit \eqref{csf-type condition}, the condition on $a$; 
\begin{align*}
2a\left(u_y,u(y),y,t\right)+2a\left(u_x,u(x),x,t\right)&\ge \frac {2A}{\max\left(|u_x|,|u_y|\right)^2}
\\&\ge \frac {2A}{(|\psi'|+2\gamma|x+y-2\beta|)^2}.
\end{align*}
This is greater than $1/4c{\psi'}^2$ whenever $(\sqrt{8Ac}-1)|\psi'|\ge4\gamma$, so we use \eqref{lower bound on psi}, the lower bound on $|\psi'|$, and choose $$c\ge\frac{(4\gamma+P)^2}{8AP^2} =\left(\frac{4M}{(1-2|\beta|)^2}+P\right)^2\frac1{8AP^2}.$$

Now $Z_t\le 0$ at interior maxima, and so the parabolic maximum principle  ensures that $Z\le0$ on $G$.

For points outside $G$, but in the rectangle $|y+x|<1$, $z_M(t)\le y-x < 1$, 
\begin{align*}
u(y,t)-u(x,t)\le M &\le \frac {M |y-x|}{z_M(t)} \\
&\le \frac{\sqrt{t}}{2M}\exp{\left(\frac {cM^2}t\right)}|y-x|.
\end{align*}
We can repeat both these estimates for a reflected region where $x>y$; putting them all together gives  
\begin{equation*}
|u(x,t)-u(y,t)|\le2M\psi\left(\frac{|y-x|}{2M},\frac t{4M^2}\right) 
+\frac{M(x-y)^2}{(1-2|\beta|)^2}+ \frac{\sqrt{t}}{2M}\exp{\left(\frac {cM^2}t\right)}|y-x|.
\end{equation*}
Now, for $|y|< 1/2$, set $\beta=y$ and let $x\rightarrow y$ to give a gradient estimate at $y$
\begin{align*}
|u'(y,t)|&\le  \psi'\left(0,\frac t{4M^2}\right)+\frac{\sqrt{t}}{2M}\exp{\left(\frac {cM^2}t\right)} \\ 
&\le \frac{t^{3/2}}{24cM^3}
\exp\left (\frac{4cM^2}t \right)+ \frac{\sqrt{t}}{2M}\exp\left( \frac{cM^2}t \right) .
\end{align*}

For the case of a general interval $\Omega=[x_1,x_2]$, we can rescale around a point $y\in\Omega$ by using scaled coordinates $\tilde x=2(x-y)/{\dist(y,\partial\Omega)}$.  We obtain the estimate 
\begin{equation*}
|u'(y,t)|\le \frac 2{\dist(y,\partial\Omega)}\left[\frac{t^{3/2}}{24cM^3}\exp{\left(\frac{4cM^2}t\right)}+ M\sqrt{t}\exp{\left(\frac{cM^2}t\right)}\right],
\end{equation*}
where $c$ depends on $A$, $P$, $M$, and $\dist(y,\partial\Omega)$.
\halmos

\section{A generalisation to fully nonlinear equations}

In this section we consider equations of the form
\begin{equation} \label{nonlinear equation in one dimensional chapter}
u_t=F(u_{xx},u_x,u,x,t),
\end{equation}
where $F:\bigR^3\times\Omega\times[0,T]\rightarrow \bigR$ is $C^1$. %

Let $u$ be a smooth solution to \eqref{nonlinear equation in one dimensional chapter}.   
As before, define
$$Z(x,y,t):={u(y,t)-u(x,t)-{\phi(|y-x|,t)}},$$ 
and suppose it first becomes non-negative  at some point $(x,y,t)$, with $y>x$.  First and second derivatives of $Z$ will satisfy \eqref{one dimensional chapter, first derivatives} and \eqref{one dimensional chapter, second derivatives}, but the evolution equation for $Z$ will be given by
\begin{align*}
\frac{\partial Z}{\partial t}&= u_t (y,t)- u_t (x,t)- \phi_t (|y-x|,t) \\
&= F(u_{yy},u_y,u(y),y,t)-F(u_{xx},u_x,u(x),x,t) - \phi_t (|y-x|,t) \\
&=\int^1_0 \pd F s (su_{yy}\smallplus(1\smallminus s)u_{xx},su_y\smallplus(1\smallminus s)u_x,su(y)\smallplus(1\smallminus s)u(x),sy\smallplus(1\smallminus s)x,t)ds  \\
&\phantom{spacespacespace}- \phi_t (|y-x|,t)  \\
&= [u_{yy}-u_{xx}]\int^1_0 \pd F r \left(su_{yy}+(1-s)u_{xx},\dots\right)%
 ds
\\ &\phantom{spacespacespace}
+ [u_{y}-u_{x}]\int^1_0 \pd F p \left(su_{yy}+(1-s)u_{xx},\dots\right)\mspace{1mu}ds 
\\ &\phantom{spacespacespace}
+ [u(y)-u(x)]\int^1_0 \pd F q \left(su_{yy}+(1-s)u_{xx},\dots\right)\mspace{1mu}ds
\\&\phantom{spacespacespace}
+ [y-x]\int^1_0 \pd F x \left(su_{yy}+(1-s)u_{xx},\dots\right)\mspace{1mu}ds  - \phi_t (|y-x|,t)
\\ 
&= \trace  \left( 
 \begin{bmatrix} a & c_1 \\ c_2 & a%
\end{bmatrix} \begin{bmatrix} Z_{xx} & Z_{xy} \\ Z_{yx} & Z_{yy} \end{bmatrix} \right)    
  + 2\phi'' \int^1_0 \pd F r \left(su_{yy}+(1-s)u_{xx},\dots\right)\mspace{1mu}ds  
\\ &\phantom{spacespacespace}
-(c_1+c_2)\phi'' - \phi_t  
+ \phi \int^1_0 \pd F q \left(su_{yy}+(1-s)u_{xx},\dots\right)\mspace{1mu}ds
\\&\phantom{spacespacespace} 
+ [y-x]\int^1_0 \pd F x \left(su_{yy}+(1-s)u_{xx},\dots\right)\mspace{1mu}ds,
\end{align*}  
where we have used that $u_y=u_x$ at a spatial maximum, and have abbreviated\begin{equation*}
 {\int^1_0 \pd F r \left(su_{yy}+(1-s)u_{xx},\dots\right)\mspace{1mu}ds}=a\end{equation*}
and where we have added and subtracted $(c_1+c_2)Z_{xy}$ for some $c_1$, $c_2$.  If we choose $$c_1=c_2= -a %
$$ the first matrix above is positive semi-definite.  
Since the matrix of second derivatives is negative semi-definite, we have
\begin{align*}
\frac{\partial Z}{\partial t}&\le 4 \phi'' \int^1_0 \pd F r \left(su_{yy}+(1-s)u_{xx},\dots\right)\mspace{1mu}ds -\phi_t  
\\&\phantom{spacespacespace} 
+ \phi \int^1_0 \pd F q \left(su_{yy}+(1-s)u_{xx},\dots\right)\mspace{1mu}ds
\\&\phantom{spacespacespace} 
+ [y-x]\int^1_0 \pd F x \left(su_{yy}+(1-s)u_{xx},\dots\right)\mspace{1mu}ds.
\end{align*}

If we make quite harsh restrictions on $F$, then we can use our explicit barrier to find an analogue of Theorem \ref{theorem for periodic gradient bounds} for periodic nonlinear equations.

\begin{theorem}[Nonlinear version of Theorem \ref{theorem for periodic gradient bounds}]
Let $u:\bigR\times(0,T)\rightarrow\bigR$ be a $C^2$ solution of 
\begin{gather*} u_t=F(u_{xx},u_x,u,t) \\
u(\cdot,0)=u_0 
\end{gather*}
where
$u_0$ is continuous 
and periodic, \mbox{$u_0(x+L)=u_0(x)$},  
(and therefore $u$ is also periodic); $\osc u(\cdot,t) \le M$; where we can find positive constants $A$ and $P$ such that 
\begin{equation*}
\pd F r (r,p,q,t)p^2\ge  A \text{ for all $|p|\ge P$}; \label{condition for F}
\end{equation*} 
and where $\pd F q \le 0$ %
.  

Then there is a $T'>0$ such that for $t\in (0,T']$, $$|u_x|\le C_1 \sqrt t\left(1+t\right) \exp(C_2/t),$$
where $T'$, $C_1$ and $C_2$ are dependent on $M$, $A$ and $P$.

\end{theorem}

\proof  The proof is the same as that of Theorem \ref{theorem for periodic gradient bounds}, including the choice $\phi(z,t)=2M\psi(z/2M,t/4M^2)$, but instead of using inequality \eqref{equation from discussion} for interior maximum points, we have 
\begin{align*}
\frac{\partial Z}{\partial t}&\le 4 \phi'' \int^1_0 \pd F r (\dots) ds
+ \phi \int^1_0 \pd F q (\dots) ds
-\phi_t  \\
&\le 4 \phi'' \frac A{{\phi'}^2} - \frac{\phi''}{4c{\phi'}^2} \\
&\le 0,
\end{align*}
where the omitted argument of the derivatives of $F$, denoted by $(\dots)$, is $(su_{yy}+(1-s)u_{xx},su_y+(1-s)u_x,su(y)+(1-s)u(x),sy+(1-s)x,t)$.
\halmos

 \newchapter{An existence result for a parabolic equation in one space dimension}{Existence for one space dimension }
\label{short-time chapter}

Although this is a standard result (see Theorem 12.25 of \cite{li:parabolic}), for completeness we sketch a   short time existence result in the one-dimensional case, where the spatial domain is $\Omega=(x_0,x_1)\subset{\bigR}$, and the initial and boundary data is continuous.

The  parabolic equation is
\begin{align}
u_t&=a(u_x,u,x,t)u_{xx} \label{yetanotherlabel},
\end{align} 
with initial and boundary data  prescribed by 
\begin{equation} \label{boundary and initial data}
u(x,t)=u_0(x,t) \text{ for $(x,t)\in{\cal P}(\Omega\times[0,T])$}.
\end{equation}

We require that $a>0$ is in $H_{\alpha}(\cal K)$ for all bounded ${\cal K}\subseteq {\mathbb R}\times {\mathbb R}\times \Omega\times [0,T]$ and some $\alpha\in(0,1)$.

This implies that for every such $\cal K$ we can find 
positive  $\lambda_{\cal K}$ and  $\Lambda_{\cal K}$ such that 
\begin{equation} \lambda_{\cal K}\le a(p,q,x,t)\le\Lambda_{\cal K}, \text{ when } (p,q,x,t)\in{\cal K}. \label{definition of lambdas}\end{equation}

When we can find bounds of this type that depend only on the gradient, we will write
\begin{equation} \lambda(K)\le a(p,q,x,t)\le\Lambda(K), \text{ for  } |p|\le K. \label{more general upper and lower bound for a}
\end{equation} \label{page ref for K}

Suppose also  that 
there are positive constants $A$ and $P$ such that 
\begin{equation}
a(p,q,x,t)p^2\ge A >0,{\text{ for }}|p|\ge P.  \label{degeneracy control}
\end{equation}

The first part of this chapter is a survey of the main steps needed to find the existence result for the Cauchy-Dirichlet problem with $H_{1+\beta}$ initial and boundary data.  We follow the treatment in Lieberman \cite{li:parabolic}.

These results mean that when we approximate continuous initial data by smooth initial data, a solution will exist for the approximate initial data.      In the later parts of the chapter, we  use the gradient estimate established in Chapter \ref{onedimensionalchapter} to find uniform gradient estimates for $t>0$.  This will gives us a solution for $t>0$; in order to show that this approaches the initial data as $t\rightarrow 0$, we will need some displacement estimates which limit the distance a function can travel in a given time.

\section{Existence of solutions with $H_{1+\beta}$ initial and boundary data}

\begin{theorem} \label{first main existence} Consider the Cauchy-Dirichlet problem given by \eqref{yetanotherlabel} and \eqref{boundary and initial data}, where  $0<\beta<1$.  

Suppose that $u_0$ is defined on the parabolic boundary ${\cal P}\left(\Omega\times[0,T]\right)$ and
$u_0\in H_{1+\beta}({\cal P})$. 
 Also,  suppose that either $u_0$ is time-independent, or else 
there are constants $A$ and $P$ such that \eqref{degeneracy control} is satisfied.   

Then there is a smooth solution 
$u\in C^{2+1}\left(\Omega\times(0,T)\right)\cap C\left(\overline\Omega\times[0,T]\right)$.

This solution has a gradient bound $|u|_{1+\delta,\delta/2}\le C$ where $C$ depends on $|u_0|_{1+\beta,\beta/2}$, 
$\beta$,  $\lambda_{\cal K}$,  $\Lambda_{\cal K}$ and $\diam \Omega$.  
\end{theorem}

The proof of this result follows a standard pattern for showing existence ---
a bound on $\sup|u|$; 
a bound on $\sup|Du|$; 
a H\"older gradient bound $|Du|_\alpha$; 
and then the application of a fixed point theorem. 
These steps are sketched by the following results.

We begin by 
 using the comparison principle to bound $|u|$. 
\begin{lemma}[A bound on $\sup |u|$]  \label{1d existence- bound for |u|}
If $u$ is a smooth solution of \eqref{yetanotherlabel}, \eqref{boundary and initial data} in $\Omega\times [0,T]$, then  
\begin{equation*}
\sup_{\Omega\times[0,T]}|u(x,t)|\le \sup |u_0|. 
\end{equation*}
\end{lemma}

{\smallskip {\noindent\textbf{Proof idea:}}\quad} 
Set 
 $k=\sup {u_0}^+$
and apply  
 the comparison principle (Theorem \ref{comparison principle}) to $u$ and $k$ on $E:=\{\,(x,t)\in\Omega\times[0,T]:u(x,t)>0\,\}$.  
Since 
 $k\ge u$ on $\partial E$, it follows that $k\ge u$ on all of $E$ and so on all of $\Omega\times[0,T]$. 

Similar steps can be followed to find that $\inf_{\Omega \times [0,T]}u\ge \inf {u_0}^-$, completing the result.
\halmos

We begin our gradient estimates with a boundary gradient estimate. 
 
\begin{lemma}[Boundary gradient estimate]
Let $0<\beta\le1.$  If $u$ is a smooth solution of \eqref{yetanotherlabel}, \eqref{boundary and initial data}  in $\Omega\times [0,T]$, and either $u_0$ is time-independent or else $a$ satisfies the condition \eqref{degeneracy control}, then
\begin{equation*}
\sup_{\substack{  
(x,t)\in\partial\Omega\times[0,T] \\ 
(y,s)\in\Omega\times[0,t] }}
                \frac{|u(x,t)-u(y,s)|}{|(x,t)-(y,s)|}\le L,
\end{equation*}
where $L$ depends only on $\osc u_0$
, $|u_0|_{{1+\beta,\beta/2}}$
and $\beta$.
\end{lemma}

In fact, we can relax the regularity requirements on the initial and boundary data and still find a continuity estimate on the boundary.

A {\emph{modulus of continuity}}   \label{page where we define mod of cont}
is a concave, continuous function $\omega:\bigR^+\rightarrow\bigR^+$, with $\omega(0)=0$.  This $\omega$ is a modulus of continuity for a function $g$ at $y$ if  \nomenclature[omega]{$\omega$ }{ a modulus of continuity\refpage}
\begin{equation*}
|g(x)-g(y)|\le \omega(|x-y|)
\end{equation*}
for all $x$ in the domain of $g$.  It is a modulus of continuity for $g$ if the above relationship also holds for all $y$ in the domain of $g$.  

A modulus of continuity can be defined for every continuous function on a closed   bounded set.

\begin{lemma}[Boundary continuity estimate] \label{boundary continuity lemma} Let $u$ 
satisfy \eqref{yetanotherlabel}, \eqref{boundary and initial data},
where  $u_0$ has modulus of continuity $\omega$, and suppose there are positive constants $\mu$ and $P$ so that
\begin{equation}\label{structure condition}|p|\Lambda(p,q,x,t)+1\le\mu 
a(p,q,x,t)p^2\end{equation}
whenever $|p|\ge P$.   

Then $u$ has a modulus of continuity {\it{on the boundary}} 
$$ 
{|u(x,t)-u(y,t)|}\le \omega^*({|x-y|})
$$
for $x\in\Omega$ and $y\in\partial\Omega$, where $\omega^*$ can be determined by $\omega$,  $\sup_{\cal P }|u_0|$, 
$\Omega$, and $a$.
\end{lemma}

Equipped with the boundary gradient estimate, we can now find a global gradient estimate.  
In this one-dimensional case, the global gradient estimate is the result of Kru{\v{z}}kov mentioned   in Chapter \ref{onedimensionalchapter}.

\begin{lemma}[Global gradient estimate] 
If $u$ is a smooth solution of \eqref{yetanotherlabel}, \eqref{boundary and initial data} in $\Omega\times[0,T]$ with an oscillation bound $\osc u=M$, and a Lipschitz estimate on the parabolic boundary 
$|u(x,t)-u(y,t)|\le L|x-y|$ for all $(x,t)\in{\cal P}(\Omega\times[0,T])$ and $y\in\Omega$, then
\begin{equation*}
\sup_{\Omega\times[0,T]}|u_x|\le 2L.
\end{equation*}

\end{lemma}

\begin{lemma}[Global H\"older gradient estimate]
Suppose that $u$ satisfies \eqref{yetanotherlabel}, \eqref{boundary and initial data} in $\Omega\times [0,T]$, %
where there are positive constants  $\lambda_{\cal K}$ and $\Lambda_{\cal K}$ such that whenever 
$(p,q,x,t)$ in the set  ${\cal K}:=\lbrace\, (p,q,x,t): |p|\le K, |q|\le M, x\in\Omega,t\in[0,T]\,\rbrace$,
$$\lambda_{\cal K}\le a(p,q,x,t) \le \Lambda_{\cal K}.$$

If $u\in C^{2+1}(\Omega\times[0,T])\cap C(\overline\Omega\times[0,T])$, %
set $M=\sup|u|$ and $K=\sup|Du|$.   
Then there are positive constants $\alpha$ and $C$ determined by $\beta$, $\lambda_{\cal K}$, $\Lambda_{\cal K}$ and $\diam \Omega$ such that 
\begin{equation*}
[Du]_\alpha\le C\left(\sup|u|+\sup|Du|+|u_0|_{1+\beta,\beta/2} %
\right).
\end{equation*}
\end{lemma}

Now that we have bounds for $|u|_{1+\alpha,\alpha/2}$, we can apply the following existence theorem, which is derived from a fixed point theorem.  %

\begin{lemma}[Existence theorem]
Let $u_0$ %
be in $H_{1+\delta}$ for some $\delta\in(0,1)$ %

If there is a constant $M_\delta$ independent of $\epsilon$ such that any solution of \eqref{yetanotherlabel}, \eqref{boundary and initial data} on $\Omega\times[0,\epsilon)$ satisfies 
\begin{equation*}
|u|_{1+\delta,\delta/2} \le M_\delta,
\end{equation*}
then there is a solution of the Cauchy-Dirichlet problem \eqref{yetanotherlabel}, \eqref{boundary and initial data} in $\Omega\times[0,T]$.  

\end{lemma}

\section{Displacement estimates}
\label{displacement estimates}
The following estimates for the displacement suffered  in a given interval of time by a function moving under a parabolic flow apply to any strictly parabolic operator satisfying bounds of the form \eqref{definition of lambdas} or \eqref{more general upper and lower bound for a}.

\begin{lemma}[Displacement estimate 
for Lipschitz initial data] \label{displacement}

Let $u:\Omega\times[0,T]\rightarrow \bigR$ satisfy \eqref{yetanotherlabel}, where $\Omega\subseteq\bigR$, and $a$ has bounds of the form \eqref{more general upper and lower bound for a}.

Suppose that $u$ has initial data whose graph lies below a cone centred at some point $h$
\begin{equation*}
u(x,0)\le L|x-h|
\end{equation*}
and, in the case that $\Omega\not=\bigR$, whose boundary data lies below the same cone
\begin{equation*}
u(x,t)\le L|x-h|, \quad x\in\partial\Omega.
\end{equation*}
Then, at later times,
\begin{equation}
u(x,t)\le L(x-h)\erf{\frac {x-h}{2\sqrt{\Lambda t}}}+2L\sqrt{\frac{\Lambda t}\pi}\exp\left({-\frac{(x-h)^2}{4 \Lambda t}}\right), \label{Lipschitz displacement}
\end{equation}
where $\Lambda=\Lambda(L)$ is given by \eqref{more general upper and lower bound for a}.
\end{lemma}

\proof   
For some small $\epsilon>0$, set 
\begin{equation*}
v(x,t):=%
L(x-h)\erf{\sqrt{\frac c {t+\epsilon}}(x-h)} + {L}\sqrt{\frac{t+\epsilon}{c\pi}}\exp\left(-\frac{c(x-h)^2}{t+\epsilon}\right),
\end{equation*}
which satisfies the heat equation
\begin{equation*}
v_t=\frac{v_{xx}}{4{c}}
\end{equation*}
 and approaches the cone of gradient $L$ centred at $h$ as $t+\epsilon\rightarrow 0$.
\begin{center}
\begin{figure}[h]
\centering
\epsfig{file=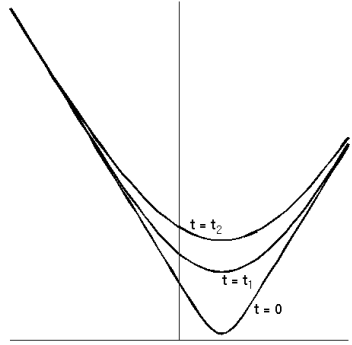, width=5cm}
\caption{The barrier $v(\cdot,t)$}
\end{figure}
\end{center}

Note that $v(x,0)>L|x-h|$, that $|v_x(x,t)|=L\left|\erf{x\sqrt{c/({t+\epsilon})}}\right|<L$ and that $v_{xx}>0$, so 
\begin{align*}
v_t- a(v_x,v,x,t)v_{xx} &\ge   v_t - \sup_{|p|\le L} a(p,v,x,t)v_{xx} \\
&\ge v_t - \Lambda(L) v_{xx} \\
&=0 
\end{align*}
where we choose $c^{-1}={4\Lambda (L)}$.

The estimate follows by applying the comparison principle 
(Theorem \ref{comparison 1}) to show that $u(x,t)\le v(x,t)$, and then letting 
$\epsilon\rightarrow 0$.

\halmos

Now, we apply this to three different cases, firstly when $u$ initially satisfies a H\"older condition and when we have polynomial growth in $\Lambda$,  secondly when $u(\cdot,0)$ has a modulus of continuity, and thirdly when $u$ is initially bounded by a step function.

\begin{corollary}[Displacement estimate for H\"older initial data] 
\label{holder-displacement}
Let $u:\bigR\times[0,T]\rightarrow \bigR$ satisfy \eqref{yetanotherlabel}.  Suppose that $a$ not only satisfies  
\eqref{more general upper and lower bound for a},  but more specifically has at most polynomial growth, so that  
\begin{equation}
a(p,q,x,t)\le \bar \Lambda (1+ K^m) \text{ for $|p|\le K$} \label{polynomial upper bound for a, again} 
\end{equation} where $\bar\Lambda$ and $m$ are positive constants.
Also, suppose that $u(\cdot,0)$ %
satisfies a H\"older condition around some point $h$
\begin{equation*}
|u(h,0)-u(x,0)|\le L|x-h|^\alpha, \quad 0<\alpha<1.
\end{equation*}
Then, at later times,
\begin{equation*}
|u(h,0)-u(h,t)|\le c(\alpha, m,L,{\bar\Lambda})t^{\frac \alpha{2+m(1-\alpha)}}.
\end{equation*}
\end{corollary}

\proof
\begin{figure}[h]
\begin{center}
\epsfig{file=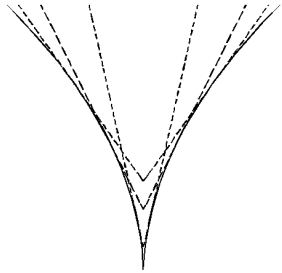, width=4cm}
\caption{The bounding cusp is itself bounded above by cones}
\end{center}
\end{figure}
For simplicity, assume $h=0$ and $u(h,0)=0$.  The initial data is bounded above by cones centred at $h$ and indexed by $k$, the (positive) $x$-coordinate of the point of contact with the bounding cusp $L|x|^\alpha$, so
\begin{equation*}
u(x,0)\le L|x|^\alpha\le \alpha L  k^{\alpha-1} |x| 
 +L (1-\alpha)k^\alpha. 
\end{equation*}

The estimate \eqref{Lipschitz displacement} taken at $x=0$ is then
\begin{align*}
u(0,t)%
&\le 2L\alpha k^{\alpha-1}
\left(1+\left(L k^{\alpha-1}\right)^m\right)^{1/2}
\sqrt{\frac{{\bar\Lambda} t}\pi}  
+L(1-\alpha)k^\alpha 
\intertext{ and optimizing over $k$ gives} 
u(0,t)& \le c(\alpha, m,L,{\bar\Lambda})t^{\frac \alpha{2+m(1-\alpha)}}.
\end{align*}
\halmos

\afterpage{\clearpage} %

\begin{corollary}[Displacement estimate for continuous initial data] 
\label{modulus-displacement}
Suppose that \linebreak $u: \bigR \times [0,T] \rightarrow \bigR$ satisfies \eqref{yetanotherlabel} and  
\eqref{more general upper and lower bound for a},
 where $u$ has initial data with a modulus of continuity $\omega$ at a point $h$
\begin{equation*}
|u(h,0)-u(x,0)|\le \omega\left(|x-h|\right).
\end{equation*}

Then
\begin{equation*}
|u(h,0)-u(h,t)|\le c(t)
\end{equation*}
where $c$ is dependent on $\omega$ and $\Lambda$, and where $c(t)\rightarrow 0$ as $t\rightarrow0$.
\end{corollary}

\proof
For simplicity, assume $h=0$ and $u(h,0)=0$.  Consider the cones 
$$C_k(x)=c_k\left(|x|-k\right)+ \omega(k)$$
indexed by $k>0$, the (positive) $x$-coordinate of a point of contact with $ \omega$.  As $ \omega$ is concave it has both left and right derivatives, and we can choose the slope of the cone $c_k= \omega'_-(k)$.  Then 
\begin{equation*}
u(x,0)=u(x,0)-u(0,0)\le  \omega(|x|)\le  \omega'_-(k)\left(|x|-k\right)+ \omega(k)=C_k(x).
\end{equation*}

Now we have a cone as an upper boundary, we can use estimate \eqref{Lipschitz displacement} at $x=0$
\begin{align*}
u(0,t)%
&\le  2 \omega'_-(k)\sqrt{\frac{\Lambda t}\pi}- \omega'_-(k)k+ \omega(k),
\end{align*}
where $\Lambda=\Lambda( \omega'_-(k))$ is given by \eqref{more general upper and lower bound for a}. 
Minimize this over $k$ to get the displacement bound 
$$c(t):=\inf_{k>0}\left(2 \omega'_-(k)\sqrt{\frac{\Lambda t}\pi}- \omega'_-(k)k+ \omega(k)\right).$$
 
In order to show that $c(t)\rightarrow 0$ as $t\rightarrow0$, let $\delta>0$.  As $ \omega$ is concave and positive, it has positive left derivative and for $k>0$ we have 
$$0\le  \omega'_-(k) k\le  \omega(k). $$
And as $ \omega$ is continuous,    
$$
0\le \lim_{k\rightarrow 0}  \omega'_-(k)k \le \lim_{k\rightarrow 0} \omega(k)=0.
$$
Choose $k=k_\delta$ so that $ \omega(k_\delta)- \omega'_-(k_\delta)k_\delta<\delta.$  Choose $\tau$ so that 
$$\sqrt\tau\le \frac{\delta\sqrt{{\Lambda\left( \omega'_-(k_\delta)\right)}}}{2 \omega'_-(k_\delta)\sqrt\pi},$$  
then for all $t\le\tau$, $$c(t%
)\le 2 \omega'_-(k_\delta)\sqrt{\frac{\Lambda t}\pi}- \omega'_-(k_\delta)k_\delta+ \omega(k_\delta)\le 2\delta,$$ and so $c(t)\rightarrow 0$.
\halmos

Set $\sigma$ to be the maximal monotone graph 
\begin{equation} 
\label{first instance of step function f}
\sigma(x)=\begin{cases} +1, &x>0  \\ [-1,1], & x=0 \\-1, & x<0\end{cases}
\end{equation}
which we will refer to as the {\emph{step ``function''}}.

\begin{corollary}[Displacement estimate for step functions]
\label{step displacement}
If $u$ satisfies \eqref{yetanotherlabel} and \eqref{more general upper and lower bound for a}, and is initially bounded by a step function  
\begin{equation*}
u(x,0)\le c\sigma(x),
\end{equation*}
then for $x<0$ 
\begin{equation*}
u(x,t)\le \min\left\lbrace \frac{4c}{|x|}\sqrt{\frac{\Lambda t}\pi}-c,c\right\rbrace,  %
\end{equation*}
where $\Lambda =\Lambda (2c/|x|)$ as in \eqref{more general upper and lower bound for a}.
\end{corollary}
\proof  
Near some point $h<0$, $u(\cdot,0)$ satisfies a Lipschitz condition 
\begin{equation*}
u(x,0)\le L^h|x-h|-c
\end{equation*}
where $L^h=2c/|h|$.  
\begin{center}
\begin{figure}[h]
\centering\epsfig{file=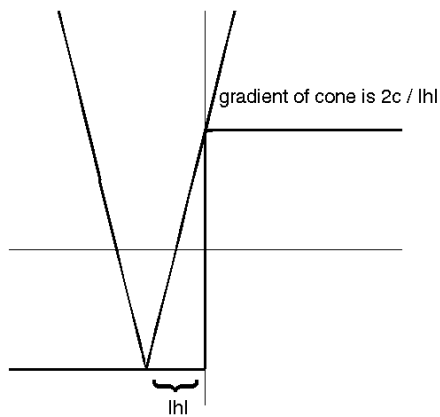, width=6cm}
\caption{Cone bounding the step function}
\end{figure}
\end{center}
Lemma \ref{displacement} then gives that 
\begin{equation*}
u(x,t)\le \inf_{h<0} \frac {2c}{|h|} \erf{\frac{x-h}{2\sqrt{\Lambda t}}}+\frac {4c}{|h|}\sqrt{\frac{\Lambda t}\pi}\exp\left(-\frac{(x-h)^2}{4\Lambda t}\right)-c
\end{equation*}
and if we let $h=x$ then we find that for $x<0$,  
\begin{equation*}
u(x,t)\le \frac{4c}{|x|}\sqrt{\frac{\Lambda t}\pi}-c.
\end{equation*}

The final result is found by comparison to the constant function $c$.
\halmos  

\remark If $a$ satisfies the condition \eqref{polynomial upper bound for a, again}, then 
\begin{equation*}
u(x,t)\le C(c,\bar \Lambda ) c^{1+m/2}\sqrt {t}|x|^{-1-m/2}-c. 
\end{equation*}
\halmos

\section{Existence of solutions with continuous initial data}

\begin{theorem} \label{second main existence} Consider the Cauchy-Dirichlet problem given by 
\eqref{yetanotherlabel} and \eqref{boundary and initial data}. 
If $u_0\in C\left({\cal P}(\Omega\times[0,T])\right)$ %
and 
if there are constants $A$ and $P$ such that \eqref{degeneracy control} holds,   
then \eqref{yetanotherlabel}, \eqref{boundary and initial data} has a solution $u\in C^{2+1}\left(\Omega\times(0,T)\right)\cap C\left(\overline\Omega\times[0,T]\right)$.
\end{theorem}

The first step in the proof of the above is to approximate $u_0$ %
by $u_0^\epsilon$  
in $C^\infty$,
so that \linebreak[3]
$\sup_{x\in\Omega}{|u_0^\epsilon-u_0|}<\epsilon$.%

\begin{lemma}[Existence of solutions with approximate boundary data] \quad For all \linebreak
$\epsilon>0$,
 there exist solutions $u^\epsilon:\Omega\times[0,T]\rightarrow \bigR$ to \eqref{yetanotherlabel} with boundary data $u_0^\epsilon$ %
These solutions are in $C^{2+1}\left(\Omega\times[0,T]\right)\cap C\left(\overline\Omega\times[0,T]\right)$.
\end{lemma}
\proof
As $u_0^\epsilon$ is %
in the H\"older space $H_{1+\beta}$, this is a consequence of Theorem \ref{first main existence}. 
\halmos

\begin{lemma}[Existence of uniform oscillation bound] 
For all $\epsilon>0$, $$\osc u^\epsilon \le 4\left(\sup|u_0| %
\right).$$  
\end{lemma}
\proof 
For any fixed $\epsilon$, set $k=\sup {u_0^\epsilon}^+$
and apply  %
 the comparison principle (Theorem \ref{comparison principle}) to $k$ and $u^\epsilon$ on $E=\{\,(x,t)\in\Omega\times[0,T]:u^\epsilon(x,t)>0\,\}$.  
Since %
 $k\ge u^\epsilon$ on $\partial E$, it follows that $k\ge u^\epsilon$ on all of $E$ and hence on all of $\Omega\times[0,T]$, and so 
$$\sup_{\Omega\times[0,T]}|u^\epsilon(x,t)|\le \sup|u_0^\epsilon| %
\le 2\left(\sup| {u_0}| %
\right),$$
where the last inequality will hold for small enough $\epsilon$.  

This leads to a uniform oscillation bound for $u^\epsilon$, which we denote by $M$ ---
\begin{equation*}
\osc u^\epsilon\le 4\sup |u_0| %
=:M. 
\end{equation*}
\halmos

Theorem \ref{interior 1d estimate} gives a uniform gradient bound on interior sets, up to some time $T'>0$.   %
For  $t_0\in(0,T'/2)$, 
\begin{equation*}
|{u^\epsilon}_x|_{\Omega'\times (t_0,T')} \le C_1 \sqrt{t_0}(1+t_0) \exp{\left(\frac{C_2}{t_0}\right)}=:L(t_0),
\end{equation*}
where $T'$, $C_1$ and $C_2$ are dependent on $A$, $P$, $M$ and $\dist(\Omega',\partial\Omega)$.

\begin{lemma}[Higher regularity on interior sets] 
\label{Higher regularity on interior sets}
On interior sets $\Omega'\times(2t_0,T')$ we can estimate  higher derivatives \begin{equation*}
|{u^\epsilon}|_{2+k+\alpha
}\le C%
\end{equation*}
where $C$ depends on $\dist(\Omega',\partial\Omega)$, $\diam(\Omega)$, $t_0$, $A$, $P$, $|a|_\alpha$ and $M$.
\end{lemma}

\proof  
Once we have an oscillation bound $M$ and a gradient bound $L(t_0)$, \eqref{definition of lambdas} implies uniform parabolicity.  A uniform H\"older gradient bound on interior sets results from Theorem 12.2 of \cite{li:parabolic}.  In particular, on interior sets and when $T'/2>t_0>0$,
\begin{equation*}
[{u^\epsilon}_x]_{\alpha;\Omega'\times(2t_0,T')}\le C
\min\left\lbrace \dist(\Omega',\partial\Omega),\sqrt{t_0}\right\rbrace^{-\alpha},
\end{equation*}
where  both $\alpha$ and $C$ depend on $\lambda_{\cal K}$ and $\Lambda_{\cal K}$, given by \eqref{definition of lambdas}, with 
$${\cal K}=\lbrace \, (p,q,x,t): |p|\le L(t_0), |q|\le M, x\in \Omega, \text{ and } t\in[0,T]\,\rbrace $$
and $C$ also depends on $L(t_0)+M$, and $\diam\Omega$.

Equipped with a H\"older gradient bound, we can treat the equation as a uniformly parabolic equation with H\"older continuous coefficients, and use standard results, such as Theorem \ref{interior1}, to find that $u^\epsilon$ is uniformly bounded in $H_{2+\alpha}\left(\Omega'\times (2t_0,T')\right)$.  

From here, it is possible to use the bootstrapping method to obtain interior estimates for all higher derivatives.

\halmos

\begin{corollary} %
On any interior set $\Omega'\times(t_0,T')$, there exists a subsequence converging to some $u$ that also solves the partial differential  equation \eqref{yetanotherlabel}.
\label{corollary asserting convergence}
\end{corollary}
 
In order for this $u$  to be a solution of the Cauchy-Dirichlet problem, we need to show that $u$ attains the initial and boundary data.  

\begin{lemma}[Convergence to initial data]
On any spatially interior set $\Omega'$,
\begin{equation*}
\sup_{x\in\Omega'}|u(x,t)-u_0(x)|\rightarrow 0 {\text{ as }}t\rightarrow0.
\end{equation*}
\end{lemma}
\proof
Let $x$ be any point in $\Omega'$.  Let $\omega$ be a modulus of continuity for $u_0$.

We can off-set $u^\epsilon$ by defining
$$w^\epsilon(y,t):=u^\epsilon(y,t)-u^\epsilon(x,0)+u_0(x),$$ so that
$w^\epsilon(x,0)=u_0(x)$. 
Let $u$ be the limit of a subsequence $u^\epsilon$, as in Corollary \ref{corollary asserting convergence}. 
\begin{align*}
|u(x,t)-u_0(x)|&=  \lim_{\epsilon\rightarrow 0}|u^\epsilon(x,t)-u_0(x)|\\
&=\lim_{\epsilon\rightarrow0}|w^\epsilon(x,t)+u^\epsilon(x,0)-2u_0(x)|\\
&\leq \lim_{\epsilon\rightarrow 0} \left(|w^\epsilon(x,t)-u_0(x)| +   
|u^\epsilon(x,0)-u_0(x)|\right) \\
&= \lim_{\epsilon\rightarrow 0} |w^\epsilon(x,t)-w^\epsilon(x,0)| +  \lim_{\epsilon\rightarrow 0} |u^\epsilon(x,0)-u_0(x)|.
\end{align*}
The second of these terms is zero.  To estimate the first term, note that the approximations $u^\epsilon(\cdot,0)$ satisfy the same the same continuity condition as $u_0$,  and therefore so does $w^\epsilon(\cdot,0)$, with
$|w^\epsilon(0,x)-w^\epsilon(0,y)|\leq \omega(|x-y|)$,
for all $x,y\in\Omega$.  Corollary \ref{modulus-displacement} then gives the estimate
\begin{equation*}
|w^\epsilon(x,t)-w^\epsilon(x,0)|\le c(t),
\end{equation*} 
where $c$ depends only on the exact forms of $\omega$ and $\Lambda_{\cal K}$ (given in \eqref{definition of lambdas}).  In particular, $c$ is independent of $x$ and $\epsilon$, and as $c(t)\rightarrow 0$ as $t\rightarrow0$, the result follows.   
\halmos 

More specific continuity-in-time estimates are given by the continuity of the initial data and the upper growth bound of $a$.  If, for example, the initial data is H\"older continuous $$|u_0(x)-u_0(y)|\le L|x-y|^\alpha$$ 
and $a$ has polynomial growth in the gradient term, satisfying \eqref{polynomial upper bound for a, again}
for constants  $\bar \Lambda$ and $m$, then Corollary \ref{holder-displacement} indicates that $c(t)=C(\alpha,m,L,\bar\Lambda)|t|^\frac\alpha{2+m(1-\alpha)}$.

\begin{lemma}[Convergence to boundary data]
We can continuously extend $u(\cdot,t)$, defined on the interior of $\Omega$ at time $t$, to $\overline\Omega$.  Moreover,  $u=u_0$ on the boundary.
\end{lemma}

\proof
We need to show that for $y\in\partial\Omega$, $\lim_{x\rightarrow y}u(x,t)=u_0(y,t)$.  

We note that our parabolic equation satisfies condition \eqref{structure condition}, since 
\begin{equation*}
|p|\Lambda(p,q,x)+1%
= |p||a(p,q,x)|+1 \le \frac 2 A a|p|^2
\end{equation*}
for $|p|\ge P$, using \eqref{degeneracy control}.  

Let $\omega$ be a modulus of continuity for $u_0$.  As each $u_0^\epsilon$ has at least the same modulus of continuity as $u_0$, Lemma \ref{boundary continuity lemma} gives us an estimate uniform in $t$ and $\epsilon$,  
\begin{equation*}
{|u^\epsilon(x,t)-u_0^\epsilon(y,t)|}\le \omega^*({|x-y|}).
\end{equation*}

Then for a point $y\in\partial\Omega$ and fixed $t$,
\begin{align*}
\sup_{B_r(y)\cap\Omega}|u(x,t)-u_0(y,t)|
&= \sup_{B_r(y)\cap\Omega}
|\lim_{\epsilon\rightarrow 0}u^\epsilon(x,t)-u_0^\epsilon(y,t)+u_0^\epsilon(y,t)-u_0(y,t)| \\
&\le \sup_{B_r(y)\cap\Omega}
 \omega^* \left(|x-y|\right) \\
&=\omega^*(r)
\end{align*}
so as $|x-y|\rightarrow 0$, $\omega^*(|x-y|)\rightarrow 0$ and $u(x,t)\rightarrow u_0(y,t)$  --- that is,   we can continuously extend $u$ to $u_0$ on $\partial\Omega$ for $t>0$.
\halmos

\section{Existence of entire solutions with stepped initial conditions}
 
\label{step function existence section}

Consider equation \eqref{yetanotherlabel}, under the conditions on $a$ given by \eqref{definition of lambdas} and \eqref{degeneracy control}.

\begin{lemma}  There exist entire solutions to this equation with the periodic, crenellated initial data
\begin{equation*}
g_R(x)= M \sigma\left(\sin(\pi x/R)\right),
\end{equation*}
where $\sigma$ is given by \eqref{first instance of step function f}.  
\end{lemma}
\begin{center}
\begin{figure*}[h]
\centering
\includegraphics[width=12cm]{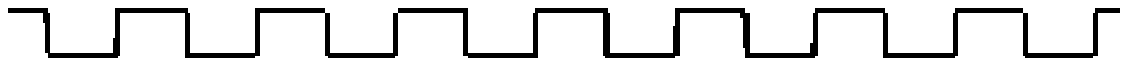}
\end{figure*}
\end{center}

\proof
If we let $g^\epsilon$ be the smooth mollification of $g_R$, then for $|x|<R$,  
$$M\sigma(x-\epsilon)\le g^\epsilon(x)\le M\sigma(x+\epsilon).$$
 
Theorem \ref{first main existence} ensures that there is a smooth solution $u^\epsilon$ to \ref{yetanotherlabel} with initial condition $g^\epsilon$, with a H\"older gradient bound dependent on $\epsilon$. The gradient bound in Theorem \ref{theorem for periodic gradient bounds} is independent of $\epsilon$; for $t\in(0,T')$,
\begin{equation*}
|u^\epsilon_x|\le C_1 \sqrt t\left(1+t\right) \exp(C_2/t)
\end{equation*}
where $T'$, $C_1$ and $C_2$ are dependent on $M$, $A$ and $P$, but not $R$.  Higher gradient bounds for $t>0$ follow from the interior estimate \eqref{Higher regularity on interior sets} and we can find a subsequence converging to $u_R$ which also solves the equation on $[t,T']$.

To show convergence of $u_R$ to the initial data, suppose that $-R/2<x<0$.  As in Section \ref{displacement estimates}, we can bound the initial data $g^\epsilon$ by cones centred at $h\in(-R/2,-\epsilon)$ ---
$$g^\epsilon(x)\le \frac M {|h+\epsilon|}|x-h|-M.$$
Applying Lemma \ref{displacement} to this, and setting $h=x$, we find that for $-R/2<x<-\epsilon$, 
\begin{equation*}
u^\epsilon(x,t)\le \frac{2M}{|x+\epsilon|}\sqrt{\frac{\Lambda t}\pi} -M,
\end{equation*}
where $\Lambda=\Lambda\left(M/{|x+\epsilon|}\right)$ as in \eqref{definition of lambdas} and so we have the estimate
\begin{align*}
|u_R(x,t)-g_R(x)|&=|\lim_{\epsilon\rightarrow 0}u^\epsilon(x,t)+M| \\
&= \frac{2M}{|x+\epsilon|}\sqrt{\frac{\Lambda t}\pi}. 
\end{align*}
A similar estimate holds for all $x\not= nR$, and so for all $\mu>0$, we can find $t$ (dependent on $x$) such that $|u_R(x,t)-g_R(x)|\le \mu$.     
\halmos

\begin{corollary}
\label{existence of entire step functions} 
There exists an entire solution to this problem with the initial data
\begin{equation*}
u_0(x)= M \sigma(x).
\end{equation*}
This solution has a gradient estimate for $t<T'$:
\begin{equation*}
|u_x|\le C_1 \sqrt t\left(1+t\right) \exp(C_2/t), 
\end{equation*}
where $T'$, $C_1$ and $C_2$ are dependent on $M$, $A$ and $P$.

\end{corollary}

\proof  Take the limit of the solutions $u_R$ given by the previous lemma as $R\rightarrow \infty$.  \halmos

  \newchapter{Gradient estimates for parabolic equations 
 in higher dimensions}{Gradient estimates in higher dimensions}

\label{higher dimensional chapter}

In this chapter we extend the methods of Chapter \ref{onedimensionalchapter} to higher dimensions.  %

Consider a smooth solution $u:{\mathbb R}^n\times[0,T]\rightarrow {\mathbb R}$ to
\begin{equation}
u_t=a^{ij}(Du,t)D_{ij}u+b(Du,t), \label{evolution equation in hd chapter}
\end{equation}
where $A(p,t)=[a^{ij}(p,t)]$ is a symmetric, positive semi-definite  $n\times n$ matrix that is smoothly dependent on $(p,t)\in{\mathbb R}^n\times[0,T]$. 

Define \begin{equation}
\alpha(p,t):=|p|^2\inf_{v\in S^n, v\cdot p\not=0} \frac{v^TA(p,t)v}{(v\cdot p)^2}.\label{defn of alpha}
\end{equation}

Compare this definition to that of the {\emph{Bernstein $\cal E$ function}}, (see Chapter 10 of \cite{trudinger})
\begin{equation*}
{\cal E} (p,q,x,t)=a^{ij}(p,q,x,t)p_ip_j.
\end{equation*}
Clearly, $\alpha(p)|p|^2\le {\cal E}$ and if $\lambda$, $\Lambda$ are the smallest and largest eigenvalues of $A$, then for $p\not=0$,
\begin{equation*}
\lambda\le \alpha(p)\le \frac{\cal E}{|p|^2}\le \Lambda.
\end{equation*}

The middle inequality here becomes an equality when $p$ is an eigenvector of $A$.

Our aim is to reduce the $n$-dimensional problem to a parabolic equation in one space dimension; we can do this if  $\alpha(p)$ is bounded below by a positive function of $|p|$.  We will call this 
\begin{equation}\tilde \alpha (s):= \inf_{p\in\bigR^n:|p|=s} \alpha(p).\label{alpha depends only on |p|}\end{equation}

For the existence of specific barriers we will require a control on the degeneracy of $A$ --- the existence
 of positive constants $A_0$ and $P$ such that 
\begin{equation} 
\tilde\alpha(s)s^2\ge A_0 \text{ for $s\ge P $}. \label{alpha-degeneracy}
\end{equation}

{\par\smallskip{\noindent\textbf{Example 1:}}\quad}  If $p$ is an eigenvector of $A(p)$, then $\alpha(p)$ is the associated eigenvalue. 

{\par\smallskip{\noindent\textbf{Example 2:}}\quad} As a specific example of the above, if $a^{ij}$ is of the form 
\begin{equation} \label{special form:  of mce-type?}
a^{ij}(p)=a_\infty(p)\left(\delta^{ij}-\frac{p_ip_j}{|p|^2}\right)+a_0(p)\frac{p_ip_j}{|p|^2}
\end{equation}
for functions $a_\infty,a_0:\bigR^n\rightarrow \bigR$,    with $a_0> 0$ and $a_\infty\ge0$, then $\alpha(p)=a_0(p)$.  

In the mean curvature flow case, $a_\infty=1$ and $$a_0(p)=\alpha(p)=\frac1{1+|p|^2}.$$ 

{\par\smallskip{\noindent\textbf{Example 3:}}\quad}
In the most general situation, if $v_k(p)$ are the non-null eigenvectors of $A(p)$ with eigenvalues $\lambda_k(p)>0$, then
\begin{equation*}
\alpha(p)=\begin{cases}
|p|^2\left(\sum_k \frac{(v_k\cdot p)^2}{\lambda_k}\right)^{-1} &  \text{if $p\in \left(\text{Null}A(p)\right)^\perp$} \\
0  &\text{otherwise.} 
\end{cases}
\end{equation*}

{\par\smallskip{\noindent\textbf{Example 4:}}\quad}
 In the case that $A$ is positive definite, all eigenvalues are positive  and $$\alpha(p)=\frac{|p|^2}{pA^{-1}p}.$$  As Example 2 shows, $A$ need not be positive definite.

{\par\smallskip{\noindent\textbf{Example 5:}}\quad} An elliptic operator $[a^{ij}]$ is called of \emph{mean curvature      type} if there are positive constants $\lambda$, $\Lambda$ so that 
\begin{equation*}
\lambda m^{ij}(p)\xi_i\xi_j\le a^{ij}(p,q,x)\xi_i\xi_j\le \Lambda m^{ij}(p)\xi_i\xi_j,\end{equation*}
where $m^{ij}$ are the coefficients    of 
mean curvature flow \cite{finn,simon,trudinger}.  For such equations, one can (under some conditions, particularly on the shape of the boundary)  find \emph{ apriori } estimates on $|Du|$ in terms of $|u|$.  It   is therefore interesting to note that if $[a^{ij}]$ satisfies only the lower inequality above, then it also satisfies \eqref{alpha-degeneracy}.

{\par\smallskip{\noindent\textbf{Example 6:}}\quad}  If one may be forgiven for referring to a future section, note that if the flow is the \emph{anisotropic mean curvature flow} \eqref{anisotropic curvature flows} of Section \ref{AMCF section}, then  Lemma \ref{lemma connecting amcf and alpha} implies that $\tilde \alpha(|p|)\ge A_0|p|^2$ when $\bar F(p)\ge P$.  The function $\bar F$ will be positive and homogeneous of   order one, so that $c_1|p|\le \bar F(p) \le c_2 |p|$.   

Thus, anisotropic mean curvature flow satisfies conditions \eqref{alpha depends only on |p|} and \eqref{alpha-degeneracy}.

\section{Reduction to a one-dimensional problem}

Let $u:\bigR^n\times(0,T]\rightarrow \bigR$ satisfy \eqref{evolution equation in hd chapter}, where $A=[a^{ij}]$ is a symmetric, positive semi-definite matrix with $\alpha>0$.
 
In the following, we generalise the calculations of Section \ref{general discussion section of 1d chapter} to higher dimensions.   As in the one-dimensional case, we begin our discussion by defining 
\begin{equation*}
Z(x,y,t):=u(y,t)-u(x,t)-\phi(|y-x|,t),
\end{equation*}
where $\phi:\bigR\times[0,T]\rightarrow \bigR$ is a $C^2$ function that will be chosen later.

At an internal maximum point of $Z$, the first derivative conditions are 
\begin{equation} \label{first derivatives in hd chapter}
\begin{split}
&D_{x^i}Z=-D_iu(x,t)-D_{x^i}\phi(|y-x|),t)=-D_iu(x,t)+\phi'\frac{y^i-x^i}{|y-x|}=0 \\
&D_{y^i}Z=\phantom{-}D_iu(y,t)-D_{y^i}\phi(|y-x|,t)=\phantom{-}D_iu(y,t)-\phi'\frac{y^i-x^i}{|y-x|}=0.\\
\end{split}
\end{equation}

The second derivatives of $Z$ are  {\allowdisplaybreaks{
\begin{align*}
D_{x^ix^j}Z&=-u_{ij}(x,t)-D_{x^ix^j}\phi(|y-x|,t) 
\\&= -u_{ij}(x,t)
        -\phi''\frac{(y^i-x^i)(y^j-x^j)}{|y-x|^2}
        -\frac{\phi'}{|y-x|}
        \left(\delta_{ij}-\frac{(y^i-x^i)(y^j-x^j)}{|y-x|^2}\right)\\
D_{y^iy^j}Z&=u_{ij}(y,t)-D_{y^iy^j}\phi(|y-x|,t)
\\&=
u_{ij}(y,t)
        -\phi''\frac{(y^i-x^i)(y^j-x^j)}{|y-x|^2}
        -\frac{\phi'}{|y-x|}
        \left(\delta_{ij}-\frac{(y^i-x^i)(y^j-x^j)}{|y-x|^2}\right)\\
D_{x^iy^j}Z&=-D_{x^iy^j}\phi(|y-x|,t)
\\&= \phi''\frac{(y^i-x^i)(y^j-x^j)}{|y-x|^2}
        +\frac{\phi'}{|y-x|}
        \left(\delta_{ij}-\frac{(y^i-x^i)(y^j-x^j)}{|y-x|^2}\right)
\end{align*}}}
and at a maximum point, the matrix $[D^2Z]$ must be negative semi-definite.

The evolution equation for $Z$ is 
{\allowdisplaybreaks[1]{\begin{align*} 
\frac{\partial Z}{\partial t}&= u_t(y,t)-u_t(x,t)-\phi_t \\* 
&= a^{ij}(D_yu,t)u_{ij}(y,t)+b(D_yu,t)-a^{ij}(D_xu,t)u_{ij}(x,t)-b(D_xu,t)
        -\phi_t\\*
& =  a^{ij}(D_yu,t) \left[D_{y^iy^j}Z + D_{y^iy^j}\phi(|y-x|,t)\right]  +b(D_yu,t)
\\* &\phantom{spacespace}
        -a^{ij}(D_xu,t)\left[-D_{x^ix^j}Z- D_{x^ix^j}\phi(|y-x|,t)\right] -b(D_xu,t) 
\\* & \phantom{spacespace} -\phi_t +2c^{ij}D_{x^iy^j}Z 
+2c^{ij}D_{x^iy^j}\phi(|y-x|,t) \\
& =  a^{ij}(D_yu,t) \left[D_{y^iy^j}Z +\phi''\frac{(y^i-x^i)(y^j-x^j)}{|y-x|^2}  %
+\frac{\phi'}{|y-x|} \left(\delta_{ij}-\frac{(y^i-x^i)(y^j-x^j)}{|y-x|^2}\right)\right]  \\*
& \phantom{spacespace}
        -a^{ij}(D_xu,t)\left[-D_{x^ix^j}Z
        -\phi''\frac{(y^i-x^i)(y^j-x^j)}{|y-x|^2}\right. 
\\*&\phantom{spacespacespacespacespace}
 \left. -\frac{\phi'}{|y-x|}\left(\delta_{ij}-\frac{(y^i-x^i)(y^j-x^j)}{|y-x|^2}\right)\right]  \\*
& \phantom{spacespace} +b(D_yu,t)-b(D_xu,t)-\phi_t +2c^{ij}D_{x^iy^j}Z \\*
&\phantom{spacespace}
+2c^{ij}\left[
        -\phi''\frac{(y^i-x^i)(y^j-x^j)}{|y-x|^2}
        -\frac{\phi'}{|y-x|}
        \left(\delta_{ij}-\frac{(y^i-x^i)(y^j-x^j)}{|y-x|^2}\right)\right]  ,
\end{align*}}}
where we add and subtract cross derivative terms with yet-to-be-chosen coefficients $c^{ij}$.  If we write $\xi=(y-x)/|y-x|$, and assume that we are at an internal maximum, then \eqref{first derivatives in hd chapter} implies that $D_iu(x,t)=\phi'\xi_i=D_iu(y,t)$, and we can continue the calculation:
\begin{align*} 
\pd Z t& =  a^{ij}(\phi'\xi,t)D_{y^iy^j}Z + a^{ij}(\phi'\xi,t)D_{x^ix^j}Z
        +2c^{ij}D_{x^iy^j}Z \\
&\phantom{=sp}
+\phi''\left(2a^{ij}(\phi'\xi)\xi_i\xi_j-2c^{ij}\xi_i\xi_j\right)
        +2\frac{\phi'}{|y-x|}\left(a^{ij}(\phi'\xi)-c^{ij}\right)
   \left(\delta_{ij}-\xi_i\xi_j \right) \\
&\phantom{spacespacespace} +b(\phi'\xi,t)-b(\phi'\xi,t)-\phi_t\\
&=\text{trace}
        \left(\begin{bmatrix}A(\phi'\xi) & C \\ C^T & A(\phi'\xi) \end{bmatrix} D^2Z \right)
 +2\phi''\left(\xi^T A \xi -\xi^T C \xi\right) \\
& \phantom{spacespacespace } 
+ 2\frac{\phi'}{|y-x|}\left(\text{trace}A - \xi^T A \xi -\text{trace}C 
        + \xi^T C \xi \right) -\phi_t  
\end{align*}
The idea now is to choose the off-diagonal block $C=[c^{ij}]$ in such a way that the $2n\times 2n$ matrix 
\begin{equation*}
A'= \begin{bmatrix} A & C \\ C^T & A \end{bmatrix}
\end{equation*}
is positive semi-definite,  leaves the coefficient of $\phi''$ positive and sets the coefficient of $\phi'/|y-x|$ to zero.

The first and third of these requirements imply that $C$ is given by  
$c^{ij}
= a^{ij}(\phi'\xi,t)-c\,\xi_i\xi_j$
for some $c>0$.   We can check that  
$$\xi^TA\xi-\xi^TC\xi=c>0$$
and that 
\begin{equation*}
\begin{split}
\trace C-\xi^T C\xi&= \trace[a^{ij}-c\,\xi_i\xi_j]-\xi_j(a^{ij}-c\,\xi_i\xi_j)\xi_i \\
&=\trace A-\xi^T A \xi.
\end{split}
\end{equation*} 

If we set $c=2\alpha(\phi'\xi,t)$, with $\alpha$ defined by \eqref{defn of alpha}, this maximizes the coefficient of $\phi''$, while keeping $A'$ positive semi-definite. 
For any $v,w\in {\mathbb R}^n$,
\begin{align*}
(v^T,w^T)A'\begin{pmatrix}v \\ w\end{pmatrix}&=v^TAv+w^TAw+2v^TCw \\
&=v^TAv+w^TAw+2v^TAw-4\alpha(\xi\cdot v)(\xi\cdot w)\\
&=(v+w)^TA(v+w) -\alpha\left[(\xi\cdot v+w)^2 - (\xi\cdot v-w)^2\right] \\
&\ge(v+w)^TA(v+w) -\alpha(\xi\cdot v+w)^2 \\ 
&\ge0.
\end{align*}

At the maximum point of $Z$, we find that
\begin{align*}
\frac{\partial Z}{\partial t}&\le 4\alpha(\phi'\xi)\phi'' -\phi_t. 
\end{align*}

In this way, we have reduced our problem to finding $\phi$ that satisfies the above equation, %
or, if we can find a lower bound on $\alpha$ dependent only on $|p|$, as in \eqref{alpha depends only on |p|}, then
\begin{equation*}4\tilde\alpha(\phi')\phi'' -\phi_t\le 0,\end{equation*}

We will use the results of Chapter 
\ref{onedimensionalchapter}
to do this.

\section{Estimates for periodic solutions}
\label{section periodic n>1}
The following theorem is an analogue of Theorem \ref{theorem for periodic gradient bounds}  for higher dimensions.    In the special case of mean curvature flow, this is joint work with Ben Andrews.  

Let $\alpha$ be defined by \eqref{defn of alpha}.  

\begin{theorem}
\label{periodic higherdimensional theorem}
Let $u:{\mathbb R}^n\times[0,T]\rightarrow{\mathbb R}$ be a smooth solution to 
\begin{gather*}
\frac {\partial u}{\partial t}=a^{ij}(Du,t)u_{ij}+b(Du,t) \\
u(x,0)=u_0(x)  
\end{gather*}
where $u_0$ is smooth with oscillation bound $\osc u_0\le M$,  and  is also spatially periodic, $u_0(x)=u_0(x+\Gamma)$, for some lattice $\Gamma$.

Suppose that $\alpha(p)=\tilde \alpha(|p|)>0$ for all $p$.   %

If $\varphi:{\mathbb R}^+\times [0,T]$ is a smooth solution to the auxilliary one-dimensional equation
\begin{equation} \label{conditions and equation for varphi}
\varphi_t=4\tilde\alpha(|\varphi'|,t)\varphi'', \end{equation}
and satisfies the boundary conditions 
\begin{equation} \label{conditions for varphi}
\begin{split}
&\lim_{t\rightarrow0} \varphi(z,t)=1 {\text{ for }} z\not=0  \\
&\varphi(0,t)=0 {\text{ for all }} t>0\\
&\lim_{z\rightarrow \infty} \varphi(z,t)\rightarrow 1 \text{ for all }t>0\\
\end{split}\end{equation}
then %
\begin{equation*}
|u(y,t)-u(x,t)|\le M\varphi\left(\frac{|y-x|}M,\frac t {M^2}\right).
\end{equation*}

\end{theorem}
\begin{corollary} \label{periodic higherdimensional corollary}
If there are positive constants $A_0$ and $P$ so that 
\begin{equation} 
\tilde \alpha(|p|)|p|^2\ge A_0 \text{ for $|p|\ge P $} \label{degeneracy for alpha}, 
\end{equation}
 then there is a $T'>0$ such that for $t\in(0,T'],$ 
\begin{equation*}
|Du|\le C_1\sqrt{t}(1+t)\exp(C_2/t),
\end{equation*}
where $T'$, $C_1$ and $C_2$ depend on $n$, $M$, $A_0$, and $P$.
\end{corollary}

\noindent{\textbf{Proof of Theorem \ref{periodic higherdimensional theorem}.} } The proof of this is substantially the same as the proof of Theorem \ref{theorem for periodic gradient bounds}, the gradient estimate for periodic, one-dimensional equations.  

As in the previous pages, let $$Z(x,y,t):=u(y,t)-u(x,t)-\phi(|y-x|,t),$$ and choose $\phi(z,t)=M\varphi\left(z/M,t/{M^2}\right)$, so that $Z(x,y,0)\le0$.   

As $u$ is periodic over the lattice $\Gamma=(L_1,\dots,L_n)$, $Z$ is periodic over regions
$$\lbrace \,(x,y,t)\in\bigR^{2n}\times[0,T]:2nL_i-x_i\le y_i\le2(n+1)L_i-x_i\,\rbrace. $$
On any one of these regions, note that $Z(y,y,t)=0$ and that $Z(x,y,t)\le M-\phi\rightarrow 0$ as $|y_i-x_i|\rightarrow\infty$, so $Z$ attains a spatial maximum on the region (and hence on the entire domain $\bigR^{2n}$).  

If there is a maximum point $(x,y)$ at some $t_0\in(0,T')$ with $x\not=y$, then at this point $Z$ is smooth and  
\begin{equation*}
\frac{\partial Z}{\partial t}\le 4\tilde\alpha(\phi'\xi)\phi'' -\phi_t\le
 4\tilde\alpha(|\varphi'|)\frac{\varphi''}M -\frac{\varphi_t}M=0.
\end{equation*} 
If $x=y$ at the maximum point, then here $Z(x,x,t)=0$ and in either case, $Z\le0$. %
The estimate for $|u(y,t)-u(x,t)|$ follows.
\halmos

\noindent{\textbf{Proof of Corollary \ref{periodic higherdimensional corollary}:}}
If $\tilde\alpha$ satisfies the degeneracy condition \eqref{degeneracy for alpha}, then for small times the gradient of $\varphi$ may be estimated by Corollary \ref{particular varphi'(0) estimate}.  Letting $x\rightarrow y$ gives that  
\begin{align*}
|Du(y,t)|&\le n \varphi'(0,t)\\
&\le   C_1\sqrt{t}(1+t)\exp(C_2/t).
\end{align*}
\halmos

\begin{draftcomment}

\end{draftcomment}

\section{Estimates for boundary value problems}
In the special case of mean curvature flow the following theorem is joint work with Ben Andrews.  

\begin{theorem}[Neumann problem] \label{Neumann theorem}
Let $\Omega\subset{\mathbb R}^n$ be a convex domain with $C^2$ boundary, and let $u$ be a smooth solution of\begin{gather*}
\frac {\partial u}{\partial t}=a^{ij}(Du)u_{ij}+b(Du,t) \\
D_\nu u(x,t)=0, {\text{ for }}x\in\partial\Omega, t>0.
\end{gather*} %
 If $a^{ij}$ and $\varphi$ satisfy the same conditions as in Theorem \ref{periodic higherdimensional theorem},
 then for any $x$ and $y$ in $\overline\Omega$, 
\begin{equation*} 
|u(y,t)-u(x,t)|\le \varphi(|y-x|,t),
\end{equation*} 
where $\osc u_0=M$.

Furthermore if $a^{ij}$ satisfies the degeneracy condition \eqref{alpha-degeneracy}, then for $t\in(0,T')$ a short-time gradient bound holds:
\begin{equation}
|Du(x,t)|\le C_1\sqrt{t}(1+t)\exp(C_2/t) {\text{ for $(x,t)\in\Omega\times(0,T']$}} 
\end{equation}
where $T'$, $C_1$ and $C_2$ depend on $n$, $M$, $A_0$, and $P$.
\end{theorem}

\proof As in the previous proof, set
$$Z:=u(y,t)-u(x,t)-M\varphi\left(|y-x|/M,t/M^2\right).$$   Note that $Z\le 0$
when $t=0$.

 For any $t>0$, suppose that $(x,y)$ is a spatial maximum of $Z$.  We will consider the possibility that $x$ and $y$ are both interior points, that $y$ is a boundary point and so is $x$, or that $y$ is a boundary point while $x$ is not (the converse follows without loss of generality). 

If both $x$ and $y$ are interior points, then the arguments of Theorem
\ref{periodic higherdimensional theorem} apply and $Z\le 0$ at this point.  

Consider the case that $y$ is on the boundary $\partial \Omega$.  If we take derivatives at $y$ that are in directions $\mu_y$ that have no component normal to the boundary,  then as before $D_{\mu_y} Z(x,y,t) =0$.  On the other hand,  let $\nu_y$ be the outward
unit normal at $y$.  The outwards-pointing derivative of $Z$ here is
\begin{align*}
\left.\dfrac d{ds} Z(x,y+s\nu_y,t)\right|_{s=0}&=\nu_y\cdot Du(y,t)-\varphi'
\frac{y-x}{|y-x|}\cdot\nu_y \\
&= 0-\varphi'\frac{y-x}{|y-x|}\cdot\nu_y \\
&\le 0,
\end{align*}
{\nopagebreak[3]{where we have used the boundary condition $D_{\nu_y} u(y,t)=0$ and that as $\Omega$
is convex,  $(y-x)\cdot\nu\ge 0$.}}

This inequality cannot be strict, for if it is, then there is a small
$s>0$ such that
\begin{equation*}
Z(x,y-s\nu_y,t)>Z(x,y,t)
\end{equation*}
which would contradict that $(x,y)$ is a maximum of $Z$.   Therefore $$\left.\dfrac d{ds} Z(x,y+s\nu_y,t)\right|_{s=0}=0,$$
and indeed $D_yZ(x,y,t)=0$. 

Now consider the position of  $x$.  If it is on the boundary, let $\nu_x$ be the
outward unit normal at $x$, and so
\begin{equation*}
\left.\dfrac d{ds} Z(x+s\nu_x,y,t)\right|_{s=0}\,=\,-\nu_x\cdot Du(x,t)+\varphi'\frac{y-x}{|y-x|}\cdot\nu_x %
\,\le\, 0,
\end{equation*}
Again, this inequality cannot be strict if $(x,y)$ is to be a maximum of
$Z$, so the outward derivative $\left.\dfrac d{ds} Z(x+s\nu_x,y,t)\right|_{s=0}=0$.  As before, other non-normal derivatives are also zero, so $D_x Z(x,y,t)=0$.

So, when both $x$ and $y$ are boundary points, $DZ=0$ and 
$[D^2 Z]$ is negative semi-definite.  We can argue as before that $Z_t\le 0$.

In the case that $x$ is an interior point, $D_{x} Z=0$ and so
$DZ=0$, $[D^2 Z]$ is negative semi-definite here,  and $Z_t\le 0$.

It follows that $Z\le 0$ for all $t>0$.  \halmos

The highly geometric nature of mean curvature flow allows us to relax the conditions on the convexity of the boundary.  In the following theorem we consider domains that are merely mean-convex.  This means that at every point on the boundary, the sum of the principal curvatures of $\partial\Omega$ is positive:
\begin{equation*}
\sum_{i=1}^{n-1}\kappa_i\ge 0.
\end{equation*}

Under the assumption of convexity (rather than mean-convexity), the following theorem is joint work with Ben Andrews.  

\begin{theorem}[Dirichlet problem for mean curvature flow] \label{mcf-dirichlet-theorem}
Let $\Omega\subset{\mathbb R}^n$ be a mean-convex domain with a $C^2$ boundary, and let $u$ be a smooth solution of the mean curvature flow for graphs
{\allowdisplaybreaks{\begin{gather*}  
\frac{\partial u}{\partial t}=
\left(\delta_{ij}-\frac{D_iu D_ju}{ 1+|Du|^2}\right)D_{ij}u, \\
\intertext{with prescribed boundary values}
u(x,t)=0\quad\text{for $x\in\partial\Omega$, $t>0$}, %
\\  
u(\cdot,0)=u_0. %
\end{gather*}}}

If $\varphi$ is a smooth solution to curve shortening flow
\begin{equation*}
\varphi_t=\frac{\varphi''}{1+(\varphi')^2},
\end{equation*}
with boundary conditions given by  \eqref{conditions for varphi},
then there is an estimate
\begin{equation*}
|u(y,t)-u(x,t)|\le 2M \varphi\left(\frac{|y-x|}{2M},\frac t{4M^2}\right),
\end{equation*}
where $M=\sup |u_0|$.
\end{theorem}
\proof
We find a boundary gradient estimate by defining a new $Z_B$ on $\Omega\times(0,T)$ which incorporates the distance to the boundary 
$$Z_B(y,t):=u(y,t)-2M\varphi\left(\frac{d}{2M},\frac t{4M^2}\right),$$
where 
$d(y)=\dist(y,\partial\Omega)$ is a $C^2$ function in the neighbourhood of the boundary      
$$\Omega\setminus\Omega^R:=\lbrace\, y\in\Omega :  \dist(y,\partial\Omega)\le R\,\rbrace.$$
Here, $R=(\sup_{\partial\Omega} \kappa_i)^{-1}$ and $\kappa_i$ are the principal curvatures  of $\partial\Omega$.  This close to the boundary,  each point $y$ has a unique closest point $x\in\partial\Omega$. 

Choose $T'>0$ so that for  $0<t<T'$,  if $0<d<R$ then $\varphi'\left(d/(2M),t/(4M^2)\right)\ge 0$, and if $d\ge R$ then $\varphi\left(d/(2M),t/(4M^2)\right)\ge 1$.   

At $t=0$, $Z_B\le 0$ for all $y\in\overline\Omega$.  For $y$ on the boundary, $Z_B(y,t)=0-2M\varphi\le 0$.  
 For  $t<T'$ and points at least distance $R$ from the boundary, $y\in\Omega^R$, $Z_B(y,t)\le u(y,t)- M\le 0.$

 We can find spatial derivatives for $Z_B$:  
{\allowdisplaybreaks{\begin{gather*}
D_{y^i}Z_B=D_{y^i} u-\varphi' D_i d  \\
D_{y^iy^j}Z_B=D_{y^iy^j} u - \frac{\varphi''}{2M}D_idD_jd-\varphi' D_{ij} d.
\end{gather*}}}

Now suppose that $y$ is an interior maximum of $Z_B$ at some time $t<T'$.
At this point, $DZ_B=0$ and $[D^2Z_B]$ is negative semi-definite, so  
\begin{align*}
\pd {Z_B} t &= u_t- \frac{\varphi_t}{2M} \\
&=m^{ij}(\varphi'Dd)\left[D_{y^iy^j}Z_B + \frac{\varphi''}{2M}D_idD_jd+\varphi' D_{ij} d\right] -\frac{\varphi_t}{2M} \\
&=m^{ij} D_{y^iy^j} Z_B + \frac{1}{2M}\frac{\varphi''}{1+\varphi'^2} + \varphi' D_{ii} d -\frac{\varphi_t}{2M},
\end{align*}
where we have used that $|Dd|=1$ and $D_idD_{ij}d=0$.

As in Lemma 14.17 of \cite{trudinger},
\begin{equation*} 
\trace[D^2d(y)]=\sum_{i=1}^{n-1} \frac{-\kappa_i}{1-\kappa_i d},
\end{equation*}
where $\kappa_i$ are the principal curvatures of $\partial\Omega$ at $x$, the closest point on $\partial\Omega$ to $y$.  If $d<R$, then $\kappa_i d<1$ and
\begin{equation*} 
  \sum_{i=1}^{n-1} \frac{-\kappa_i}{1-\kappa_i d}\,\le\, -\negthickspace\negthinspace\sum_{i=1}^{n-1} \kappa_i\,\le \,0,
\end{equation*}
the last inequality resulting from the mean-convexity of $\partial\Omega$.  Then \begin{equation*}
\pd {Z_B} t =m^{ij} D_{y^iy^j} Z_B + \frac{1}{2M}\frac{\varphi''}{1+\varphi'^2} + \varphi' D_{ii} d -\frac{\varphi_t}{2M}\le 0. \end{equation*}

It follows that $Z_B\le 0$ for $t<T'$, and so for all $x\in\partial\Omega$ and $y\in B_R(x)$,  
\begin{equation}\label{dirichlet thingy}
u(y,t)-u(x,t)\le  2M\varphi\left(\frac{\dist(y,\partial\Omega)}{2M},\frac t{4M^2}\right)
\le 2M\varphi\left(\frac{|x-y|}{2M},\frac t{4M^2}\right).
\end{equation}

This gives us an estimate on the boundary.  We complete our proof by using the same $Z$ as before: 
$$Z(x,y,t)=u(y,t)-u(x,t)-2M\varphi\left(\frac{|x-y|}{2M},\frac t{4M^2}\right).$$
Once again, $Z(x,y,0)\le 0$, and when both $x$ and $y$ are on the boundary, $Z(x,y,t)=-2M\varphi \le 0$. If $x$ is a boundary point and $y$ is an interior point (or vice-versa), then
\begin{equation*}
Z(x,y,t)=u(y,t)-u(x,t)-2M\varphi\left(\frac{|x-y|}{2M},\frac t{4M^2}\right)\le 0
\end{equation*}
by \eqref{dirichlet thingy}.

Finally, if $(x,y)$ is a maximum of $Z$ at some time $t<T'$, where both $x$ and $y$ are interior, then as in Section \ref{section periodic n>1} $\pd Z t\le0$ at this point, and so $Z\le0$ for all $t<T'$.  The estimate follows.
\halmos

\remark We can use these methods to find gradient estimates for equations of more general form.

For the Dirichlet problem with conditions on $a^{ij}$ given in Theorem \ref{periodic higherdimensional theorem}, and $u=0$ on $\partial\Omega$, we can find estimates of the type in Theorem \ref{mcf-dirichlet-theorem} for convex $\Omega$.

If $a^{ij}$ has the form \eqref{special form:  of mce-type?}, then we can find estimates of this type on domains that are merely mean-convex.

\newchapter{Application of gradient estimates to the Neumann problem}{Neumann problem}
\label{neumann chapter}

In this chapter we use the gradient estimate derived previously to establish the existence of solutions to the mean curvature flow equation with Neumann boundary conditions
\begin{equation}
\begin{split} \label{Neumann1}
\frac{\partial u}{\partial t}=
\left(\delta_{ij}-\frac{D_iu D_ju}{ 1+|Du|^2}\right)D_{ij}u 
\text{ in }\Omega\times(0,T],  \\
D_\nu u(x,t)=0\text{ for $x\in\partial\Omega$ and  }t\in(0,T] \\
\end{split}
\end{equation}
\begin{equation}
u(\cdot,0)=u_0, \label{Neumann2}
\end{equation}
where $\Omega \subset \bigR^n$ is a  compact, open convex domain  with $C^{2+\alpha}$ boundary 
 $\partial \Omega $, and $u_0\in C(\overline\Omega)$.   The outward unit normal on the boundary is $\nu$.   %

This extends Huisken's result in \cite{hu:boundary} showing the existence of smooth solutions to \eqref{Neumann1} for initial data with greater regularity.

\begin{theorem}[Huisken] \label{existing neumann result}
Let $\Omega$ be a bounded domain in $\bigR^n$ with $\partial\Omega\in C^{2+\alpha}$.  If $u_0 \in C^{2+\alpha}(\overline{\Omega})$ satisfies $D_\nu u_0=0$ on $\partial\Omega$, then \eqref{Neumann1}, \eqref{Neumann2} has a smooth solution on $\Omega\times(0,T)$.%
\end{theorem}

Note that while this theorem makes no restriction on the convexity of $\Omega$, the main result of this chapter does. %

\begin{theorem} \label{neumann theorem} 
Let $\Omega\subset \bigR^n$ be a smoothly bounded, open, convex domain, and let $u_0\in C(\overline{\Omega})$.
Then the Neumann problem \eqref{Neumann1} has a smooth solution for $t>0$, which converges uniformly to $u_0$ as $t\rightarrow 0$, at a rate dependent on the modulus of continuity of $u_0$. 
\end{theorem}

{\sloppy{
\section{Some remarks about changes of coordinates that straighten boundaries} \label{boundary section}
}}
A similar discussion of boundary curvatures and the distance function may be found in Appendix 14.6 of \cite{trudinger}.  Consider a bounded domain $\Omega\subset\bigR^n$ with boundary $\partial \Omega$.  The boundary is said to be $C^{k+\alpha}$ if for each boundary point $x_0$ we can find a $C^{k+\alpha}$ mapping $f:{\bigR}^{n-1}\rightarrow\bigR$ which has the boundary in a neighbourhood of $x_0$ as its graph.  

Set 
\begin{equation} R:=\frac 12 \sup\left\lbrace\, r: {\text{If }} 
\dist(x,\partial\Omega)<r 
\text{ then }
x
\text{  has a unique closest point } x_0\in\partial\Omega\,\right\rbrace.\label{definition of R}\end{equation}  
If $\Omega$ is convex, then we can take  
$$ R:=\frac 12 \inf_{\substack{
x\in\partial\Omega \\ 1\le i\le n-1 
}} \frac1{\kappa_i(x)}, $$ 
where $\kappa_i$ are the principal curvatures  of $\partial \Omega$.   %

On balls $B_R(x_0)$ centred on the boundary, we introduce a change of coordinates $\Psi:B_R\rightarrow \partial\Omega\times [-R,R]$ such that if the new coordinates are denoted $y=(\overline y,y^n)=(y^1,\dots,y^{n-1},y^n)$, $\Psi(x)=(\overline y, y^n)$, and $I:\partial \Omega \rightarrow {\mathbb R}^n $ is the immersion of the boundary into $\mathbb R^n$, then
$d(x,\partial\Omega)=d(x,I(\overline y))=y^n$.  In other words, $I(\overline y)$ is the closest point to $x$ on  $\partial\Omega$ and $y^n$ is the {\emph{signed distance}} between $x$ and $\overline y$,  being positive if $x\notin \Omega$, zero if $x \in \partial\Omega$,  and negative otherwise. 

The inverse transformation is easier to work with, being given by {$\Psi^{-1}(\overline y,y^n)=I(\overline y)+\nu(\overline y)y^n$}, where $\nu$ is the outward-pointing unit normal to $\partial\Omega$.   

As the boundary in the new coordinates is simply $y^n=0$, this is referred to as a boundary-straightening transformation.

\begin{center}
        \epsfig{file=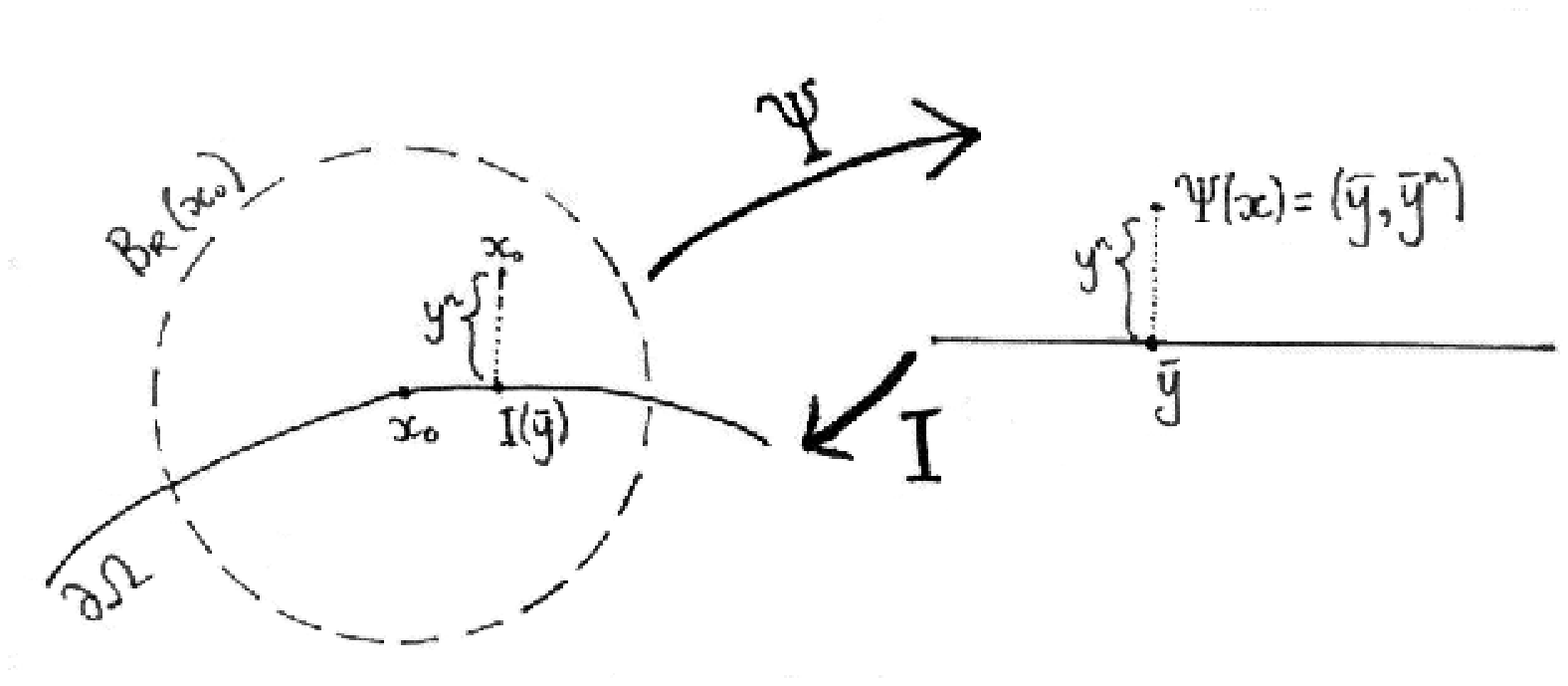, height=5cm}
\end{center}

If the graph $I(\overline y) = (y^1,\dots,y^{n-1},f(\overline y))$ is a local immersion of the boundary, then the outward unit normal is given by  \label{page where nu is second time around}
$$\nu(\overline y)=\frac{1}{\sqrt{1+|Df|^2}}\left(-\pd f {y^1},\dots,-\pd f{y^{n-1}},1\right), $$  and on $\Psi(B_R(x_0))$ we have
\begin{align}
[D\Psi^{-1}]_i^j&= \frac{\partial {\Psi^{-1}}^j}{\partial y^i} \notag \\
&=\begin{cases}
 { \delta_i^j+y^n\frac{\partial \nu^j}{\partial y^i}} &
 \text{  for }  i,j=1,\dots ,n-1 \\
 { \frac{\partial f}{\partial y^i}+y^n\frac{\partial \nu^n}{\partial y^i}} & 
 {\text{  for } j=n \text{  and }  i=1,\dots,n-1;} \\
{ \nu^j}  & \text{ for } i=n,
\end{cases} \label{Dphi^{-1} first derivatives}
\end{align}
so eigenvalues for $[D\Psi^{-1}]$ are $1-y^n\kappa_1(1+|Df|^2)^{-1/2},\dots,{1-y^n\kappa_{n-1}(1+|Df|^2)^{-1/2}}$, and lie between $\frac12$ and $\frac32$ on $B_{R}$.  The curvatures of the boundary, $\kappa_i$, are given by the eigenvalues of $[D^2f]$.

Also, second derivatives are 
$$
[D^2\Psi^{-1}]_{ik}^j= 
\begin{cases}
y^n\frac{\partial^2 \nu^j}{{\partial y^i}{\partial y^k}} &       
        \text{  for  $i,j,k\not=n$,} \\
\frac{\partial^2 f}{{\partial y^i}{\partial y^k}}
        +y^n\frac{\partial^2 \nu^n}{{\partial y^i}{\partial y^k}} 
        &\text{  for $ j=n $     and $i,k\not=n$,} \\
{\frac{\partial \nu^j}{\partial y^i}} & 
        \text{ for $i\not=n$ and $k=n$,} \\
{\frac{\partial \nu^j}{\partial y^k}} &
        \text{  for $i=n$ and $k\not=n$,} \\
0 &  \text{  for $ i=n, k=n$}. 
\end{cases}
$$
  
The smoothness of this change of coordinates is dependent on the smoothness of the boundary%
: if  $\partial\Omega$ is ${C^{k+\alpha+2}}$ then $D\Psi^{-1}$ is $C^{k+\alpha}$.  %

When $u$ is defined in the old coordinates on $\bigR ^n$, in the new coordinates we can define a new function
$$v(y,t):=u(\Psi^{-1}(y),t).$$ 

 First derivatives are related by $D_iv(y,t)=[D\Psi^{-1}]^k_iD_ku(\Psi^{-1}(y),t)$, and second derivatives by $[D^2 v]_{ij}=
 [D\Psi^{-1}]^m_j[D\Psi^{-1}]^l_i[D^2u]_{ml}+(D_mu)[D^2\Psi^{-1}]^m_{ij}$. 
 
 Putting this all together, we notice that if $u$ satisfies \eqref{Neumann1} then $v$ satisfies 
\begin{align} 
v_t(y,t)&=
m^{ij}([D\Psi]Dv)\left([D\Psi]^l_j[D\Psi]^k_i[D^2v]_{lk}+(D_sv)[D^2\Psi]^s_{ij}\right) \notag  \\
&= a^{kl}(Dv,y)[D^2v]_{lk}+b(Dv,y)
\qquad\text{  for $y\in\Psi(Q_R)$ }, \label{eqn4v}
\end{align}  
\begin{equation}
D_nv(y,t)=0  \qquad\text{  at $y^n=0$}.  \label{boundary condition part of eqn4v}
\end{equation}
 where the mean curvature operator is abbreviated as $m^{ij}(p)=\delta^{ij}-\frac{\delta^{ik}\delta^{jl}p_kp_l}{1+|p|^2}$, and we write 
\begin{equation} \label{full write-out of a and b}
\begin{split}
a^{lk}(p,y):&=m^{ij}([D\Psi]p)[D\Psi]^l_j[D\Psi]^k_i ,\\
b(p,y):&=m^{ij}([D\Psi]p)(p_s)[D^2\Psi]^s_{ij}.
\end{split}\end{equation}

Once we have straightened out the boundary, we will find it useful later on to define a reflection $\rho$ in the boundary that extends $v$ outside $\Omega$: 
$$\tilde{v}(y,t):=v(\rho(y),t) \quad\text{for $y\in B_{R}$} $$
where $\rho(y)=(y^1,\dots,y^{n-1},-|y^n|)$.

Let $Q_R$ be the intersection of a parabolic cylinder with the domain of interest:
\begin{equation*}
\left.Q_R(x_0,t_0)=\lbrace\,(x,t)\in\Omega\times[0,T] : x\in B_{R}(x_0),  t\in(t_0-R^2,t_0)\,\rbrace . \right)
\end{equation*}

When $v$ satisfies \eqref{eqn4v} on $Q_R$, $\tilde v$ will satisfy 
\begin{align} \label{eqn4v'}
\tilde v_t(y,t)&=v_t(\rho(y),t) \notag \\
&=a^{kl}(Dv(\rho(y),t),\rho(y))[D^2v(\rho(y),t)]_{kl} +b(Dv(\rho(y),t),\rho(y)) \notag\\
&=a^{kl}([D\rho]D\tilde v(y,t),\rho(y))[D\rho]^i_l[D\rho]^j_k[D^2\tilde v(y,t)]_{ij} +b([D\rho]D\tilde v(y,t),\rho(y)) \notag \\
&=\tilde a^{kl}(D\tilde v,y)[D^2\tilde v]_{kl}+\tilde b(D\tilde v,y)
\end{align}
on  $B_{R}\times(t_1-R^2,t_1)$ 
where 
$[D\rho]= \diag(1,\dots,1,-y^n/|y^n|)$.

The regularity of the coefficients of the reflected equation is estimated:
\begin{lemma} \label{lemma showing that reflected coefficients are c-alpha}
If $v$ is a $C^{1+\alpha}$ function on $Q_R$, with $D_n v(\overline y,0)=0$, then the coefficients for the reflected equation
\begin{gather*}
\tilde a^{ij}(D\tilde v,y)= a^{kl}([D\rho]D\tilde v,\rho(y))[D\rho]^i_l[D\rho]^j_k,  \\
\tilde b(D\tilde v,y)=b([D\rho]D\tilde v,\rho(y))
\end{gather*}
satisfy  H\"older estimates
\begin{gather*}
|\tilde a^{ij}|_{\alpha;B_R\times(t_1-R^2)}\le 2|a^{ij}|_{C^1;Q_R}\left(1+|Dv|_{\alpha;Q_R}\right), \\
 |\tilde b|_{\alpha;B_R\times(t_1-R^2)}\le 2|b|_{C^1;Q_R}\left(1+|D v|_{\alpha;Q_R}\right)
\end{gather*}
for  some $\,0<\alpha<1$.  
\end{lemma}

\longpage %
\proof 
In general, if a function is defined piecewise on a convex domain $U$ divided into $U_1$ and $U_2=U\backslash U_1$, 
\begin{equation*}
h(x)=\begin{cases} h_1(x) & x\in U_1 \\
h_2(x) & x\in U_2 
\end{cases}\end{equation*}
and is continuous across any shared boundary $\overline{U}_1\cap \overline{U}_2$ then if $h_1$ is $C^\alpha$ on $U_1$ and $h_2$ is $C^\alpha$ on $U_2$, it follows that $h$ is $C^\alpha$ on $U=U_1\cup U_2$. 

We can see this by letting $x_1$ and $x_2$ be in $U_1$ and $U_2$ respectively.  We can find a point $z$ in the shared boundary $\overline{U}_1\cap \overline{U}_2$   directly between the two, with $|y-x|=|y-z|+|z-x|$.  

Then 
\begin{align*}
|h(x)-h(y)|&\le|h(x)-h(z)|+|h(z)-h(y)| \\
&=|h_1(x)-h_1(z)|+|h_2(z)-h_2(y)| \\
&\le C|x-z|^\alpha + C|z-y|^\alpha  \\
&=C\left( s^\alpha |x-y|^\alpha + (1-s)^\alpha |y-x|^\alpha\right)  \\
&=2C|x-y|^\alpha 
\end{align*}
for $s=|x-z|/|x-y|<1$.  

This observation applies to both $\tilde v$ and $D\tilde v$ --- as $D_n v=0$  on the boundary, $D\tilde v=[D\rho][Dv]$ is continuous across the boundary, even though $[D\rho]$ itself is not --- and so $D\tilde v$ is $C^{\alpha}$.

It is clear that $\tilde b(D\tilde v,y)=b\left([D\rho]D\tilde v,\rho(y)\right)=b\left(Dv,\rho(y)\right)$ is continuous, and $\tilde b$ shares the same regularity as $b$. 

\longpage %

As
\begin{align*}
\tilde a^{ij}(D\tilde v,y)%
&= \begin{cases} a^{ij}\left([D\rho]D\tilde v(y,t),\rho(y)\right) & \text{ for } 1\le i,j\le n-1 \text{ or } i=j=n, \\
-\frac{y^n}{|y^n|}a^{ij}\left([D\rho]D\tilde v(y,t),\rho(y)\right), & \text{ for } i\not=n \text{ and }j=n\text{ or vice-versa, }
\end{cases}
\end{align*}
we only need to check whether the terms in the off-diagonal block  $\tilde a^{in}$ are continuous.  These are given by 
{\allowdisplaybreaks \begin{align*}
\tilde a^{in}&
(D\tilde v,y) \\
&=\frac{-y^n}{|y^n|}a^{in}\left([D\rho]D\tilde v(y,t),\rho(y)\right) \\*
&=\frac{-y^n}{|y^n|} 
m^{kl}\left([D\Psi][D\rho]D\tilde v\right)[D\Psi]^i_k[D\Psi]^n_l
\\&=  \frac{-y^n}{|y^n|} \left(  [D\Psi]^i_k[D\Psi]^n_k
\smallminus  \frac{\delta^{k\alpha}\delta^{l\beta} \left([D\Psi][D\rho]D\tilde v\right)_\alpha\left([D\Psi][D\rho]D\tilde v\right)_\beta[D\Psi]^i_k[D\Psi]^n_l}{1+ \left|[D\Psi][D\rho]D\tilde v\right|^2} \right)\\
&=\frac{y^n}{|y^n|}\left( \frac{
\delta^{k\alpha}\delta^{l\beta}
[D\Psi]^\gamma_\alpha[D\rho]^j_\gamma D_j \tilde v [D\Psi]^\mu_\beta[D\rho]^s_\mu D_s\tilde v[D\Psi]^i_k[D\Psi]^n_l} 
{1+ \left|[D\Psi][D\rho]D\tilde v\right|^2}\right) \\
&=\frac{y^n}{|y^n|} \frac1{1+ \left|[D\Psi][D\rho]D\tilde v\right|^2} 
\sum_{k,l,\alpha,\beta}\delta^{k\alpha}\delta^{l\beta}\Bigg(
\sum_{\gamma=1}^{n-1} \sum_{\mu=1}^{n-1}
[D\Psi]^\gamma_\alpha D_\gamma \tilde v [D\Psi]^\mu_\beta D_\mu\tilde v
\\*& \phantom{space}
- \frac{y^n}{|y^n|} \sum_{\gamma=1}^{n-1}  
[D\Psi]^\gamma_\alpha D_\gamma \tilde v [D\Psi]^n_\beta D_n\tilde v 
- \frac{y^n}{|y^n|} \sum_{\mu=1}^{n-1}
[D\Psi]^n_\alpha D_n \tilde v [D\Psi]^\mu_\beta D_\mu \tilde v
\\*&\phantom{spacespacespacespacespacespacespace}
+ [D\Psi]^n_\alpha D_n \tilde v [D\Psi]^n_\beta D_n\tilde v
\Bigg) [D\Psi]^i_k[D\Psi]^n_l,
\end{align*} (here there is no summation over $n$).} 

Between the third and the fourth line, we have used that $[D\Psi^{-1}]^{\cdot}_i$ is tangent to the boundary while $[D\Psi^{-1}]^{\cdot}_n$ is normal to the boundary (see equation \eqref{Dphi^{-1} first derivatives}), so  for $i\not=n$, we have $\sum_k [D\Psi^{-1}]^k_i[D\Psi^{-1}]^k_n=0$ .  It follows that $[D\Psi]^i_k[D\Psi]^n_k=0$.

In the last step, the second, third and fourth terms are zero (and so continuous) on the boundary, 
as $D_n\tilde v(\overline y,0)=0$.    The first term is zero due to the presence of 
$$\sum_l \left([D\Psi]^\mu_l D_\mu v\right)[D\Psi]^n_l=0,$$
since $\mu\not=n$.

So, both $\tilde a^{ij}$ and $\tilde b$ are continuous across the boundary.

It follows from the first observation that the H\"older constant of $\tilde a^{ij}$ on $B_R(y_1)\times(t_1-R^2,t_1)$ is the same as that of  $a^{ij}$ on $Q_R(y_1,t_1)$; and if we consider $a^{ij}(Dv(y),y)$ as a function of $y$, then we find that
\begin{align*}
|\tilde a^{ij}(D\tilde v(\cdot),\cdot)|_\alpha &\le  2  | a^{ij}(D v(\cdot),\cdot)|_\alpha \\
&= 2\sup_{z_1,z_2}\frac{|a^{ij}(Dv(z_1),z_1)-a^{ij}(Dv(z_2),z_2)|}{|z_1-z_2|^\alpha} \\
&\le 2\sup_{z_1,z_2}\frac1{|z_1-z_2|^\alpha}\Big[{|a^{ij}(Dv(z_1),z_1)-a^{ij}(Dv(z_1),z_2)|}
\\&\phantom{spacespacespacespacespace}+{|a^{ij}(Dv(z_1),z_2)-a^{ij}(Dv(z_2),z_2)|} \Big] \\
&\le2\sup_{z_1}|a^{ij}(Dv(z_1),\cdot)|_\alpha 
\\&\phantom{spacespace}+ \sup_{z_2}|a^{ij}(\cdot,z_2)|_{C^1} \sup_{z_1,z_2}\frac1{|z_1-z_2|^\alpha}|Dv(z_1)-Dv(z_2)| \\
&\le 2|a^{ij}|_{C^1}\left(1+|Dv|_\alpha\right).
\end{align*}

Similarly, if we consider $\tilde b(D\tilde v,y)$ as a function of $y$, we find that
\begin{align*}
|\tilde b (D\tilde v(\cdot),\cdot)|_\alpha&\le  2  | b (D v(\cdot),\cdot)|_\alpha \\
&\le 2|b|_{C^1}\left(1+|D v|_\alpha\right).
\end{align*}
\halmos

 \section{Existence of solutions with continuous initial data and Neumann boundary conditions} \label{the section where neumann existence is proved} 

We begin our proof of Theorem \ref{neumann theorem} by approximating the continuous initial data by mollified functions that will satisfy the requirements of Theorem \ref{existing neumann result}, being smooth and satisfying the Neumann boundary condition.  

\begin{lemma} \label{neumann existence of approximating sequence} There exists an approximating sequence $u_0^\epsilon\in C^\infty(\overline\Omega)$ with $u_0^\epsilon\rightarrow u_0$ in $C(\Omega)$, $D_\nu u_0^\epsilon=0$ on $\partial\Omega$, and $u^{\epsilon_1}>u^{\epsilon_2}>u_0$ whenever $\epsilon_1>\epsilon_2$.   
\end{lemma}

\proof Let $B_R$ be a ball centred on the boundary.
We work in the new coordinates on $\Psi(B_{R})$, and write $v_0(y)=u_0(\Psi^{-1}(y))$.

Remembering that $\tilde v_0$ denotes the extension by reflection of $v_0$, define the mollified function
$$ 
v_0^\epsilon(y):=\left.\eta_\epsilon * \tilde v_0\right. (y) = \int_{z\in\Psi(B_R)}\eta_\epsilon(y-z){\tilde v_0}(z) dz, $$
where we use the usual mollifier  
\begin{equation*}
 \eta(z)= 
\begin{cases}
c \exp\left(\frac1 {|z|^2-1}\right), &  |z|\le 1 \\
 0 & \text{otherwise, } 
\end{cases}
\end{equation*}
and $\eta_\epsilon(z)=\frac1{\epsilon ^n} \eta(\frac z \epsilon)$. 

This approximation has all the usual qualities of mollifications: $v_0^\epsilon \in C^\infty(\Psi(B_R)^\epsilon)$, 
where $\Psi(B_R)^\epsilon=\{\,y\in\Psi(B_R):\dist(y,\partial\Psi(B_R))>\epsilon\,\}$
; and since ${\tilde{v}_0}\in C(\Psi(B_R))$, $v_0^\epsilon\rightarrow \tilde v_0 $ uniformly on compact subsets of $\Psi(B_R)$. 

In addition, each $v_0^\epsilon$ satisfies the Neumann condition $D_nv_0^\epsilon(y)=0$ for $y\in\partial\Omega\cap \Psi(B_R)^\epsilon$, since 
\begin{align*}
D_nv_0^\epsilon(y)&=D_n\int_{z\in B_\epsilon(0)}\eta_\epsilon(z)\tilde v_0(y-z)dz \\
&=\int_{\substack{
 z\in B_\epsilon(0)
\\ z^n\ge0
 }}\eta_\epsilon(z)D_n\tilde v_0(y-z)dz + \int_{\substack{
z\in B_\epsilon(0)
\\ z^n<0
 }}\eta_\epsilon(z)D_n\tilde v_0(y-z)dz 
\end{align*} 
Recalling the relationship between the reflected and original functions, $\tilde v_0(y-z)=v_0(y^1-z^1,\dots,-|y^n-z^n|)$, we observe that when $y\in\partial\Omega$, $y^n=0$ and $D_n\tilde v_0(y-z)= \frac{z^n}{|z^n|}D_n v_0(y^1-z^1,\dots,-|z^n|)$.   Consequently,
\begin{align*}
D_nv_0^\epsilon(y)&=\int_{\substack{
 z\in B_\epsilon(0)
\\ z^n>0
 }}
\eta_\epsilon(z)D_nv_0(y^1-z^1,\dots,-|z^n|)dz \\
&\phantom{spacespacespace}- \int_{\substack{
 z\in B_\epsilon(0)
\\ z^n<0
 }}
\eta_\epsilon(z)D_n v_0(y^1-z^1,\dots,-|z^n|)dz
=0 
\end{align*}
as the mollifier has the symmetry $\eta_\epsilon(z^1,\dots,z^n)=\eta_\epsilon(z^1,\dots,-z^n)$.  

This is only a \emph{local} approximation, but in the next step we extend it to the entire domain, taking care to preserve the Neumann boundary condition.

Let the set of boundary-centred balls $\lbrace \, B_{R}(x_i)\,\rbrace_{i=1,N}$ be a finite cover of the boundary $\partial\Omega$ with the property that the set of balls of half the radius $\lbrace\,B_{R/2}(x_i)\,\rbrace_{i=1,N}$ is also a cover.  On each ball $B_R(x_i)$ we can define the approximation  ${{v^\epsilon_{0,i}}}:=v_0^\epsilon$ as described above.

Now, define a new cover of $\Omega$ by the sets 
$$W^i:=\{\,x\in B_R(x_i):\text{ if }\Psi(x)=(\overline y,y^n),\text{ then }\overline y\in \Psi(B_{R/2}(x_i))\text{  and }|y^n|\le \frac R2\,\}.$$    
The cover is completed by $W^0:=\{\,x\in \Omega:\text{dist}(x,\partial\Omega)>R/4\,\}$.   

Note that $B_{R/2}(x_i)\subseteq W^i \subseteq B_R(x_i)$ and so this is indeed a cover; also, ${{v^\epsilon_{0,i}}}$ is defined on $W^i$.  On $W^0$ we define the usual mollification with no reflection, which we call ${u^\epsilon_{0,0}}$.

Let $\overline\xi_i$ be a partition of unity with respect to the sets $\lbrace B_{ R/2}(x_i)\cap \partial\Omega\rbrace$  which cover the boundary; that is, $0\le\overline\xi_i(x)\le 1$ for $x\in B_{ R/2}(x_i)\cap\partial\Omega$, $\overline\xi_i\in C^\infty_0(B_{ R/2}\cap \partial \Omega)$ (that is, compactly supported with respect to $\partial\Omega$), and $\sum \overline\xi_i(x)=1$ for all $x\in\partial \Omega$.

In the new coordinates 
on $B_R(x_i)$, 
we could write $\overline\xi_i=\overline\xi_i(y^1,\dots,y^{n-1})$, since $\overline\xi_i$ is defined only on the boundary.   We can extend $\overline\xi_i$ to {\emph{all}} of $B_{R/2}(x_i)$ by setting 
$\xi_i(x):=\overline\xi_i(\Psi(x)^1,\dots,\Psi(x)^{n-1})$.  

Let $\tilde \zeta:\bigR\rightarrow\bigR$ be a smooth cut-off function satisfying \begin{equation*}
\tilde\zeta(d)=\begin{cases}
1 & |d|< \frac R4  \\
0 & |d| \ge \frac R2.
\end{cases}
\end{equation*}
We will set $\zeta(x):=\tilde\zeta(d(x,\partial\Omega))$ 
where $d(\cdot,\partial\Omega)$ is the signed distance function.  

Now, we claim that the functions $\xi_i(x)\zeta(x)$ and $1-\zeta(x)$ are a partition of unity with respect to the sets $\lbrace W^i\cap\Omega                                                                                                                             \rbrace$ and $W^0$.  Firstly, all functions are smooth and compactly supported on their respective domains (but they are not zero on the external boundary $\partial \Omega$).   Secondly,  if $x\in\Omega$, then
$$\sum_{i}\xi_i(x)\zeta(x)+(1-\zeta(x))=1.$$  
 This is because if  $\dist(x,\partial\Omega)\ge R/2$ then $\zeta(x)=0$, while if $\dist(x,\partial\Omega)< R/2$ then $x$ has a unique closest boundary point $x_0$.  In the latter case $\xi_i(x)=\overline \xi_i(x_0)$ and $\sum_{i}\overline\xi_i(x)\zeta(x)+(1-\zeta(x))=\zeta(x)\sum_{i}\overline\xi_i(x_0)+(1-\zeta(x))=1$ as $\overline\xi_i$ is a partition of unity on the boundary.

This construction ensures that $D_\nu\left(\zeta(x)\xi_i(x)\right)=0$ when $x\in\partial\Omega$ and so if we define our global approximation as 
$$u_0^\epsilon(x):=\sum_{i=1}^N \xi_i(x)\zeta (x){{v^\epsilon_{0,i}}}(\Psi(x)) + (1-\zeta(x))u^\epsilon_{0,0}(x),$$ we find that $u_0^\epsilon\rightarrow u_0$ uniformly in $C(\Omega)$, and each $u_0^\epsilon \in C^\infty(\overline{\Omega})$ satisfies $D_\nu u_0^\epsilon=0$ on $\partial \Omega$.

We can ensure that this sequence is monotone in $\epsilon$, in the sense that $u_0^{\epsilon_1}(x)<u_0^{\epsilon_2}(x)$ whenever $\epsilon_1<\epsilon_2$ by restricting to a subsequence and off-setting if necessary.  %
\halmos

The result of Huisken mentioned at the start of this chapter now implies that there is a smooth solution $u^\epsilon$ to \eqref{Neumann1} with $u^\epsilon(\cdot,0)=u^\epsilon_0$.  

\begin{lemma}
The approximate solutions $u^\epsilon$ have a uniform height bound
$$\sup_{\Omega\times[0,T]} |u^\epsilon|\le \sup_\Omega |u_0|. $$ 
\end{lemma}
\proof As the mollification $u^\epsilon_0$ is created by a local averaging of $u_0$, 
$$\sup_\Omega |u^\epsilon(\cdot,0)|\le \sup_\Omega |u_0|.$$
Suppose at some time $t>0$ and point $x_1$, $u^\epsilon$ equals $|u_0|$.   From the Comparison Principle (Theorem \ref{comparison 1}), $x_1$  can be assumed to be a boundary point. The Neumann condition $D_\nu u^\epsilon=0$ implies that $Du^\epsilon=0$ and so $[D^2 u^\epsilon]$ is negative semi-definite at this point; it follows that $\pd {u^\epsilon}{t}\le0$ and so $u^\epsilon$ is not increasing at this point.    
\halmos

This height estimate is of course also an oscillation bound 
$$|u^\epsilon(x,t)-u^\epsilon(z,t)|\le 2|u_0|.$$

We are now in a position to use the gradient estimate of Theorem \ref{Neumann theorem}.
For some $T>0$, there is  a gradient bound
\begin{equation}|Du^\epsilon(x,t)|\le L(t) \text{ for $t\in(0,T)$ and $x\in\overline\Omega$} \label{main grad bound in neumann chapter}, \end{equation}
where $L(t)$ and $T$ are dependent on $n$ and $\sup_\Omega |u_0|$.

\begin{lemma} \label{higher derivatives bounded on interior} Higher derivatives of $u^\epsilon$ are uniformly bounded on the interior, with 
\begin{equation}|D^{k}u^\epsilon|_{\Omega^r\times\lbrace t_1\rbrace}\le c\left(n,k,L(t_0)\right)
\left(\frac1{%
r^2}+\frac1{t_1-t_0}\right)^{\frac{k-1}2} \label{interior higher derivatives estimate}
\end{equation}
for $t_1>t_0>0$ and all $k=1,2,\dots$, where 
$\Omega^r$ is the interior set $\lbrace\,x\in\Omega:\dist(x,\partial\Omega)>r\,\rbrace$. 
\end{lemma} 

\proof This is an application of the Ecker-Huisken interior curvature estimate described in Theorem \ref{interior-higher-estimate}, originally in the paper \cite{eh:interior}.

We  apply it to the interior of $\Omega$ (with $\theta=0$ and $k=m+2$) to find bounds on all higher derivatives.%
\halmos

This estimate provides no information as we approach the boundary.  However, our uniform gradient bound $L(t_0)$ ensures that the evolution equation is uniformly parabolic, since for $t>t_0$,
\begin{equation*}
m^{ij}(Du^\epsilon)\xi_i\xi_j
\,\ge\,\frac1{1+|Du^\epsilon|^2}|\xi|^2 
\,\ge\,\frac1{1+L(t_0)^2}|\xi|^2. 
\end{equation*}

As we have uniform parabolicity for strictly positive times,  extending regularity up to the boundary  is a routine  application of known results.  This is the subject of Lemma \ref{linking lemma} -- Lemma \ref{lemma 6.9}.  

We begin by showing that a function with a H\"older estimate on the boundary of a region, and a strictly interior gradient estimate, has a global H\"older estimate.   We plan to apply this to finding a H\"older estimate for the gradient $Du$.

\begin{lemma} \label{linking lemma}Let $\Omega$ be a convex domain.  If $f:\Omega\times[0,T]\rightarrow\bigR$ has a H\"older oscillation bound on the boundary
\begin{gather*}
\osc_{Q_r(x_0,t_0)}f \le C({t_0})r^\alpha \text{ for all  }x_0\in\partial\Omega \text{ and } {t_0}> r^2, \\
\intertext{ where $C(t)$ is non-increasing in $t$; and gradient bounds on the interior} 
|Df(x,t)|\le c\left(\frac 1 {\dist(x,\partial\Omega)^2}+ \frac{1}t\right)^{1/2}; \intertext{and} 
|f_t(x,t)|\le c\left(\frac 1 {\dist(x,\partial\Omega)^2}+ \frac{1}t\right);
\end{gather*}
then we can find an $\alpha'>0$ such that for all $x,y\in\Omega$ and $ s,t\in(0,T]$,
\begin{equation*}
|f(x,t)-f(y,s)|\le C(|x-y|+{|t-s|}^{1/2})^{\alpha'} 
\end{equation*}
where $C$ depends on $\min\lbrace t,s\rbrace$, $\diam{\Omega}$, $c$ and $\alpha$, and $\alpha'$ depends on $\alpha$.  
\end{lemma}
\proof  We split the difference in the obvious way 
\begin{equation}|f(x,t)-f(y,s)|\le|f(x,t)-f(y,t)|+|f(y,t)-f(y,s)|,
\label{linking}
\end{equation} 
and look at the first term.

Without loss of generality, set $d=\dist(y,\partial\Omega)\le \dist(x,\partial\Omega)$, and $y_0$ to be the closest point to $y$ in $\partial\Omega$, so that $|y-y_0|=d$.

If we are close to the boundary, so that $d\le |x-y|$, then 
\begin{align*}
|f(x,t)-f(y,t)|&\le|f(x,t)-f(y_0,t)|+|f(y_0,t)-f(y,t)| \\
&\le \osc_{Q_{|x-y_0|}(y_0,t+|x-y_0|^2)} f +  \osc_{Q_{|y-y_0|}(y_0,t+|y-y_0|^2)} f \\ 
&\le C(t+|x-y_0|^2)|x-y_0|^\alpha +  C(t+|y-y_0|^2)|y-y_0|^\alpha \\
&\le C(t)\left[|x-y|^\alpha +|y-y_0|^\alpha\right]  + C(t)|y-y_0|^\alpha \\
&\le 2 C(t)|x-y|^\alpha + C(t)|y-x|^\alpha \\
&\le  c_1  |x-y|^\alpha
\end{align*} 
where we have used that $|y-y_0|=\dist(y,\partial\Omega)\le |x-y|$, and have set $c_1=3C(t)$. 

If we are further from the boundary, so that $d>|x-y|$, then for some $\epsilon\in(0,1)$,
\begin{align*}
|f(x,t)-f(y,t)|&= |f(x,t)-f(y,t)|^\epsilon \,  |f(x,t)-f(y,t)|^{1-\epsilon} \\
&\le  \left[ \osc_{Q_{|x-y_0|}(y_0,t+|x-y_0|^2)} f \right]^\epsilon 
\left[ c\left(\frac{1}{\dist(y,\partial\Omega)^2}+\frac1t\right)^2|x-y|\right]^{1-\epsilon}\\ 
&\le \left[ C(t+|x-y_0|^2) |x-y_0|^\alpha\right]^\epsilon  
 c^{1-\epsilon} \left(\frac{1}
{d^2}+\frac1t\right)^{\frac{1-\epsilon}2}|x-y|^{1-\epsilon} \\
&\le C(t)^\epsilon \left[ |x-y|^\alpha +|y-y_0|^\alpha\right]^\epsilon  
 c^{1-\epsilon}\left({d^{\epsilon-1}}+{t^{(\epsilon-1)/2}}\right)|x-y|^{1-\epsilon}\\
&\le  2^\epsilon C(t)^\epsilon \left({d^{\alpha\epsilon+\epsilon-1}}+d^{\alpha\epsilon}{t^{(\epsilon-1)/2}}\right)|x-y|^{1-\epsilon} \\
&\le c_2|x-y|^{1-\epsilon},
\end{align*}
setting ${\epsilon=1/{(\alpha+1)}}$ so that ${{\alpha\epsilon+\epsilon-1}=0}$,
 and %
$c_2 = 2^\epsilon C(t)^\epsilon \left(1+(\diam{\Omega})^{\alpha\epsilon}{t^{(\epsilon-1)/2}}\right)$.

Now consider the second term of \eqref{linking}, and suppose without loss of generality that $s<t$.  As before, set $d=\dist(y,\partial\Omega)$ and $y_0\in\partial\Omega$.

If we are close to the boundary, so that $d\le\sqrt{t-s}$, then
\begin{align*}
|f(y,t)-f(y,s)|&\le \osc_{Q_{\sqrt {t-s}}(y_0,t)} f\\
& \le C(t)|t-s|^{\alpha/2}  \\
&= c_3|t-s|^{\alpha/2}
\end{align*}

Otherwise, if $d>\sqrt{t-s}$, then
\begin{align*}
|f(y,t)-f(y,s)|&= |f(y,t)-f(y,s)|^\mu \,  |f(y,t)-f(y,s)|^{1-\mu} \\
&\le  \left[ \osc_{Q_d(y_0,t)} f \right]^\mu 
\left[c\left(\frac{1}{d^2}+\frac1s\right)|t-s|\right]^{1-\mu}\\ 
&\le \left[ C(t) d^\alpha\right]^\mu  
 c^{1-\mu} \left(\frac{1}
{d^2}+\frac1s\right)^{{1-\mu}}|t-s|^{1-\mu} \\
&\le C(t)^\mu c^{1-\mu} \left(d^{\alpha\mu-2(1-\mu)}+ \frac{d^{\alpha\mu}}{s^{1-\mu}}\right) |t-s|^{1-\mu} \\
&\le c_4%
|t-s|^{(2-2\mu)/2}
\end{align*}
where $\mu=2/(\alpha+2)$ so that $\alpha\mu-2(1-\mu)=0$, and  $c_4\smallequals C(s)^\mu c^{1-\mu}\mspace{-0.9mu}\left(1\smallplus(\diam{\Omega})^{\alpha\mu}s^{\mu-1}\right)$.  

We find the final estimate by choosing $C=\sup\lbrace c_1,c_2,c_3,c_4\rbrace$ and $\alpha'=\min\lbrace 1-\epsilon ,\linebreak[2] 2-2\mu \rbrace=\min\lbrace \alpha/(\alpha+1), 2\alpha/(\alpha+2)\rbrace=\alpha/(\alpha+1).$
\halmos

\begin{lemma} \label{holder gradient bound lemma} The gradient $Du^\epsilon$ is H\"older continuous, with bound 
$$|Du^\epsilon|_{\alpha,\alpha/2}\le C$$
on $\Omega\times[t,T]$, for $t>t_0>0$ and some $\alpha>0$, where $C=C(n,L(t_0),|t-t_0|,\Omega,|\partial\Omega|_{C^2})$.
\end{lemma}

\proof We can use the H\"older gradient estimate near a flat boundary from  Theorem \ref{Lieb}, but we will need to work locally with $v^\epsilon(y,t)=u^\epsilon(\Psi^{-1}(y),t)$ in the flat-boundary coordinates.

The gradient bound \eqref{main grad bound in neumann chapter}
for $u^\epsilon$ implies that %
$$|Dv^\epsilon|_{Q_r}\,\le\, L(t_0) |D\Psi^{-1}|_{B_r}\,\le \,L(t_0)|\partial\Omega|_{C^2}$$
on the cylinder $Q_{r}(y_1,t_1)$ for  some $r<R$ to be chosen later, and
where $y_1\in\partial\Omega$ and $t_1>t_0+{r}^2$. 

On $Q_r$, $v^\epsilon$ satisfies the evolution equation \eqref{eqn4v}.  
To check that the coefficients $a^{lk}$
and $b$
(given by \eqref{full write-out of a and b})
satisfy the conditions of Theorem \ref{Lieb}, we note that:  {\allowdisplaybreaks[1]{$a^{lk}$ is uniformly parabolic, since for $\xi\in\mathbb R^n$,
\begin{align*}
a^{lk}\xi_l\xi_k
&=m^{ij}([D\Psi]Dv^\epsilon)[D\Psi]^l_j[D\Psi]^k_i\xi_l\xi_k \\*
&\ge  \frac1{1+|[D\Psi]Dv^\epsilon|^2} |[D\Psi]\xi|^2 \\*
&\ge \frac1{1+\Lambda^2_{[D\Psi]}}\lambda^2_{[D\Psi]}|\xi|^2\\*
&\ge \frac1{1+4L(t_0)^2}\left(\frac23\right)^2|\xi|^2 \\*
&\ge c(L(t_0))|\xi|^2
\end{align*}}}
where $\lambda_{[D\Psi]}$ and $\Lambda_{[D\Psi]}$ are the smallest and largest eigenvalues of $[D\Psi]$, which are bounded between $2/3$ and $2$ on $B_R$; 
$|a^{lk}|$ is bounded above, as
$$|m^{ij}([D\Psi]Dv^\epsilon)[D\Psi]^l_j[D\Psi]^k_i| \,\le\, |\Lambda_{[D\Psi]}|^2;%
$$
and $a^{lk}$ has bounded derivative with respect to the gradient, for if we write
 $q=[D\Psi]p$, then   $a^{lk}(p,y)=m^{ij}(q)[D\Psi]^l_j[D\Psi]^k_i$ and
\begin{align*} 
\frac{\partial a^{lk}(p,y)}{\partial p^\alpha}&=
[D\Psi]^l_j[D\Psi]^k_i \frac{\partial m^{ij}(q)}{\partial q^\beta}\frac{\partial q^\beta}{\partial p^\alpha} \\
&=[D\Psi]^l_j[D\Psi]^k_i \left[-\frac{(\delta_{i}^\beta q_j+\delta_{j}^\beta q_i)}{1+|q|^2}+\frac{2q_\beta q_i q_j}{(1+|q|^2)^2}\right][D\Psi]^\beta_\alpha \\
&\le 4 |D\Psi|^3 \\&\le 4 |\partial\Omega|^3_{C^2}. 
\end{align*}         

It is straightforward to bound the lower-order term in the equation ---  
 $$|m^{ij}([D\Psi]Dv^\epsilon)[D^2\Psi]^s_{ij}D_sv^\epsilon|\,\le\, L(t_0)|D^2\Psi|\,\le\, L(t_0)|\partial\Omega|_{C^3}.$$

Finally, we need an oscillation bound smaller than $\sigma$ for $a^{lk}(p,\cdot)$ on $Q_r$, but since
\begin{align*}
\frac{\partial a^{lk}(p,y)}{\partial y}
&=\frac{\partial m^{ij}([D\Psi]p)}{\partial q^l}\frac{\partial([D\Psi]p)^l}{\partial y}[D\Psi]^l_j[D\Psi]^k_i
+m^{ij}\frac{\partial}{\partial y}\left([D\Psi]^l_j[D\Psi]^k_i\right)\\
&\le C|D^2\Psi||D\Psi|\\
&\le C|\partial\Omega|^2_{C^3},
\end{align*} 
we find that $\osc_{Q_r(y_1,t_1)}a^{ij}(p,\cdot)\le r C |D^2\Psi||D\Psi|$.  By choosing $r$ small enough, we can ensure that this is less than $\sigma$.

Now Theorem \ref{Lieb} implies that for all $s<r$ there is some $\alpha'>0$ so that
\begin{align*} 
\osc_{Q_s} D v^\epsilon 
& \le c\left(\frac s{r}\right)^{\alpha'}
    \left(\osc_{Q_{r}} Dv^\epsilon + L(t_0)|\partial\Omega|_{C^2} r\right)\\
&\le c\left(\frac s{r}\right)^{\alpha'}
    \left(L(t_0) + L(t_0)|D^2\Psi| r\right)\\
& \le c s^{\alpha'} r^{-\alpha'}(1+r).
\end{align*}
where $c(L(t_0),n,|D\Psi|,|D^2\Psi|)$.   

This boundary oscillation estimate for $Dv^\epsilon$ on $Q_s(y_1,t_1)$ for all $s<r<R$ and $t_1>t_0+r^2$, together with the interior gradient bounds given by Lemma \ref{higher derivatives bounded on interior}, means we can use Lemma \ref{linking lemma} to give a global H\"older bound for $Dv^\epsilon$ and hence for $Du^\epsilon$. 
\halmos

\begin{lemma}\label{lemma 6.9} We can find bounds for  $u^\epsilon$  in  $H_{2+\alpha}\left(\Omega\times(t,T]\right)$ 
$$|u^\epsilon|_{2+\alpha,1+\frac\alpha2}\le C, $$  
for $t>t_0>0$, where $C=C(n,L(t_0),|t-t_0|,|\partial\Omega|_{C^2},\alpha)$. 
\end{lemma}

\proof To establish this, it is possible to use boundary estimates for the Neumann problem, but instead our approach is to use the reflection  $\tilde v^\epsilon$ on \mbox{$B_{R}(y_1)\times[t_0,T]$ ---} a domain that extends beyond $\Omega$ --- which satisfies the reflected evolution equation \eqref{eqn4v'}, and apply the interior estimate from Theorem \ref{interior1}.

We need to check that equation \eqref{eqn4v'} satisfies the conditions of Theorem \ref{interior1}.
 
In Lemma \ref{lemma showing that reflected coefficients are c-alpha}, we showed that the coefficients in equation \eqref{eqn4v'} have regularity estimates
\begin{gather*}
|\tilde a^{ij}|_{\alpha;B_r\times(t_1-r^2)}\le 2|a^{ij}|_{C^1;Q_r}\left(1+|Dv^\epsilon|_{\alpha;Q_r}\right) \\
 |\tilde b|_{\alpha;B_r\times(t_1-r^2)}\le 2|b|_{C^1;Q_r}\left(1+|Dv^\epsilon|_{\alpha;Q_r}\right);
\end{gather*}
and in Lemma \ref{holder gradient bound lemma} we found a uniform global bound for $|D v^\epsilon|_{\alpha,\alpha/2;\Omega\times(t,T)}$ for $t>t_0$.

Our gradient estimate ensures that $\tilde a$ is uniformly  parabolic for $t>t_0$.

Applying Theorem \ref{interior1} results in the bound 
\begin{equation*}
|\tilde v^\epsilon|_{2+\alpha,1+\alpha/2 ; Q_{ {R}/3}}\le C |\tilde v^\epsilon|_{0;Q_{2R/3}},\end{equation*}
where $C$ is dependent on the dimension $n$, $R$ (which is determined by $|\partial\Omega|_{C^2}$),  the ellipticity constant of $\tilde a^{ij}$ (dependent on $L(t_0)$), 
the H\"older exponent $\alpha$, and the bound on the $\alpha$ norm of the coefficients (which is bounded by $|Dv^\epsilon|_{\alpha,\alpha/2}$).   That is, $C=C(n,L(t_0),|\partial\Omega|_{C^2},\alpha)$. 

We can repeat this over the entire boundary, and together with the interior estimate \eqref{interior higher derivatives estimate}, this gives us the claim.
\halmos

\begin{lemma}
The sequence of approximate solutions $u^\epsilon$ converges
\begin{equation*}
\lim_{\epsilon\rightarrow 0}  |u^\epsilon-u|_{2+\alpha',1+\alpha'/2}\rightarrow 0\end{equation*}
 to some $u\in H_{2+\alpha'}$  on $\Omega\times (t_0,T)$, for all $t_0>0$.  
 \end{lemma}

\proof
The uniform 
$H_{2+\alpha}(\overline\Omega\times(t_0,T))$ bounds on the $u^\epsilon$ 
ensure that there is a convergent subsequence (for a slightly smaller $\alpha'<\alpha$); the disjointness of initial data $u^\epsilon_0$ (and hence the disjointness of $u^\epsilon(\cdot,t)$) implies that the entire sequence must converge, and so this limit is unique.   The limit $u$ is in $H_{2+\alpha'}(\Omega\times(t_0,T))$ (and is $C^\infty$ on the interior, by virtue of the interior estimate in Lemma \ref{higher derivatives bounded on interior}).

 It also satisfies $D_\nu u=0$ on the boundary, and so is a solution to the Neumann problem given by \eqref{Neumann1}.%
\halmos

\begin{lemma} \label{convergence to initial time}This solution $u$ converges to $u_0$  in $C(\Omega)$ as $t\rightarrow 0$.   The convergence is uniform in time, and if $u_0\in C^\alpha(\Omega)$ then $u\in H_\alpha(\Omega\times[0,T])$.  
\end{lemma}

\proof  Let $\delta>0$ be fixed.  Our aim is to show that we can find $t_\delta$ so that for all $t\in(0,t_\delta)$,
$$\sup_{x\in\Omega}|u(x,t)-u_0(x)|\le \delta.$$

The modulus of continuity for $u_0$ is
$$\omega(r):=\sup_{|x-y|=r} |u_0(x)-u_0(y)|.$$

Recall the approximate solutions $u^\epsilon$, converging uniformly
$$|u^\epsilon(\cdot,t)-u(\cdot,t)|\rightarrow 0 \text{ as }\epsilon\rightarrow 0$$
for all $t>0$.  
The approximations at the initial time have (at least) the same modulus of continuity as $u_0$:
\begin{equation*}
|u^\epsilon(x,0)-u^\epsilon(y,0)|\le \omega(|x-y|).
\end{equation*}

\begin{center}
\begin{figure}[h]
\centering
\epsfig{file=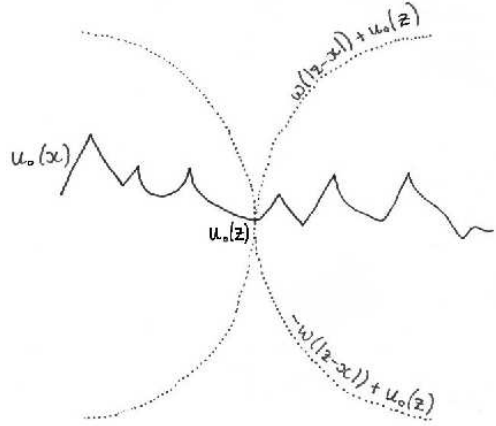, height=6cm}
\caption{The modulus of continuity bounds $u_0$}
\end{figure}
\end{center}

Fix $z$ to be any point in the interior of $\Omega$.    We can define a new solution to \eqref{Neumann1}  by off-setting $u^\epsilon$ around $u_0(z)$: 
$$w^\epsilon(x,t):=u^\epsilon(x,t)-u^\epsilon(z,0)+u_0(z),$$ so that
$w^\epsilon(z,0)=u_0(z)$. 

Then, for $t>0$, we have
\begin{align}
|u(z,t)-u_0(z)|&= \lim_{\epsilon\rightarrow 0}|u^\epsilon(z,t)-u_0(z)| \notag \\ 
&=\lim_{\epsilon\rightarrow 0}|w^\epsilon(z,t)+u^\epsilon(z,0)-2u_0(z)| \notag \\
&\leq \lim_{\epsilon\rightarrow 0} \left(|w^\epsilon(z,t)-u_0(z)| +   
|u^\epsilon(z,0)-u_0(z)|\right). \label{significant eqn in final lemma of neumann chapter}
\end{align} 

To estimate the first term of this, we
observe that every $w^\epsilon(\cdot,0)$ is inside the `envelope' given by the continuity condition, $$-\omega(|z-x|)+u_0(z)\leq w^\epsilon(x,0) \leq
\omega(|z-x|)+u_0(z).$$

\begin{center}
\begin{figure}[h]
\centering
\epsfig{file=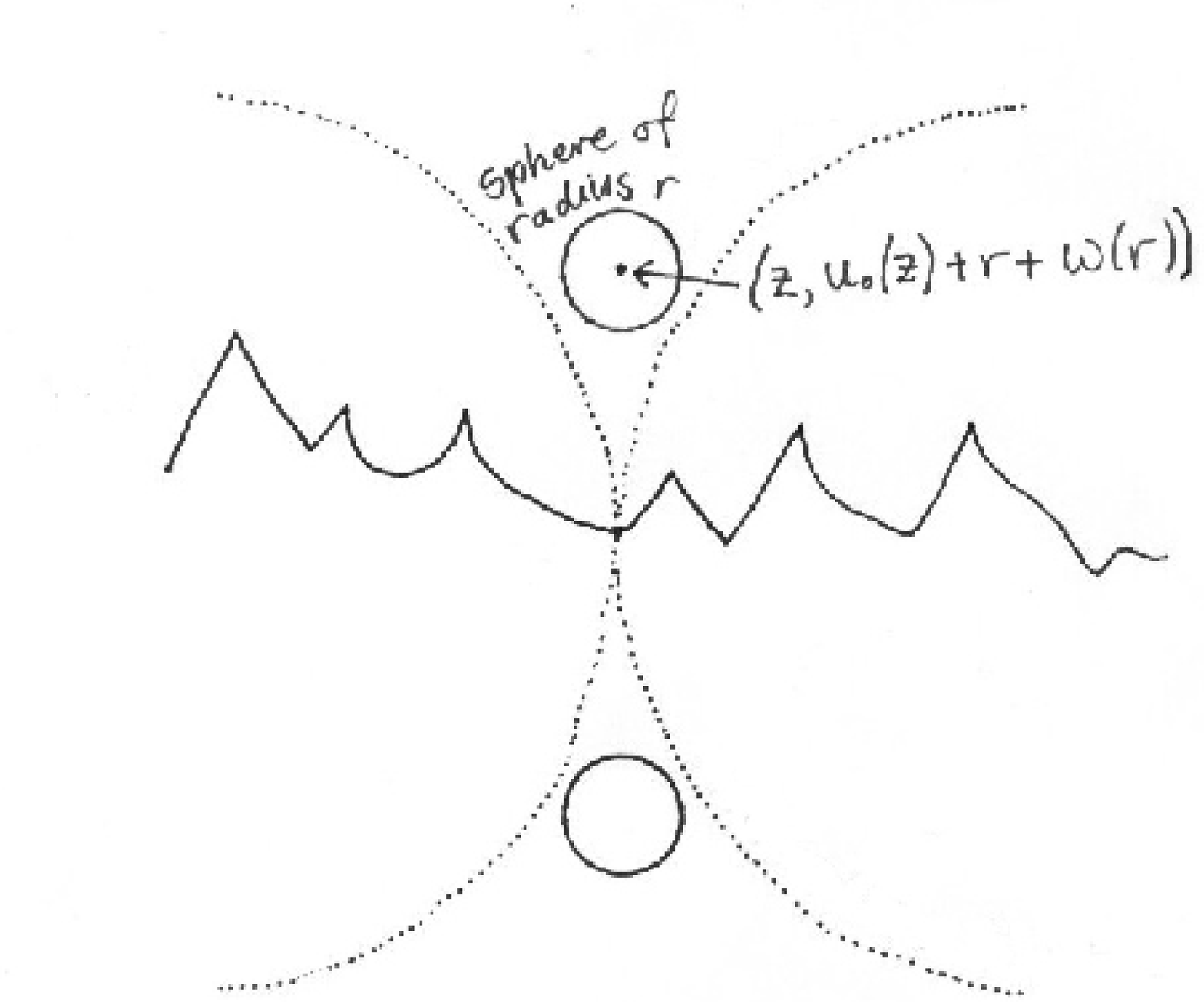, height=6cm}
\caption{Spheres above and below $w^\epsilon$ or $u_0$}
\end{figure}
\end{center}

Above and below this envelope, we can place two spheres of radius $r$ centred at \linebreak[3] $\left(z,u_0(z)\pm [r+\omega(r) ] \right)$.   At $t=0$, the spheres and the graph of  $w^\epsilon(\cdot,0)$ are completely disjoint.   The spheres are also completely disjoint from the graph of $u_0$.

The evolution of spheres under mean curvature flow is well-known --- the centre remains fixed and the radius shrinks from the initial radius $r(0)=r_0$, with 
$$r(t)=\sqrt{r_0^2-2nt}$$ 
until the sphere disappears at time 
$t={r_0^2}/{2n}$.

The parts of these spheres closest to the graph of $w^\epsilon$ --- the lower part of the upper sphere and the upper part of the lower sphere ---  are 
\begin{gather*}
S^+(x,t):=u_0(z)+ r_0 +\omega(r_0) -\sqrt{r(t)^2-|x-z|^2}, \\  
S^-(x,t):=u_0(z)- r_0 -\omega(r_0) +\sqrt{r(t)^2-|x-z|^2}. 
\end{gather*}

Suppose that one of these spheres and $w^\epsilon$ first touch at some time $t>0$.   From the Comparison Theorem \ref{comparison 1}, we know that at this time there must be an intersection occurring on the boundary of $\Omega$, say at $x_1\in\partial\Omega$ (this doesn't rule out other intersections occurring simultaneously on the interior).  

This intersection on the boundary is an extreme point of $S^{+/-}(\cdot,t)-w^\epsilon(\cdot,t)$ (either a minimum of $S^+-w^\epsilon$ or a maximum of $S^- - w^\epsilon$).  

Therefore the sign on the outward derivative of the intersecting sphere at this point is known --- either 
\begin{equation*}
D_\nu \left(S^+(x_1,t)-w^\epsilon(x_1,t)\right)\le0 \text{ and so } D_\nu S^+(x_1,t)\le0, 
\end{equation*}
or else 
\begin{equation*}
D_\nu \left(S^-(x_1,t)-w^\epsilon(x_1,t)\right)\ge0 \text{ and so }  D_\nu S^-(x_1,t)\ge 0,
\end{equation*}
 where we have used that $w^\epsilon$ satisfies the Neumann condition $D_\nu w^\epsilon=0$ on the boundary.

On the other hand, we can explicitly calculate the gradients of the \mbox{spheres ---} 
\begin{equation*}
D_i S^+(x,t)=\frac{2(x^i-z^i)}{\sqrt{r^2-|x-z|^2}},  \qquad D_i S^-(x,t)=\frac{-2(x^i-z^i)}{\sqrt{r^2-|x-z|^2}}.  
\end{equation*} 

The convexity of $\Omega$ means that for any $z\in\overline\Omega$ and $x_1$ on the boundary, $\nu\cdot(x_1-z)\ge 0$, and the inequality is strict if $z$ is in the interior of $\Omega$.     Hence the sign of the normal derivative is known --- 
\begin{equation*}
D_\nu S^+(x_1,t)\ge 0,
  \qquad D_\nu S^-(x_1,t)\le 0.  
\end{equation*} 

Between these two observations, it must be the case that the intersecting sphere has a flat normal gradient --- 
\begin{equation*}
D_\nu S^+(x_1,t)= 0
  \quad\text{ or }\quad D_\nu S^-(x_1,t)= 0,  
\end{equation*} 
and so $\nu\cdot(x_1-z)= 0$, which in turn implies that $z$ is on the boundary of $\Omega$, contradicting our original assumption that $z$ was an interior point.  It follows that such spheres, centred on interior points, never touch the graph of $w^\epsilon$ for the duration of the spheres' existence, until $t=r_0^2/2n$.  

In particular, above the point $z$ the surfaces move by no more than $r+\omega(r)$ in the time $t\in(0,{r}^2/2n)$.  We can choose $r>0$ so that $\delta=r+\omega(r)$, and a corresponding
$t_\delta={r^2/ {2n}}$, ensuring that 
\begin{equation*}
|w^\epsilon(z,t)-u_0(z)|\leq r+f(r)=\delta \text{ for } 0\leq t \leq t_\delta, \end{equation*}
where $t$ is dependent on $\delta$ and $\omega$ alone.

This estimate is
independent of $z$ and $\epsilon$, so 
\begin{equation*}
\lim_{\epsilon\rightarrow 0} |w^\epsilon(z,t)-u_0(z)|\leq \delta. \end{equation*}

We still need to estimate the second part of equation \eqref{significant eqn in final lemma of neumann chapter},  $|u^\epsilon(z,0)-u_0(z)|$.  However the convergence here is uniform (for $z\in\Omega$), so that 
\begin{equation*} \sup_{z\in \Omega}|u(z,t)-u_0(z)|\leq \delta + \sup_{z\in\Omega}\lim_{\epsilon\rightarrow 0}
|u^\epsilon(z,0)-u_0(z)|= \delta \end{equation*}
for $0< t\leq t_\delta$, and so $u$ is in $C([0,T);C(\Omega))$, with $u(\cdot,0)=u_0$.

This is also an estimate for the smoothness in time of the convergence; we can consider $\delta=\delta(t)$, with 
\begin{equation*}
\sup_{z\in\Omega} |u(z,t)-u_0(z)|\leq \delta(t)
\end{equation*}
by setting $t=r^2/2n$, so that $r=\sqrt{2nt}$ and thus $\delta(t)=\sqrt{2nt}+\omega(\sqrt{2nt})$.  The convergence in time that this gives is  \emph{at best} like $t^{1/2}$ --- which is in concordance with the result for initial data with a Lipschitz bound.  In that case,  $\omega(r)=|r|$ and the convergence is $C^{0+1/2}$ in time (Theorem 3.5 of \cite{hu:local}).   

 In the case that the initial data has a H\"older gradient bound, $\omega(r)=|r|^\alpha$, then the convergence to the initial data is as $t^{1/2}+t^{\alpha/2}\sim t^{\alpha/2}$. 
\halmos

\remark    While we have mined the rich theory arising from mean curvature flow to find this result, there are similar results for other equations of the type studied in Chapter \ref{higher dimensional chapter}, and we expect to be able to find similar existence results. 

In particular, one can find short-time existence results for anisotropic mean curvature flow with a zero Neumann boundary condition and continuous initial data on $\Omega\times[0,T]$, for convex $\Omega$. 
\newchapter{Existence of solutions to the Dirichlet problem for mean curvature flow }{Dirichlet problem}
\label{dirichlet chapter}

In this chapter we use the gradient estimate to establish existence of solutions to the Cauchy-Dirichlet problem with zero boundary data.
\begin{theorem}[Existence of solutions to the Dirichlet problem] \label{dirichlet}
Let $\Omega$ be a domain in $\mathbb R^n$, with $C^2$ boundary $\partial\Omega$ that has non-negative mean curvature.  If $u_0\in C^0(\overline\Omega)$ and $u_0=0$ on $\partial\Omega$, then the problem   
\begin{align}
\frac{\partial u}{\partial t}&=
\left(\delta_{ij}-\frac{D_iu D_ju}{ 1+|Du|^2}\right)D_{ij}u, \label{mcf in dirichlet chapter} \\
u(x,t)&=0\quad\text{for $x\in\partial\Omega$, $t>0$}, \notag 
\\
u(\cdot,0)&=u_0, \notag 
\end{align}
has a smooth solution for $t>0$ which converges uniformly to $u_0$ as $t\rightarrow 0$.
\end{theorem}

The existence of solutions to the mean curvature flow problem with prescribed boundary values was considered by Lieberman in \cite{li:first-bvp} (and by Huisken in \cite{hu:boundary}, where the long-time behaviour of solutions was also studied).    

Lieberman considered time-dependent boundary data $u_0\in H_{1+\alpha}\left({\cal P}(\Omega\times[0,T])\right)$, with a Lipschitz bound (in time) on $\partial\Omega\times [0,T]$.     The following theorem may be found as Theorems 12.10 and 12.18 of \cite{li:parabolic}.  

\begin{theorem}[Lieberman] \label{dirichletexistence} Let $\Omega\subset\bigR^n$ be a domain with $C^2$ boundary,  and  let $u_0$ be a function defined on the boundary, with $u_0\in H_{1+\alpha}\left({\cal P}(\Omega\times[0,T])\right)$ for some $0<\alpha<1$.  
Then if the mean curvature $H$ of $\partial\Omega$ is  non-negative, there exists a solution to \eqref{mcf in dirichlet chapter} with initial and boundary data
$$ u(x,t)=u_0(x,t) \text{ for } (x,t)\in{\cal P}(\Omega\times[0,T]).$$

Moreover, such a solution satisfies
\begin{equation*}
[Du]_{\beta}\le c,
\end{equation*}
where $c$ and $\beta$ depend on $|u_0|_{1+\alpha,\alpha/2}$.    
\end{theorem}

The proof of Theorem \ref{dirichlet} is very similar to that for the Neumann problem, Theorem \ref{neumann theorem}.  We will use the boundary-straightening change of coordinates described in Section \ref{boundary section}, and the corresponding $v_0(y,t):=u_0(\Psi^{-1}(y),t).$  

\begin{lemma} \label{approx-drich} There exists an approximating sequence $u^\epsilon_0\in C^\infty(\overline\Omega)$ with $u^\epsilon_0=0$ on $\partial\Omega$ and $u^\epsilon_0\rightarrow u_0$ in $C(\Omega)$.
\end{lemma}  
\proof  Let $B_R$ be a boundary centred ball, where $R$ is given by \eqref{definition of R}.  We will define a local approximation on $B_R$ and then put similar local approximations together to give a global one.

In Section \ref{boundary section} we defined $\tilde v_0$ to be a reflection across the boundary; this time, we let $\tilde v_0$ be the \emph{odd} reflection over the boundary 

\begin{equation} \tilde v_0(y):=-\frac{y^n}{|y^n|} v_0\left(\rho( y)\right)=
\begin{cases}u_0\left(\Psi^{-1}(\overline y,y^n)\right)& \text{if $y^n<0$,}\\
0 & \text{ if $y^n=0$,} \\
-u_0\left(\Psi^{-1}(\overline y,-y^n)\right)& \text{ if $y^n>0$ .}
\end{cases}
\end{equation}

Mollifying this in the standard way
\begin{equation}
v^\epsilon_0(y):=\eta_\epsilon*\tilde v_0 =\int_{z\in \Psi(B_R)} \eta_\epsilon(y-z)\tilde v_0(z) dz,
\end{equation}
we note that $v^\epsilon_0\rightarrow v_0$ uniformly on subsets of $\Psi(B_R)$, and we can check that if $y^n=0$, then $v^\epsilon(y)=0$. Returning to the original coordinates, set $u^\epsilon_1(x):=v^\epsilon_0(\Psi(x))$ on $\overline\Omega\cap B_R(x_1)$.

Now, cover $\partial\Omega$ by $N$ such boundary-centred balls $B_R(x_1),\dots,B_R(x_N)$ on which are defined approximations $u^\epsilon_1,\dots,u^\epsilon_N$.  Complete the cover of $\Omega$ by the set $\Omega^{R/2}:=\lbrace x\in\Omega: \dist(x,\partial\Omega)>R/2 \rbrace.$  On this interior set let $u^\epsilon_{N+1}$ be the usual mollification of $u_0$.  

If $\{\xi_i\}_{i=1,N+1}$ is a partition of unity with respect to these sets, then the sum $u_0^\epsilon=\sum_{i=1}^{N+1}\xi_iu_i^\epsilon$ converges uniformly to $u_0$ on $\overline\Omega$ and also has $u^\epsilon_0=0$ on $\partial\Omega$. 

We can restrict this to a subsequence $u^{\epsilon_i}_0$, where $|u_0-u^{\epsilon_i}_0|< 2^{-i} $
 (but retain the notation $u^\epsilon_0$ for the subsequence).  If we off-set each member of the subsequence, by replacing $u^{\epsilon_i}$ by $u^{\epsilon_i}+3(2^{-i})$, then this is a completely disjoint subsequence that still converges to $u_0$ as $i\rightarrow\infty$.  
\halmos

As the boundary values 
\begin{gather*}
u^\epsilon=u^\epsilon_0 \text{ on } \Omega\times\lbrace 0 \rbrace,\\
u^\epsilon=0 \text{ on }\partial\Omega\times[0,T], 
\end{gather*}
are in $H_{1+\alpha}$ on ${\cal P}(\Omega\times[0,T])$,   Theorem \ref{dirichletexistence} ensures that a solution with these boundary values for \eqref{mcf in dirichlet chapter} exists.   Denote these approximate solutions by $u^\epsilon$.

The gradient estimate derived for the Dirichlet problem in Theorem  \ref{mcf-dirichlet-theorem}
 implies that $$|Du^\epsilon(\cdot,t)|\le C_1\sqrt{t}\left(1+t\right)\exp(C_2/t):=L(t),$$ for constants $C_1$ and $C_2$ dependent only on $\osc u^\epsilon_0\le\osc u_0$.

\begin{lemma}
The approximate solutions $u^\epsilon$ have a H\"older gradient bound 
\begin{equation*}
|u^\epsilon(\cdot,t)|_{{C^\alpha}} \le C,
\end{equation*}
for $t>t_0>0$, where $C$ is dependent on $L(t_0)$,$|t_1-t_0|$, $n$, $|\partial\Omega|_{C^2}$, and $\diam\Omega$. 
\end{lemma}
\proof
If we revert to $v^\epsilon$ in the straightened-boundary coordinates satisfying \eqref{eqn4v} on $Q_R$, %
Theorem \ref{side-estimate} gives an oscillation bound for the gradient on the boundary-centred cylinder $Q_r(x_1,t_1)$ for $t_1>t_0+r^2$ and $r<R$ --- 
\begin{equation*}
\osc_{Q_r}Dv^\epsilon\le cr^\alpha,
\end{equation*}
where $c$ and $\alpha$ depend on $n$, $|a^{ij}|_{C^1}$, $|b|$, $\lambda_{a^{ij}}$ and $\Lambda_{a^{ij}}$, where $a^{ij}$ and $b$ are the coefficients of \eqref{eqn4v}, as in 
 \eqref{full write-out of a and b}.   These last four are in turn dependent on $L(t_0)$ and $|\partial\Omega|_{C^2}$.  
 
On the interior  we have the bounds on higher derivatives given by Lemma \ref{higher derivatives bounded on interior} 
\begin{equation}
|D^{m+2}u^\epsilon(x,t)|\le c\left(n,m,L(t_0)\right)
\left(\frac1{\dist(x,\partial\Omega)^2}+\frac1{t-t_0}\right)^{\frac{m+1}2} \label{interior bound in dirichlet chapter}. 
\end{equation} 
With $m=0$, this is a gradient bound for $Du$.

These can be linked together using Lemma \ref{linking lemma} to find a global H\"older bound
\begin{equation*}
|Du|_{\alpha';Q_r(x_1,t_1)}\le C,
\end{equation*}
where $\alpha'=\alpha/(\alpha+1)$ and $C$ additionally depends on $|t_0-t_1|$ and $\diam \Omega$. 
\halmos

Now, if we consider the approximate solutions to begin at some time $t_1>t_0$  with the initial data $u^\epsilon(\cdot,t_1)$, the uniform H\"older gradient bound on $u^\epsilon(\cdot,t_1)$ means that Theorem \ref{dirichletexistence} gives a H\"older gradient bound %
\begin{equation*}
|u^\epsilon|_{1+\beta;\Omega\times[t_1,T]}\le C,
\end{equation*}
for all $t_1>t_0>0$, where $C$ and $\beta$ depend on $L(t_0)$,$|t_1-t_0|$, $n$, $|\partial\Omega|_{C^2}$, and $\diam\Omega$.   
  
With these uniform estimates for positive times, there must be a convergence subsequence in $H_{1+\beta'}(\Omega\times[t_1,T])$ for some $\beta'<\beta$.   If we off-set  the initial data as mentioned in Lemma \ref{approx-drich},  the convergence of the subsequence implies  the convergence of the entire sequence to a limit $u$, defined on all $\Omega$ for $t>0$.  Finally, the interior bounds \eqref{interior bound in dirichlet chapter} mean that on the interior, $u$ is smooth.

We now need to show that $u(\cdot,t)\rightarrow u_0$ as $t\rightarrow 0$.

\begin{lemma} As $t\rightarrow 0$, 
$$\sup_{x\in\Omega}|u(x,t)-u_0(x)|\rightarrow 0.$$
Furthermore, $u$ has a modulus of continuity in time dependent only on that of $u_0$.  
\end{lemma}
\proof  As in the proof of Lemma \ref{convergence to initial time}, fix $z$ to be any point in the interior of $\Omega$ and set   
$$w^\epsilon(x,t):=u^\epsilon(x,t)-u^\epsilon(z,0)+u_0(z),$$ so that
$w^\epsilon$ is a solution to \eqref{mcf in dirichlet chapter} with 
$w^\epsilon(z,0)=u_0(z)$ and
$w^\epsilon(x,t)=-u^\epsilon(z,0)+u_0(z)$ on the boundary.  

Let $\omega$ be a modulus of continuity for $u_0$ and hence for $w^\epsilon(\cdot,0)$, so that 
\begin{equation*} -\omega(|z-x|)+u_0(z)\leq w^\epsilon(x,0) \leq \omega(|z-x|)+u_0(z). \end{equation*}

Above and below these two bounds, we can place two spheres of radius $r$ centred at $\left(z,u_0(z)\pm [r+\omega(r) ] \right)$.   At $t=0$, the spheres and the graph of  $w^\epsilon(\cdot,0)$ are completely disjoint.   The spheres are also completely disjoint from the graph of $u_0(\cdot)$.

As they evolve under mean curvature, the parts of these spheres closest to the graph of $w^\epsilon$ --- the lower part of the upper sphere and the upper part of the lower sphere ---  are 
\begin{gather*}
S^+(x,t):=u_0(z)+ r +\omega(r) -\sqrt{r^2-2nt-|x-z|^2}, \\  
S^-(x,t):=u_0(z)- {r} -\omega({r}) +\sqrt{{r}^2-2nt-|x-z|^2}. 
\end{gather*}

Initially $S^+-w^\epsilon$ is positive:  suppose that $t\in(0,{r}^2/2n)$ is the first time that $S^+-w^\epsilon$ decreases to zero.  The comparison principle (Theorem \ref{comparison 1}) means that this occurs at some point $x_1$ on the boundary.   On the boundary,  
\begin{equation*} \dfrac d {dt}\left(  S^+-w^\epsilon\right)\,=\, \frac n {\sqrt{{r}^2-2nt-|x_1-z|^2}} -0 \,>\, 0 
\end{equation*}
(since $w^\epsilon$ is constant on the boundary) and so this cannot be the first zero point; it follows that $S^+-w^\epsilon>0$ for the duration of the sphere's existence.   

The same argument shows that $S^- -w^\epsilon<0$.  

We can conclude that  the $w^\epsilon$ move at most by $r+\omega(r)$ in the time $t\in(0,r^2/2n)$.  

This estimate is
independent of $z$ and $\epsilon$, so this implies that $u^\epsilon(\cdot,t)\rightarrow u_0$ and that 
$u\in C\left(\Omega\times[0,T]\right)$.  Furthermore, if $u_0\in C^\alpha(\Omega)$, then $u\in H_\alpha\left(\Omega\times[0,T]\right)\cup C^\infty\left(\Omega\times(0,T]\right)$. \halmos

\newchapter{Gradient estimates found by counting intersections}{Gradient estimates found by counting intersections}

\label{zero counting chapter}

In the paper  \cite{ang:zeroset}, Angenent proved a series of results regarding the finiteness and non-proliferation of the zeroes of a parabolic equation in one space dimension.  

A {\emph{zero}} of $v(\cdot,t)$ is simply a point $x$ where $v(x,t)=0$.  A {\emph{multiple zero}}  is a point where both $v$ and $v_x$ vanish.  In contrast to earlier results, Angenent did not exclude multiple zeroes from the zero set, defining the zero set as 
$$z(t)=\lbrace\, x\in\bigR : v(x,t)=0\,\rbrace.$$
In the following,  $z(t)$ is often used as shorthand for the counting measure ${\cal H}^0\left(z(t)\right)$.

These zero-counting results have been influential in many different areas, and have been used for geometric flows by Angenent himself, in \cite{angenent:nodal}, \cite{angenent:curves-II}, \cite{angenent:formation}, and \cite{aag:rotation}, the last with Altshuler and Giga.  Many others working in the area have also used these results. 

Unlike approaches that depend more explicitly on the maximum principle, this technique seems limited to equations in one dimension.      The gradient estimates found do not depend on the initial gradient, but do depend explicitly on the height:   the smallest gradient estimates are found for when the height is largest.

This work originates in an idea of Ben Andrews; also, this approach to finding gradient estimates has been independently  used by Nagase and Tonegawa in the forthcoming paper \cite{yoshi:grad-estimates}.

\section{Counting zeroes}

\label{explaining angenent section}

The estimates in this chapter rely on Theorem D of Angenent's paper:

\begin{theorem}[Angenent] \label{angenent's theorem} Let $v:[x_0,x_1]\times [0,T]\rightarrow \bigR$ be a solution of  
\begin{equation*}
v_t=a(x,t)v_{xx}+b(x,t)v_x+c(x,t)v  \label{ang equation 2}
.\end{equation*}
 such that there are no zeroes on the boundaries
\begin{equation*} 
v(x_i,t)\not = 0, \qquad i=0,1.
\end{equation*}
Let $a$, $b$, $c$ satisfy %
\begin{align*}
&a {\text{ positive; }} \\
&a,\text{ }a^{-1},\text{ } a_t, \text{ }a_x, {\text{ and }} a_{xx} {\text{ bounded; }} \\
&b, \text{ }b_t {\text{ and }} b_x {\text{ bounded; }}\\
&c {\text{ bounded. }}
\end{align*}

Then if $v_t$, $v_x$ and $v_{xx}$ are continuous on $(x_0,x_1)\times[0,T]$, 
\begin{itemize}
\item for $t>0$, $z(t)$ is finite
\item if $\tilde x $ is a multiple zero of $v$ at $\tilde t$ then  for all $t_1<\tilde t <t_2$ we have $z(t_1)>z(t_2)$. %
\end{itemize}
\end{theorem} 

Consider a fully nonlinear equation on a domain $\Omega\times[0,T]$, where $\Omega$ is a connected subset of $\bigR$, 
\begin{equation}
u_t=F(u_{xx},u_x,u,x,t),  \label{ang nonlinear}
\end{equation}
where $F$ is parabolic, by which we mean that   
\begin{equation*}
\frac {\partial}{\partial r} F(r,p,q,x,t)>0
\end{equation*}
for all $(r,p,q,x,t)\in \bigR^3 \times \overline\Omega\times [0,T]$. 

Suppose that $u$ and $\varphi$ are smooth solutions of \eqref{ang nonlinear},
with
\begin{gather*}
|u|,|{u}_x|,|{u}_{xx}|,|{u}_{t}|
\le C_1 \intertext{ and } 
|\varphi|,|{\varphi}_x|,|{\varphi}_{xx}|, |{\varphi}_{t}|
\le C_1.
\end{gather*}

{\allowdisplaybreaks[1] Then we can form the difference $w:=u-\varphi$  satisfying the evolution equation
\begin{align}
w_t&=u_t-\varphi_t \notag \\*
&= F(u_{xx},u_x,u,x,t)-F(\varphi_{xx},\varphi_x,\varphi,x,t) \notag \\
&=\int_0^1 \frac d {ds} F(su_{xx}+(1-s)\varphi_{xx},su_x+(1-s)\varphi_x,su+(1-s)\varphi,x,t) ds \notag\\
&=\int_0^1 \frac\partial{\partial r}F(\dots)
ds \phantom{,} (u_{xx}-\varphi_{xx}) \notag 
+ \int_0^1 \frac\partial{\partial p}F(\dots)
ds \phantom{,}(u_x-\varphi_x) \notag 
\\*&\phantom{spacespacespace} 
+ \int_0^1 \frac\partial{\partial q}F(\dots)%
ds \phantom{,} (u-\varphi) \notag \\*
&= A(x,t)w_{xx} +B(x,t)w_x +C(x,t)w, \label{angenent intersection eqn}
\end{align}}
where the omitted argument of the derivatives of $F$, denoted by $(\dots)$, is always $\left(su_{xx}+(1-s)\varphi_{xx},su_x+(1-s)\varphi_x,su+(1-s)\varphi,x,t\right)$.  In the last line, 
\begin{gather}
A(x,t):=\int_0^1 \frac\partial{\partial r}F(su_{xx}+(1-s)\varphi_{xx},su_x+(1-s)\varphi_x,su+(1-s)\varphi,x,t)ds  \label{label one} \\
B(x,t):=\int_0^1 \frac\partial{\partial p}F(su_{xx}+(1-s)\varphi_{xx},su_x+(1-s)\varphi_x,su+(1-s)\varphi,x,t)ds  \label{label two} \\ \intertext{and}
C(x,t):= \int_0^1 \frac\partial{\partial q}F(su_{xx}+(1-s)\varphi_{xx},su_x+(1-s)\varphi_x,su+(1-s)\varphi,x,t)ds. \label{label seven}
\end{gather}

In order to use Angenent's theorem, we need to establish that:
\begin{itemize}
\item $A$, $A^{-1}$, $A_t$, $A_x$, and $A_{xx}$ are bounded,
\item $B$, $B_t$ and $ B_x $ are bounded 
\item and $C$  is bounded on $\Omega\times[0,T]$.
\end{itemize}

Let ${\cal{K}}=\left\lbrace\,(r,p,q,x,t)\in\bigR^3\times\Omega\times[0,T]:|r|,|p|,|q|\le  C_1\,\right\rbrace$.  
\label{page where cal S is defined}

If $\pd F r$ is continuous, then there are positive constants $\lambda_{\cal K}$ and $\Lambda_{\cal K}$ for which 
\begin{equation}
0<\lambda_{\cal K}\le \pd F r \le \Lambda_{\cal K} \text{ for all } (r,p,q,x,t)\in{\cal K}. \label{parabolicity for nonlinear F}
\end{equation}
Bounds on $A$ and $A^{-1}$ follow from this:
\begin{equation}
A%
 \le \Lambda_{\cal K}, \quad\quad\quad A^{-1} \le \lambda_{\cal K}^{-1}.  \label{upper and lower bounds on A} 
\end{equation}

$A_t$ is given by 
\begin{align*} 
A_t&=
\int_0^1  \pD {F} r  (\dots)
\left[su_{xxt}+(1-s)\varphi_{xxt}\right] 
+ \pdd F r p   
(\dots)
\left[su_{xt}+(1-s)\varphi_{xt}\right] 
\\ & \phantom{spacespacespace} 
+ \pdd F r q 
(\dots)
\left[su_{t}+(1-s)\varphi_{t}\right] 
+ \pdd F r t 
(\dots)
 ds,
\end{align*}
and so
\begin{equation}
|A_t|\le \bigg|\pd F r\bigg|_{C^1({\cal K})}\left(C_1+ |u_{xxt},u_{xt},
\varphi_{xxt},\varphi_{xt}|_{\Omega\times[0,T]}\right). \label{upper bound for A_t}  
\end{equation} %

Similarly, 
\begin{equation*} 
\begin{split} A_x=&\int_0^1  \pD F r  (\dots)
\left[su_{xxx}+(1-s)\varphi_{xxx}\right] 
+ \pdd F r p   
(\dots)
\left[su_{xx}+(1-s)\varphi_{xx}\right] 
\\ & \phantom{spacespace} 
+ \pdd F r q 
(\dots)
\left[su_{x}+(1-s)\varphi_{x}\right] 
+ \pdd F r x 
(\dots)
 ds, 
\end{split}
\end{equation*} so
\begin{equation}
|A_x|\le \bigg|\pd F r \bigg|_{C^1({\cal K})}(C_1+ |u_{xxx}
,\varphi_{xxx}|_{\Omega\times[0,T]}). \label{upper bound for A_x}  
\end{equation} 

$A_{xx}$ is given by 
\begin{equation*} 
\begin{split} 
A_{xx}=&\int_0^1  \pDDD F r  (\dots)
\left[su_{xxx}+(1-s)\varphi_{xxx}\right]^2 \\
& \phantom{space}+ 2\pd {^3F} {^2r\partial p }     
(\dots)
\left[su_{xxx}+(1-s)\varphi_{xxx}\right]\left[su_{xx}+(1-s)\varphi_{xx}\right] \\
& \phantom{spacespace }+ 2\pd {^3F} {^2r\partial q }    
(\dots)
\left[su_{xxx}+(1-s)\varphi_{xxx}\right]\left[su_{x}+(1-s)\varphi_{x}\right] \\
& \phantom{spacespacespace }+ 2\pd {^3F} {^2r \partial x } 
(\dots) 
\left[su_{xxx}+(1-s)\varphi_{xxx}\right] + R_1 \\ 
& \phantom{spacespacespacespace} + \pD F r  (\dots)\left[su_{xxxx}+(1-s)\varphi_{xxxx}\right]  + R_2\,ds,
\end{split}
\end{equation*}
where $R_1$ and $R_2$ are combinations of terms involving second and first derivatives of  $\pd F r$  respectively.  Consequently,
\begin{multline}  \label{upper bound for A_xx}
|A_{xx}|\le \bigg|\pd F r\bigg|_{C^2(\cal K)}\left(|u_{xxx},
\varphi_{xxx}|_{\Omega\times[0,T]}
+C_1\right)^2 \\
+\bigg|\pd F r\bigg|_{C^1(\cal K)}\left(|u_{xxxx}
,\varphi_{xxxx},
u_{xxx}
,\varphi_{xxx}|_{\Omega\times[0,T]}+C_1\right). 
\end{multline}

The bounds for $B$, its derivatives, and $C$ follow in a similar manner:
\begin{equation*} 
\begin{split} B_t=&\int_0^1  \pdd F p r  (\dots)
\left[su_{xxt}+(1-s)\varphi_{xxt}\right] %
+ \pD F  p   
(\dots)
\left[su_{xt}+(1-s)\varphi_{xt}\right] 
\\& \phantom{spacespace} 
+ \pdd F p q 
(\dots)
\left[su_{t}+(1-s)\varphi_{t}\right] 
+ \pdd F p t 
(\dots)
 ds, 
\end{split}
\end{equation*}
so that 
\begin{equation}
|B_t|\le \bigg| \pd F p \bigg|_{C^1(\cal K)}\left(C_1%
+ |u_{xxt},
\varphi_{xxt}|_{\Omega\times[0,T]}\right); \label{upper bound for B_t}  
\end{equation}
while  \begin{equation*} 
\begin{split} B_x=&\int_0^1  \pdd F p r  (\dots)
\left[su_{xxx}+(1-s)\varphi_{xxx}\right] 
+ \pD F  p   
(\dots)\left[su_{xx}+(1-s)\varphi_{xx}\right] 
\\&\phantom{spacespace}
+ \pdd F p q 
(\dots)
\left[su_{x}+(1-s)\varphi_{x}\right]
+ \pdd F p x 
(\dots)
 ds, \\
\end{split}
\end{equation*} so that
\begin{equation}
|B_x|\le \bigg| \pd F p \bigg|_{C^1(\cal K)}\left(C_1
+ |u_{xxx}
|_{\Omega\times[0,T]} + |
\varphi_{xxx}|_{\Omega\times[0,T]}\right).  \label{upper bound for B_x}  
\end{equation}
and finally 
\begin{align}
|C(x,t)|&= \left| \int_0^1 \frac\partial{\partial q}F(su_{xx}+(1-s)\varphi_{xx},su_x+(1-s)\varphi_x,su+(1-s)\varphi,x,t)ds \right| \notag \\
&\le \bigg| \pd F q \bigg|_{0;{\cal K}}.  \label{upper bound for C}
\end{align}
It is clear that we will be able to apply  Angenent's result to a smooth solution of a nonlinear parabolic equation $u_t=F$, when  $F$ satisfies the parabolicity condition \eqref{parabolicity for nonlinear F} and both $\pd F r$ and $F$ are $C^2$ on the bounded domain ${\cal K}$.  

These conditions are not optimal --- for example, in estimate \eqref{upper bound for C}  above, it is sufficient if $\pd F q$ is  $L^1$ along line segments in $\cal K$ --- however, they are enough to allow a theorem for \emph{intersections of two solutions} rather than \emph{zeroes} of one solution. 

\begin{theorem}[Intersection-counting theorem] 
 Let $u$ and $\varphi:[x_0,x_1]\times [0,T]\rightarrow \bigR$ be solutions of  
\begin{equation*}
u_t=F(u_{xx},u_x,u,x,t),
\end{equation*}
which do not intersect on the boundaries
\begin{equation*} 
u(x_i,t)\not=\varphi(x_i,t) , \quad i=0,1,\quad t\in[0,T].
\end{equation*}
If $u$ and $\varphi$ are $C^2$ on $(x_0,x_1)\times[0,T]$
$$ |u|_{C^2}, |\varphi|_{C^2} \le c_1,$$
and  if $F$ is parabolic
$$\pd { } r F(r,p,q,x,t)>0$$
and if both $F$ and $\pd F r$ are $C^2$ on
$${\cal{K}}=\left\lbrace\,(r,p,q,x,t)\in\bigR^3\times\Omega\times[0,T]:|r|,|p|,|q|\le  c_1\,\right\rbrace,$$
then
for $t>0$ the number of intersections of $u$ and $\varphi$ are finite; and if $\tilde x$ is an intersection of $u$ and $\varphi$ at $\tilde t$ then for all $t_1<\tilde t<t_2$, the number of intersections at $t_1$ is strictly less than the number of intersections at $t_2$.  
\end{theorem} 

\proof We apply Theorem \eqref{angenent's theorem} to the difference $w=u-\varphi$ which satisfies equation \eqref{angenent intersection eqn}.  \halmos

\section{Gradient estimates for equations in one space dimension}

In this chapter we seek interior estimates for bounded solutions of parabolic equations on  connected domains $\Omega\times[0,T]$.   When two functions intersect at a single point, then the gradient of the one that is smaller on the left of the intersection will dominate the gradient of the one that is smaller on the right. 

The main idea is that optimal regularity of $u(x,t)$ is found by comparison to the solution of the same parabolic equation with initial data $C \sigma$, where $\sigma$ is the maximal monotone graph 
\begin{equation} \label{step function f}
\sigma(x)=\begin{cases} +1, &x>0  \\ [-1,1], & x=0 \\-1, & x<0\end{cases}
\end{equation}  \nomenclature[sig]{$\sigma $      }{  the step ``function''\refpage}
\label{page for step function f}which we will refer to as the {\emph{step ``function''}}.

\begin{center}
\begin{figure}[ht]
\centering
\epsfig{file=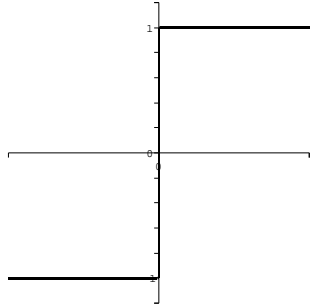, width=5cm}
\caption{The step function $\sigma$}
\end{figure}
\end{center}

The method can be broken into the following steps: 
\begin{itemize}
\item Creation of family of barriers $\lbrace\varphi^{\epsilon,s}\rbrace $ with $\varphi^{\epsilon,s}(x,t)$ approaching $C \sigma(x-s)$ as $t,\epsilon\rightarrow 0$ 
\item Show that for all $(x,t)$ in a subdomain of $\Omega\times[0,T]$, and for all $k\in[-M,M]$, we can find an $s$ such that $\varphi^{\epsilon,s}(x,t)=k$  
\item Show that $|\varphi^{\epsilon,s}|>M$ at the boundaries of $\Omega$ 
\item Then use the Angenent result to count the intersections of $u$ and $\varphi^{\epsilon,s}$.  For small enough $\epsilon$, there will be only one intersection, and so a gradient bound will follow.
\end{itemize}

\suppressfloats[t]%

\begin{figure}
\begin{center}  
\epsfig{file=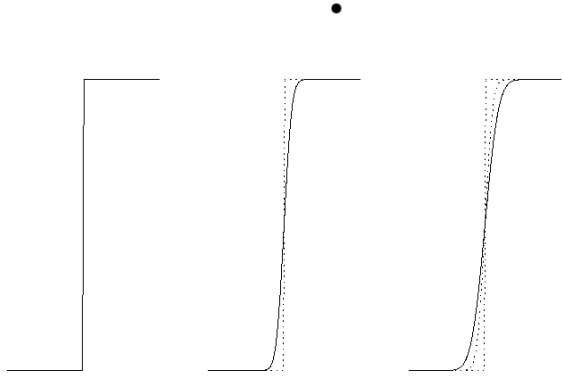, height=5cm}
\caption{Evolution of a step function under curve shortening flow}
\end{center}
\end{figure}

We begin by looking at a simple estimate for entire solutions on $\bigR$, then find more specific estimates, firstly for the  
 heat equation
\begin{equation}
u_t=\frac1 {4c} {u_{xx}}, \label{heat} 
\end{equation}
and then for a nonlinear  problem.

The following theorem says that if a solution with the step function as initial condition exists, then it will serve as a barrier for other solutions.

\begin{theorem}
Consider the parabolic equation 
\begin{equation}
u_t=F(u_{xx},u_x,u), \label{pde 1}
\end{equation}
where $F$ satisfies \eqref{parabolicity for nonlinear F}.
Let $u$ be a   solution of \eqref{pde 1} on $\bigR\times(0,T]$ that has a bound
\begin{equation*}
|u|\le M,
\end{equation*}
and has a uniform gradient bound at $t=0$.

Suppose there exists a solution to \eqref{pde 1} on $\bigR\times(0,T]$ which is smooth for $t>0$, and has initial condition
\begin{equation*}
\varphi(x,0)=(M+1)  \sigma(x)
\end{equation*}
and boundary condition
\begin{equation*}
\lim_{|x|\rightarrow\infty}\varphi(x,t)=(M+1)\sigma(x).
\end{equation*}

Then there is a gradient estimate 
\begin{equation*}
u_x(x,t)\le \varphi_x(z,t)
\end{equation*}
where $z$ is chosen so that $\varphi(z,t)=u(x,t)$.
\end{theorem}
\proof 
Let a family of barriers indexed by $(z,\tau)$ be given by $\varphi^{z,\tau}(x,t):=\varphi(x-z,t+\tau)$ for all $z\in\bigR$ and $\tau>0$.  Each of these  satisfies \eqref{pde 1}, and is smooth on $\bigR\times[0,T]$. 
     
As $u(\cdot,0)$ has a uniform gradient bound, there exists a $\tau'>0$ such that not only do $u(\cdot,0)$ and $\varphi^{z,\tau'}(\cdot,0)$ intersect only once, but also,   $u(\cdot,0)$ and $\varphi^{z,\tau}(\cdot,0)$ intersect only once for all $\tau\in(0,\tau']$.

Let $(x_1,t_1)$ be fixed. For each $\tau\le\tau'$, there exists $z$ such that $\varphi^{z,\tau}(x_1,t_1)=u(x_1,t_1)$.   

Now, apply Angenent's theorem to $w=u-\varphi^{z,\tau}$ on some region $[-R,R]\times[0,t_1]$ containing $x_1$ and which is sufficiently large enough that for all $t\in[0,t_1]$, $\varphi^{z,\tau}(R,t)\ge M$ and $\varphi^{z,\tau}(-R,t)\le -M$.  The last conditions ensure that $w$ has no zeroes on the boundary.  

As $w$ has only one zero at $t=0$, it has no more than one zero for all $t$; as $w$ is positive at $x=-R$ and negative at $x=-R$, it has exactly one zero for all $t$.  In particular the zero at $(x_1,t_1)$ is the only zero, and  
\begin{gather*}
\text{ for }x> x_1,\quad  \varphi^{z,\tau}(x, t_1)>u(x, t_1), \\
\text{ for }x<x_1, \quad \varphi^{z,\tau}(x,t_1)<u(x,t_1), 
\end{gather*}
from which we find a gradient estimate:
\begin{equation*}u_x(x_1,t_1)\le {\varphi^{z,\tau}}_x(x_1,t_1).
\end{equation*}

This holds for all $\tau\in(0,\tau_1)$ and so letting $\tau\rightarrow 0$ gives the result. \halmos

The following theorem describes an explicit barrier in the case of the heat equation.

\begin{theorem}[Gradient estimate for the heat equation] 
\label{theorem for heat equation}
Let $\Omega=[x_0,x_1]$ and $u:\Omega\times[0,T]\rightarrow\bigR$ be a smooth solution to the heat equation \eqref{heat}  with a height bound $|u|<M$ and Lipschitz bound $\Lip u(\cdot,0)<\infty$.  
 
Then for $t>0$,
\begin{equation*} 
u_x(x,t)\le  
2N \sqrt{\frac c{\pi t}}\exp\left(-\inverf{\frac u N}^2\right),
\end{equation*}
where 
\begin{equation*}
N=M\left[{\erf{\frac{\sqrt  c \phantom{i} \dist(x,\partial\Omega)}{2\sqrt t}}}\right]^{-1}.
\end{equation*}

\end{theorem}

This leads to an estimate for an entire solution.

\begin{corollary}
Let $u:\bigR\times[0,T]\rightarrow\bigR$ be a smooth solution of the heat equation \eqref{heat} with $|u|<M$.  

Then for $t> 0$, 
\begin{equation*}
u_x(x,t)\le {2M}
\sqrt{\frac c{\pi t}}\exp\left(-\inverf{\frac u M}^2\right).
\end{equation*}
\end{corollary}

\proof  Theorem \ref{theorem for heat equation} applies on any interval $[-R,R]$;  let $R\rightarrow\infty$ to find the result.   \halmos

\noindent{\textbf{Proof of Theorem \ref{theorem for heat equation}:}}
Without loss of generality, let $\Omega=[-\limbda,\limbda]$.

For some $\epsilon\in(0,\epsilon_0)$, $|s|<d$ and $N\ge M$ to be chosen later,  we can define the barrier
\begin{equation*}\varphi^{\epsilon,s}(x,t):=N\erf{(x-s)\sqrt{\frac c {t+\epsilon} }}
=\frac{2N}{\sqrt\pi}\int^{(x-s)\sqrt{\frac c{ t+\epsilon} }}_0 e^{-y^2}dy,
\end{equation*}
which  satisfies the heat equation \eqref{heat} on $\Omega\times[0,T]$. As $t+\epsilon\rightarrow 0$, $\varphi^{\epsilon,s}(x,t)\rightarrow N\sigma(x-s)$, where $\sigma$ is the step function \eqref{step function f}.%

Choose $\epsilon_0$ small enough, such that for any $|s|<d$ and $\epsilon<\epsilon_0$ there is only one intersection of $u(\cdot,0)$ and $\varphi^{\epsilon,s}(\cdot,0)$.  This is possible as $u(\cdot,0)$ has a uniform gradient bound.

Now, let $(\tilde x,\tilde t)\in \Omega\times(0,T]$ be fixed, and consider $\tilde u:=u(\tilde x,\tilde t).$  The bound on $u$ implies that $\tilde u\in(-M,M)$.  

If we choose $s=\tilde x-\sqrt{\dfrac {\tilde t+\epsilon} c}\inverf{\tilde u/N}$, then $\varphi^{\epsilon,s}(\tilde x,\tilde t)=\tilde u$.  
Also choose 
\begin{equation}N=M\left[\erf{\frac{\sqrt c(\limbda-|\tilde x|)}{2\sqrt{\tilde t+\epsilon}}}\right]^{-1}.\label{definition of n}\end{equation}
With these choices, we can check that
\begin{align*}
|s|&\le |\tilde x|+ \sqrt{\frac{\tilde t+\epsilon}c}\left|\inverf{ \tilde u / N}\right|\\
&=|\tilde x|+ \sqrt{\frac{\tilde t+\epsilon}c}\left|\inverf{\frac {\tilde u} M \erf{\frac{\sqrt c (\limbda-|\tilde x|)}{2\sqrt{\tilde t+\epsilon}}}}\right|\\
&<|\tilde x|+\sqrt{\frac{\tilde t+\epsilon}c}\left|{\frac{\sqrt c (\limbda-|\tilde x|)}{2\sqrt{\tilde t+\epsilon}}}\right|\\
&=|\tilde x|+\frac{\limbda-|\tilde x|}{2}\\
&<\limbda;
\end{align*}
and that on the boundaries $|\varphi^{\epsilon,s}|\ge M$ whenever $t<\tilde t$, since
{\allowdisplaybreaks[1]{\begin{align*}
|\varphi^{\epsilon,s}(\pm \limbda,t)|&=N\erf{|\pm \limbda -s|\sqrt{\frac c {t+\epsilon}}} \\*
&=N\erf{\left|\pm \limbda -\tilde x+\sqrt{\dfrac {\tilde t+\epsilon} c}\inverf{\tilde \tilde u/N}\right|\sqrt{\frac c {t+\epsilon}}} \\*
&\ge N\erf{\left|\pm \limbda -\tilde x\right|\sqrt{\frac c {t+\epsilon}}-\sqrt{\frac{\tilde t+\epsilon}{t+\epsilon}}\left|\inverf{\tilde u/N}\right|} \\*
&\ge N\erf{\left(\limbda -|\tilde x|\right)\sqrt{\frac c {t+\epsilon}}-\sqrt{\frac{\tilde t+\epsilon}{t+\epsilon}}\inverf{M/N}} \\
\intertext{ and if we use \eqref{definition of n} for $(\limbda-|\tilde x|)$, then this is}
&=  N\erf{2\sqrt{\frac{\tilde t+\epsilon}{t+\epsilon}}\inverf{M/N}-\sqrt{\frac{\tilde t+\epsilon}{t+\epsilon}}\inverf{M/N}} \\*
&\ge M.
\end{align*}}}

We can now apply the intersection counting Theorem \ref{angenent intersection eqn} with $w=u-\varphi^{\epsilon,s}$, with the previous calculation ensuring that there are no intersections on the boundary, and with the coefficients in equation \eqref{angenent intersection eqn} given by 
$A=1$, $B=0$ and $C=0$.
Since there is only one intersection at the initial time, there is never more than one intersection at later times (in particular, for our given $(\tilde x,\tilde t)$, there is no other intersection at time $\tilde t$ than the one at $\tilde x$).

It follows that  
\begin{gather*}
\text{ for }y>\tilde x,\quad  \varphi^{\epsilon,s}(y,\tilde t)>u(y,\tilde t), \\
\text{ for }y<\tilde x, \quad \varphi^{\epsilon,s}(y,\tilde t)<u(y,\tilde t), 
\end{gather*}
from which we find a gradient estimate:
\begin{equation*}u_x(\tilde x,\tilde t)\le {\varphi^{\epsilon,s}}_x(\tilde x,\tilde t).
\end{equation*}

This holds for any smaller $\epsilon>0$, so letting $\epsilon\rightarrow 0$ gives the final result.
\halmos

This method applies to all parabolic operators for which we can find solutions that have the step function as the initial condition.  %

When a quasilinear parabolic equation satisfies the conditions of Section
\ref{step function existence section},  we  can use the solutions with stepped initial data, whose existence was shown in that chapter.  %

Let $a>0$ be in $H_{\alpha}(\cal K)$ for all bounded ${\cal K}\subseteq {\bigR}\times {\bigR}\times \Omega\times [0,T]$, and some $\alpha\in(0,1)$.  This implies that for every such $\cal K$ we can find 
positive  $\lambda_{\cal K}$ and  $\Lambda_{\cal K}$ such that 
$ \lambda_{\cal K}\le a(p,q,x,t)\le\Lambda_{\cal K}, \text{ when } (p,q,x,t)\in{\cal K}. 
$

Since we will be looking at bounded solutions in $\Omega$, the bound on the gradient is the pertinent bound on $\cal K$; we highlight this by writing: 
\begin{equation}
\lambda(K) \le a(p,q,x,t)\le\Lambda(K), \text{ when } |p|\le K,
 \label{definition of lambdas 2}  \end{equation}
where we assume that $ |q|\le M$, $x\in\Omega$  and $t\in[0,T]$.

Also, suppose that there are positive constants $A$ and $P$ such that 
\begin{equation}
a(p,q,x,t)p^2\ge A >0,{\text{ for }}|p|\ge P.  \label{degeneracy control 2}
\end{equation} 

\begin{theorem} \label{second zerocounting theorem on omega}
Let $u:\Omega\times[0,T]\rightarrow\bigR$ be a smooth solution to 
\begin{equation}
u_t=a(u_x,u,x,t)u_{xx}, \label{yet another parabolic equation}
\end{equation}
where $a$ satisfies  \eqref{definition of lambdas 2} and \eqref{degeneracy control 2}.

Let $u$ be bounded, $|u(x,t)|<M$.

Let $\varphi^s$ solve \eqref{yet another parabolic equation} on $\bigR\times(0,T]$, with $\varphi^s(\cdot,t)\rightarrow 2M\sigma(x-s)$ as $t\rightarrow 0$, where $\sigma$ is the step function \eqref{step function f}, and where $s$ is chosen so that $u(x,t)=\varphi^s(x,t)$.

If $t\le c
\phantom{i}{ \dist(x,\partial\Omega)}^2 /\Lambda\left({cM}/{\dist(x,\partial\Omega)}\right)$, 
where $\Lambda(\cdot)$ is given by \eqref{definition of lambdas 2}, then 
\begin{equation*}
u_x(x,t)\le{\varphi^s_x}(x,t). 
\end{equation*}

That is, the gradient of $u$ is bounded by the gradient of the barriers,  at the same height.

\end{theorem}

{\smallskip{\noindent\textbf{Remark 1:}}\quad}  We can replace $\sigma(x-s)$ by $\sigma(s-x)$ here, in which case we find that 
\begin{equation*}
u_x(x,t)\ge{\varphi^s_x}(x,t). 
\end{equation*}

{\smallskip{\noindent\textbf{Remark 2:}}\quad} If $a$ has polynomial growth, so that $\Lambda(K)\le c(1+K^q)$ for some $q\ge0$, then the interior region on which we can find bounds of this form is given by $t\le c\phantom{i}\dist(x,\partial\Omega)^{2+q}M^{-q/2}$, for some constant $c$.

{\smallskip{\noindent\textbf{Remark 3:}}\quad} If $\Omega=\bigR$ in Theorem \ref{second zerocounting theorem on omega}, then the gradient bound applies for all $t\in(0,T]$.

\smallskip
\noindent{\textbf{Proof of Theorem \ref{second zerocounting theorem on omega}:}}

We initially assume that $\Omega=[-d,d]$.  We will derive a gradient bound at a single point $(0,t)$, and then generalize it to interior points on a general domain.

As $u$ is smooth there are bounds on the first derivative and on higher derivatives
\begin{equation*}
|u_x|\le c_1, \quad\quad
|u_{xx},u_t,u_{xt},u_{xxx},u_{xxxx},u_{xxt}|\le c_2.
\end{equation*}

For $0<\epsilon\le\epsilon_0<<\limbda/4$, let $\varphi^\epsilon$ be the standard mollification of $\varphi$ 
\begin{equation*}
\varphi^\epsilon:=2M \eta_\epsilon * \sigma,
\end{equation*}
where $\sigma$ is the step function \eqref{step function f}.  

Choose $\epsilon_0$ small enough so that for all $\epsilon\le \epsilon_0$, $\varphi^\epsilon$ satisfies a gradient estimate from below:
\begin{equation}
 |{\varphi^\epsilon}_x(x)|\ge c_1 \text{ whenever } |\varphi^\epsilon(x)|\le M \label{gradient bound from below 2}.
\end{equation}

Now, for fixed $\epsilon$, define a family $\lbrace \varphi^{\epsilon,s}\rbrace_{|s|\le \limbda/2}$ of barriers, each of which solves \eqref{yet another parabolic equation} on $\bigR\times[0,\tau]$ (for some $\tau$ to be decided later) with initial condition \label{page what has varphi family on it}
\begin{equation*}
\varphi^{\epsilon,s}(x,0)=\varphi^\epsilon(x-s).
\end{equation*}  
The existence of such solutions follows from Corollary \ref{existence of entire step functions}, which applies as $a$ satisfies \eqref{degeneracy control 2}.

Standard results give %
\begin{gather*}
|\varphi^{\epsilon,s}_x|\le c_3(\epsilon) \\
|\varphi^{\epsilon,s}_{xx},\varphi^{\epsilon,s}_{xxx},\varphi^{\epsilon,s}_{xxxx},\varphi^{\epsilon,s}_{t},\varphi^{\epsilon,s}_{tx},\varphi^{\epsilon,s}_{xxt}|\le c_4({\epsilon}). \label{gradient bounds for varphi 2}
\end{gather*}

To avoid intersections of $u$ and the barriers occurring on the boundary, we need to show that $|\varphi^{\epsilon,s}|\ge M$ when $x\in\lbrace-\limbda, \limbda \rbrace$.  

Each barrier in the family is initially bounded above by a step function %
\begin{equation*}
\varphi^{\epsilon,s}(x,0)\le 2M\sigma(x-s+\epsilon) 
\end{equation*} 
and so Corollary \ref{step displacement} provides an estimate for $x<s-\epsilon$
\begin{equation*}
\varphi^{\epsilon,s}(x,t)\le  
\frac{8M}{|x-s+\epsilon|}\sqrt{\frac{\Lambda t}\pi} -2M,
\end{equation*}
where $\Lambda=\Lambda\left({4M}/{|x-s+\epsilon|}\right)$.  

In particular, at $x=-\limbda$, 
\begin{equation*}
\varphi^{\epsilon,s}(-\limbda,t)\le  
\frac{8M}{|-\limbda-s+\epsilon|}\sqrt{\frac{\Lambda t}\pi} -2M. \\
\end{equation*}
As $|s|<d/2$ and $\epsilon<<d/4$,   
$$\frac1{|-\limbda-s+\epsilon|}< \frac4\limbda\quad\text{ and }\quad\Lambda\left(\frac{4M}{|-\limbda-s+\epsilon|}\right)\le \Lambda\left(\frac{16M}{\limbda}\right). $$
If we choose $ \tau=\dfrac{\limbda^2\pi}{ 32^2{{\Lambda}}}$ where $\Lambda=\Lambda(16M/d)$, then 
\begin{align*}
\varphi^{\epsilon,s}(-\limbda,t)&\le  
\frac{32M}{\limbda}\sqrt{\frac{\Lambda t}\pi} -2M \\
&\le -M,
\end{align*}
whenever $t\le \tau$.

A similar calculation for the other boundary point $x=\limbda$ gives that 
\begin{equation*}
\varphi^{\epsilon,s}(\limbda,t)\ge
M
\end{equation*}
when $t\le \tau$.

Let $\epsilon$ be fixed.

For each $s\in[-\limbda/2,\limbda/2]$, we can define $w:=u-\varphi^{\epsilon,s}$  satisfying
\begin{equation*}
w_t=Aw_{xx}+Bw_x + Cw
\end{equation*}
on $[-\limbda,\limbda]\times[0,\tau]$.  
Here, $A$ is given by setting $\pd F r (r,p,q,x,t)= a(p,q,x,t)$ in \eqref{label one}, $B$ by setting $\pd F p= r \pd { }{p}a(p,q,x,t)$ in \eqref{label two}, and $C$ by setting $\pd F q =r \pd { }q a(p,q,x,t)$ in \eqref{label seven}.

Let ${\cal K}=\lbrace\,(p,q,x,t): |p|\le c_1+c_3(\epsilon),|q|\le M,x\in\Omega,t\in[0,\tau] \,\rbrace.$

 $A,A^{-1},A_t,A_x,A_{xx},B,B_t,B_x$, and $C$ are bounded by constants dependent on $c_1$, $c_2$, $c_3(\epsilon)$, $c_4(\epsilon)$, $\lambda\left(c_1+c_3(\epsilon)\right)$,  $\Lambda\left(c_1+c_3(\epsilon)\right)$ and $|a|_{C^2({\cal K})}$,  as in \eqref{upper and lower bounds on A}--\eqref{upper bound for C}.

Since $|u|< M$, and at the boundary $|\varphi^{\epsilon,s}|\ge M$, $w$ is {\emph{never}} zero on the boundary.

In particular, $w(-\limbda,t)=u(-\limbda,t)-\varphi^{\epsilon,s}(-\limbda,t)>0$ and $w(\limbda,t)=u(\limbda,t)-\varphi^{\epsilon,s}(\limbda,t)<0$, so there is always {\emph{at least}} one zero of $w$.  The lower gradient bound \eqref{gradient bound from below 2} implies that there is {\emph{at most}} one zero of $w$ at the initial time.  

Then the intersection counting theorem (Theorem \ref{angenent intersection eqn})  implies there is {\emph{exactly}} one zero of $w$ for {\emph{all}} $t\le \tau$.

In the following lemma we show that given $(x,t)$ --- or more specifically, $(0,t)$ --- we can find $s$ such that $(0,t)$ is a zero of $w=u-\varphi^{\epsilon,s}$ .    We will then return to the proof of Theorem \ref{second zerocounting theorem on omega}.

\begin{center}
\begin{figure}
\centering
        \epsfig{file=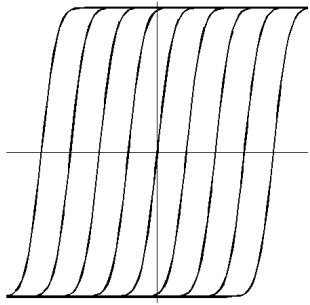, height=5cm}
\caption{A family of barriers}
\end{figure}
\end{center}

\begin{lemma}  Let $\epsilon\le\epsilon_0$ and $t\le \tau$ be fixed.  For each $k\in[-M,M]$, there exists an $s\in[-\limbda/2,\limbda/2] $ such that
\begin{equation*}
\varphi^{\epsilon,s}(0,t)=k.  
\end{equation*}
\end{lemma}
\proof
Firstly, we check that  $\varphi^{\epsilon,\limbda/2}(0,t)\le-M$ and  $\varphi^{\epsilon,-\limbda/2}(0,t)\ge M$.  Using 
Corollary \ref{step displacement}, 
\begin{equation*}
\varphi^{\epsilon,\limbda/2}(x,t)\le \frac{8M}{|x-\limbda/2+\epsilon|}\sqrt{\frac{\Lambda t}\pi}-2M
\end{equation*}
for $x<\limbda/2-\epsilon$, where $\Lambda=\Lambda(4M/|x-\limbda/2-\epsilon|)$. 
Since $|\limbda/2-\epsilon|^{-1}\le 4/\limbda$, at $x=0$ we have
\begin{align*}
\varphi^{\epsilon,\limbda/2}(0,t)&\le \frac{32 M}{\limbda}\sqrt{\frac{\Lambda t}\pi}-2M \\
&\le -M.
\end{align*}

It can similarly be shown that
$\varphi^{\epsilon,-\limbda/2}(0,t)\ge M$.  

As  $\varphi^{\epsilon,s}(\cdot,0)$ is continuous in $s$, and $a$ is a continuous operator, $\varphi^{\epsilon,s}(0,t)$ is also continuous in $s$.  In particular, $\lbrace\varphi^{\epsilon,s}(0,t)\rbrace_{|s|\le d/2}$ is onto $[-M,M]$.
\halmos

\smallskip
\noindent{\textbf{Continuation of the proof of Theorem \ref{second zerocounting theorem on omega}:}}\quad Let $t\le\tau$ be given.  From the previous lemma, %
there exists $s$ such that $\varphi^{\epsilon,s}(0,t)=u(0,t)$.  This is the only intersection point of $u$ and $\varphi^{\epsilon,s}$, and so
\begin{gather*}
\text{ for }y>0,\quad  \varphi^{\epsilon,s}(y,t)>u(y,t), \\
\text{ for }y<0, \quad \varphi^{\epsilon,s}(y,t)<u(y,t), 
\end{gather*}
from which we find the gradient estimate:
\begin{align*}u_x(0,t)&=\lim_{y\rightarrow 0}\frac{u(y,t)-u(0,t)}{y} \\
&\le \lim_{y\rightarrow 0}\frac{\varphi^{\epsilon,s}(y,t)-\varphi^{\epsilon,s}(0,t)}{y} \\
&= {\varphi^{\epsilon,s}}_x(0,t).
\end{align*}
This estimate holds for all $\epsilon\in(0,\epsilon_0]$.

If we let $\epsilon\rightarrow 0$, we firstly have that 
\begin{equation*}
\varphi^{\epsilon,s}\rightarrow \varphi^s, 
\end{equation*} where $\varphi^s$ is the solution to \eqref{yet another parabolic equation} with discontinuous initial data $2M\sigma(x-s)$; 
and secondly that
for all $t\le\tau$,
\begin{equation*}
u_x(0,t) \le\varphi^{s}_x(0,t), 
\end{equation*}
where $s$ is chosen so that $u(0,t)=\varphi^s(0,t)$.

Now we turn to the general domain with $\Omega=[x_0,x_1]$.   Given $x\in\Omega$, set $d:=\dist(x,\partial\Omega)$.  
If $t\le \dfrac {d^2\pi}{{32^2 \Lambda(24M/d)}}$ then we can repeat the same calculation on the small domain $[x-\limbda,x+\limbda]$ to find the given result.
\halmos


\newchapter{Estimates for isotropic and anisotropic mean curvature flow}{Estimates for isotropic and anisotropic mean curvature flows    }
 \label{newwork chapter}
\section{A gradient estimate for mean curvature flow} \label{section new work} This section
follows the style established in papers such as \cite{eh:mean} and \cite{eh:interior}, in particular the  local gradient estimates of Section 2 of the latter paper.    Ecker and Huisken 
consider the evolution of a hypersurface by mean curvature
 \begin{equation} \frac{d }{dt}
\mathbf{F}(p,t)=\mathbf{H}(p,t), \qquad p\in M, \label{mean curvature flow} \end{equation}
where $\mathbf{F}:M^n\times[0,T]\rightarrow \bigR^{n+1}$ is the immersion of the manifold $M^n$ at each time $t$ and
$\mathbf{H}$ is the mean curvature vector.

$M$ can be written as a graph when a fixed vector $\omega\in\bigR^{n+1}$ can be found so that for a choice of
unit normal $\nu$, \begin{equation*} \langle \nu,\omega \rangle >0 \end{equation*} everywhere.  Equivalently, $\langle \nu,\omega \rangle ^{-1}$ is bounded above.  The existence of an upper bound for this quantity and its analogue for anisotropic mean curvature flow is the subject of this chapter. 

Given the image $\mathbf{F}(p,t)$ of a point $p\in M$, its coordinate vector is $\mathbf{x}(p,t)$.  The
{\emph{height}} of $M$ above the hyperplane defined by $\omega$ is denoted by \begin{equation*}
u=\langle \mathbf{x}, \omega \rangle.  \end{equation*}

The \emph{gradient function} is given by \begin{equation*} v=\langle \nu, \omega \rangle^{-1}=\sqrt{1+|Du|^2}  \end{equation*}
We recall the evolution equations from \cite{eh:mean}: \begin{lemma}\label{evolution equation lemma} If
$M_t$ satisfies \eqref{mean curvature flow} then \begin{align*} \left(\frac d {dt} -
\Delta\right)&|\mathbf{x}|^2=-2n, \\ \left(\frac d {dt} - \Delta\right)&u=0,\\ \left(\frac d {dt} -
\Delta\right)&v=-|A|^2v-2v^{-1}|\nabla v|^2, \end{align*} where $\Delta$ is the Laplace-Beltrami operator
on $M_t$ and $A=\lbrace h_{ij} \rbrace$ is the second fundamental form. \end{lemma}

We derive a gradient estimate for periodic entire graphs,  followed by an interior estimate.   Estimates of this type have also been recently found by Colding and Minicozzi \cite{colding:2003} in the isotropic case using similar techniques, although without the explicit dependence on the height of the graph that the following estimates display.    

\begin{theorem}[Estimate for periodic mean curvature flow] Let $\mathbf{F}$ be a smooth, entire solution to mean curvature flow \eqref{mean curvature flow} which is   a periodic graph over a hyperplane, in that $u(y,t)=u(y+L,t)$ for a fixed point $L$ in the hyperplane,    and has a height bound $|u|<M$.  Then 
\begin{equation*}
v\le t^{1/2}\exp\left(\frac{c(|u|-2M)^2}{4t}\right)
\end{equation*}
for $0< t\le T'$, where $c$ and $T'$ depend on $M$.   
\end{theorem}

\proof Define a new quantity 
\begin{equation*}
Z:= v-{\varphi(u,t)}
\end{equation*}
where $\varphi$ is a smooth positive function on $[-M,M]\times(0,T')$ chosen so that $\varphi\rightarrow\infty$ as $t\rightarrow0$.  This means that $Z$ is strictly negative initially,  regardless of the initial gradient.
  
The evolution equation for $\varphi$  is given by
\begin{align}  \label{ev eqn for varphi}
&\left(\frac d {dt} - \Delta\right)\varphi=\varphi_t-\varphi_{uu}\nabsq u,
\end{align}
and we can use this with the identities from Lemma \ref{evolution equation lemma} to find that
\begin{equation*}
\left(\frac d {dt} - \Delta\right)Z=-|A|^2v-2v^{-1}|\nabla v|^2-\varphi_t+\varphi_{uu}|\nabla u|^2.  
\end{equation*}

Now, suppose $(x,t)$ is the first point at which $Z$ becomes non-negative.     Since $Z$ is periodic, this is an internal spatial maximum, and (spatial) first deriatives are zero:  
$$0=\nabla Z=\nabla v -\nabla \varphi,$$  so 
\begin{align*}
v^{-1}|\nabla v|^2&=\frac{\nabsq \varphi}\varphi=\frac{{\varphi_u^2}\nabsq u}\varphi, 
\end{align*}
and the evolution equation at this point is
\begin{align*}
\left(\frac d {dt} - \Delta\right)Z&=  -|A|^2v 
-2\frac{\varphi_u^2}\varphi \nabsq u  -\varphi_t+\varphi_{uu}|\nabla u|^2  . 
\end{align*}
   
A good choice for $\varphi$ that will allow us to make the final terms negative is $\varphi(u,t)=1/\Phi(u,t)$, where $\Phi$ solves the heat equation $\Phi_t=c \Phi''$ for some $c<1$.  In this case
\begin{gather*}
\varphi'=-\Phi^{-2}\Phi',\\
\varphi''=2\Phi^{-3}(\Phi')^2 -\Phi^{-2}\Phi'', \\
\varphi_t=-\Phi^{-2}\Phi_t \\
\intertext{and the equation satisfied by $\varphi$ is}
\varphi_t= c\varphi'' -2c\frac{\varphi'^2}\varphi.
\end{gather*}
The final three terms of the evolution equation have become
\begin{align*} 
-2\frac{\varphi_u^2}\varphi \nabsq u -\varphi_t+\varphi_{uu}|\nabla u|^2 &= 
-2\frac{\Phi'^2}{\Phi^3}\nabsq u + c\frac{\Phi''}{\Phi^2}-\frac{\Phi''}{\Phi^2}\nabsq u +2\frac{\Phi'^2}{\Phi^3}\nabsq u\\
&=\frac{\Phi''}{\Phi^2}\left(c-\nabsq u\right) \\
&=\frac{\Phi''}{\Phi^2}\left(c-1+\Phi^2\right).
\end{align*}   
In the last line we have used that, with respect to a local orthonormal frame on $M_t$, $\nabla u=\langle e_i,\omega\rangle e_i$, while $\omega$ has unit length with $1= |\omega|^2= \sum_{i=1}^n|\langle e_i,\omega\rangle e_i|^2 +|\langle \nu,\omega\rangle \nu|^2$:  it follows that      
$$\nabsq u = \sum_{i=1}^n \langle e_i,\omega\rangle ^2 = 1-\langle \nu,\omega\rangle ^2 =
\left(1-\frac1{v^2}\right)=\left(1-\Phi^2\right),$$ 
 the last equality holding only at a maximum point.

If we let $\Phi$ be a fundamental solution of the heat equation 
\begin{equation} \label{definition of p for the silly chapter}\Phi(u,t)=\frac1{\sqrt t}\exp\left(-c\frac{(u\pm2M)^2}{4t}\right),\end{equation}
we can choose $T'$ and $c$ depending only on $M$ so that $\Phi''\ge0$ and $c-1+\Phi^2\le0$  for $t<T'$. 

So, at the first interior point where $Z=0$, $Z_t\le0$ and so $Z\le0$ for $t<T'$.
\halmos

\begin{theorem}[Interior estimate for mean curvature flow]\label{interior mcf ch 9} Let $\mathbf{F}$ be a smooth solution to mean curvature flow \eqref{mean curvature flow} which is    a graph over a ball in the hyperplane $B_R(0)$.   Then we have the interior estimate
\begin{equation*}
v\le t^{q/2}\exp\left(\frac{cq(u+2M)^2}{4t}\right)(R^2-2nt-|{\mathbf{x}}|^2+u^2)^{-1}
\end{equation*}
for $0\le t\le T'$, where $q>1$, $c$,  and $T'$ depend on $M$ and $R$.   
\end{theorem}

\proof We replace $\varphi$ in our previous definition of $Z$ by $\varphi/\eta$: 
\begin{equation*}
Z:= v-\frac {\varphi(u,t)}\eta,
\end{equation*}
where a smooth positive function $\eta$ is chosen so that $Z<0$ on the boundary of a ball of shrinking radius $B_{\sqrt{R^2-2nt}}$, and, as before, $\varphi\ge0$ is chosen so that $Z<0$ at the initial time.  

In particular,  choose $\eta=R^2-2nt-|{\mathbf{x}}|^2+u^2$.
 The evolution equation for $\eta$ is given by
\begin{align*}
&\left(\frac d {dt} - \Delta\right)\eta=-2\nabsq u, 
\end{align*}
and we can use this,   the identities from Lemma \ref{evolution equation lemma}, and the evolution equation for $\varphi$ \eqref{ev eqn for varphi} to find that
\begin{multline*}
\left(\frac d {dt} - \Delta\right)Z=-|A|^2v-2v^{-1}|\nabla v|^2-\frac1\eta\left(\varphi_t-\varphi_{uu}|\nabla u|^2\right) \\ 
-2\frac\varphi{\eta^2}\nabsq u
-2\frac{\varphi_u}{\eta^2}\nabla u\cdot\nabla\eta +2\frac\varphi{\eta^3}|\nabla\eta|^2. 
\end{multline*}

Now, suppose $(x,t)$ is an internal point of this domain at which $Z$ first becomes non-negative. 

At    an internal spatial maximum of $Z$, $\nabla Z=0$ so  
\begin{align*}
|\nabla v|^2&= \nabsq{\left(\frac \varphi\eta\right)} \\
&=\frac{\varphi_u^2}{\eta^2}\nabsq u -2\frac{\varphi\varphi_u}{\eta^3} \nabla u\cdot \nabla \eta +\frac{\varphi^2}{\eta^4}\nabsq \eta. 
\end{align*}

Use this to replace the $v^{-1}\nabsq v$ term in the evolution equation,  so that at this point
\begin{align*}
\left(\frac d {dt} - \Delta\right)Z&=  -|A|^2v 
-2\frac\eta\varphi \left(\frac{\varphi_u^2}{\eta^2}\nabsq u -2\frac{\varphi\varphi_u}{\eta^3} \nabla u\cdot \nabla \eta +\frac{\varphi^2}{\eta^4}\nabsq \eta\right) \\
&\phantom{spacesp}-\frac1\eta\left(\varphi_t-\varphi_{uu}|\nabla u|^2\right) -2\frac\varphi{\eta^2}\nabsq u -2\frac{\varphi_u}{\eta^2}\nabla u\cdot\nabla\eta +2\frac\varphi{\eta^3}|\nabla\eta|^2  \\
&= -|A|^2v  -\frac1\eta\left(\varphi_t-\varphi_{uu}|\nabla u|^2+2 \nabsq u \frac{\varphi_u^2}{\varphi} \right)  
\\&\phantom{spacesp} -2 \nabsq u\frac\varphi{\eta^2} 
+  2 \frac{\varphi_u}{\eta^2} \nabla u\cdot \nabla \eta . 
\end{align*}

We can
 bound the $\nabla u\cdot\nabla \eta$ term by a $v^{-1}$ term, since (using a local orthonormal frame $\lbrace e_i\rbrace$)
\begin{align*}
\nabla u\cdot\nabla \eta&=\nabla\langle {\mathbf{x}},\omega\rangle\cdot \nabla\left(-\langle {\mathbf{x}},{\mathbf{x}}\rangle +\langle {\mathbf{x}},\omega\rangle^2\right) 
\\&=
\langle \nabla_i {\mathbf{x}},\omega\rangle e_i\cdot \left(-2\langle {\mathbf{x}},\nabla_j{\mathbf{x}}\rangle e_j+2\langle {\mathbf{x}},\omega\rangle\langle\nabla_j{\mathbf{x}},\omega\rangle e_j\right)
\\&=
\langle e_i,\omega\rangle e_i\cdot \left(-2\langle {\mathbf{x}},e_j\rangle e_j+2\langle {\mathbf{x}},\omega\rangle\langle e_j,\omega\rangle e_j\right) 
\\&=
g^{ij}\langle e_i,\omega\rangle \left(-2\langle {\mathbf{x}},e_j\rangle+2\langle {\mathbf{x}},\omega\rangle\langle e_j,\omega\rangle \right),
\end{align*}
and  {\allowdisplaybreaks{\begin{gather*} g^{ij}=\delta^{ij}-\dfrac{D_iuD_ju}{1+|Du|^2} \\  
\langle e_i,\omega\rangle=D_iu \\
\langle e_i,{\mathbf{x}}\rangle= {\mathbf{x}}_i+ uD_iu 
\end{gather*}
so that (writing $x$ for ${\mathbf{x}}-u\omega$, the position in the hyperplane) }}
 \begin{align*}
\nabla u\cdot\nabla \eta&= -2\langle Du, x\rangle +2\frac{\langle Du, x\rangle |Du|^2}{1+|Du|^2}\\
&=-2\frac{\langle Du, x\rangle }{1+|Du|^2}\\
&\le 2 \frac {|x|}{\sqrt{1+|Du|^2}} \\
&\le 2\frac {\sqrt{R^2-2nt}} v. 
\end{align*}
With this we estimate the term in the evolution equation---
\begin{equation}
 2  \frac{\varphi_u}{\eta^2} \nabla u\cdot \nabla \eta
\le 4 \frac{|\varphi_u|}{\varphi\eta} R. \label{estimate for nabla u nabla eta}%
\end{equation}
The evolution equation itself becomes
\begin{align*}
\left(\frac d {dt} - \Delta\right)Z&\le  -|A|^2v -2\nabsq u \frac\varphi{\eta^2} \\
&\phantom{space} -\frac1\eta\left(\varphi_t-\varphi_{uu}|\nabla u|^2 + 2 \nabsq u \frac{\varphi_u^2}{\varphi} -4 R%
\frac{|\varphi_u|}{\varphi}\right) .
\end{align*}

This time, we choose $$\varphi(u,t)=\Phi(u,t)^{-q}$$ for some $q>1$, where $\Phi$ still satisfies the heat equation $\Phi_t=c\Phi''$.  Then
{\allowdisplaybreaks{\label{pagereference for this part}
\begin{gather}
\varphi'=-q\Phi^{-q-1}\Phi', \notag \\
\varphi''=q(q+1)\Phi^{-q-2}(\Phi')^2 -q\Phi^{-q-1}\Phi'', \notag \\
\varphi_t=-q\Phi^{-q-1}\Phi_t, \notag \\
\intertext{so that the equation satisfied by $\varphi$ is} 
\varphi_t= c\varphi'' -c\left(1+\frac1q\right)\frac{\varphi'^2}\varphi.
\notag
\end{gather}}}
The final term in the evolution equation is
\begin{align*} 
\varphi_t-\varphi_{uu}|\nabla u|^2 + 2 &\nabsq u \frac{\varphi_u^2}{\varphi} -4  R\frac{|\varphi_u|}{\varphi}\\
&=q\frac{\Phi''}{\Phi^{q+1}}\left(-c+\nabsq u\right)+q\frac{\Phi'^2}{\Phi^{q+2}}\nabsq u\left(q-1\right)-4Rq\frac{|\Phi'|}\Phi\\
&=q\frac{\Phi''}{\Phi^{q+1}}\left(-c+1-\Phi^{2q}\right)+q\frac{|\Phi'|}{\Phi}\left((q-1)(1-\Phi^{2q})\frac{|\Phi'|}{\Phi^{q+1}}-4R\right),
\end{align*}
where we have used $\nabsq u=1-v^{-2}=1-\Phi^{2q}$.  The first term above is positive if, as in the previous case, $\Phi''\ge0$ and $c-1+\Phi^{2q}\le 0$.  The second term is positive if $q$ satisfies
\begin{equation*} 
q\ge1+\frac{4R\Phi^{q+1}}{|\Phi'|(1-\Phi^2)}.\end{equation*}

Choose  $\Phi$ to be a fundamental solution of the heat equation, as in the previous proof ---
 $$\Phi(u,t)=\frac1{\sqrt t}\exp\left(-c\frac{(u\pm2M)^2}{4t}\right).$$
Now we can choose some $c<1$ and $T'$ small so that $\Phi''>0$ and $\Phi\le(1-c)<1$ for $t<T'$.  Here, $T'$ is dependent on $M$ and $c$ only.   We can also find 
$q$ dependent on $T'$, $M$, $c$ and $R$, satisfying
\begin{equation*} 
1+\frac{4R\Phi^{q+1}}{|\Phi'|(1-\Phi^2)}\le 1+\frac{8RT'}{cM(1-(1-c)^2)}\le q.
\end{equation*}

At an internal maximum of $Z$, $Z_t\le0$ and the result follows.
\halmos

{\sloppy{
\section{Gradient estimates for anisotropic mean curvature flows}}}\label{AMCF section}
A more general case of curvature flows is that of {\emph{anisotropic mean  curvature flow}}.  This has been specifically studied by Almgren, Taylor and Wang \cite{almgren:1993}, Gurtin and Angenent \cite{angenent-gurtin:aniso}, and Andrews \cite{andrews:convexcurves,andrews:anisotropic}, among others.   The anisotropic surface energy arises in applications from materials science, such as crystalline growth and phase changes; it also arises in Finsler geometry \cite{bellettini:1996} (on a  Finsler manifold, at each point only a normed space is defined, rather than an inner product space as on a Riemannian manifold).

  In this section, we use the framework and notation of \cite{andrews:anisotropic}.

As before, we consider surfaces that can be written (either entirely or locally) as graphs, so that $M_t=\lbrace \left(x^1,\dots,x^n,u(x^1,\dots,x^n,t)\right)\rbrace=\text{graph}\,u(x,t)$.  

The equation for motion of the graph by anisotropic mean curvature is derived in \cite{andrews:anisotropic};   in the present work, we set the homogeneous degree zero ``mobility factor'' $m$ to be identically $1$, and so   
\begin{equation}
u_t=\left.F(Du)D^{ij}F\right|_{Du}u_{ij}, \label{anisotropic curvature flows}
\end{equation}
where $u_{ij}=D^2u(e_i,e_j)$, with respect to some basis for the tangent space $\lbrace e_1,\dots,e_n\rbrace$.      

The function $F:\bigR^n\rightarrow\bigR$ is defined by $F(p_1,\dots,p_n):=\bar F(p_i\phi^i-\phi^0)$, where $\lbrace \phi^0,\phi^1,\dots,\phi^n\rbrace$ is a basis for the cotangent space $V^*$, with   dual basis for $V=\bigR^{n+1}$ given by $\lbrace e_0,e_1,\dots,e_n\rbrace$, and where
$\bar F:V^*\rightarrow \bigR$   is a positive convex function that is homogeneous of degree one, $\bar F(\lambda \omega)=\lambda \bar F(\omega)$ for $\lambda>0$.  The  \emph{unit ball} of $\bar F$
$$\bar F^{-1}(1): =\lbrace \omega\in V^*: \bar F(\omega)=1 \rbrace$$  must be strictly convex.
 We also require that $\bar F$ is at least $C^3$.  
\nomenclature[F]{$\bar F$}{generating function for anisotropic mean curvature flow \refpage}

\medskip
\subsection*{Differences between the isotropic and anisotropic cases}

The introduction of the   unspecified anisotropic $\bar F$ into the flow has the effect of highlighting the special nature of the \emph{isotropic} case, when $\bar F(p_i\phi^i +p_0\phi^0)=\left(\sum_{i=0}^n p_i^2\right)^{1/2}$.  
 
One immediately notices that in the isotropic case, the term with third derivatives arising in the evolution equation is zero.  In \eqref{evolution eqn for F, first time}, this is the term   
$$\left.D\bar F\right|_z(\phi^k)  
\left.D\left(\bar F D^2\bar F\right) \right|_z(\phi^m,\hat\phi^i,\hat\phi^j)u_{mk}u_{ij}. $$

This absence of third derivatives  is apparent in the third identity of Lemma \mbox{\ref{evolution equation lemma} --- }
\begin{equation*} 
 \left(\frac d {dt} - \Delta\right)v=-|A|^2v-2v^{-1}|\nabla v|^2; \end{equation*} 
the left-hand side involves second derivatives of the gradient so we might expect to see some derivatives of curvature in the right hand side --- instead we see only curvature terms 
and first derivatives of the gradient function.

 The second difference is that there is no estimate of the form
$$\nabla u\cdot\nabla \eta\le \frac c v,$$
as in \eqref{estimate for nabla u nabla eta} for the isotropic case.  An equivalent estimate in the anisotropic context would be:   given $q=\sum_{i=1}^nq_i\phi^i$,  find $c=c(q)$ so that 
$$ \bar F(p-\phi^0)  \left.\bar F D^2\bar F\right|_{p-\phi^0}(p,q)\le c(q)$$
for all $p=\sum_{i=1}^np_i\phi^i$.   This is certainly true if we restrict $p$ to the unit ball, $\bar F(p)=1$.  If we replace $p$ by $sp$ ($s$ is a scalar), 
the putative estimate would be
\begin{equation*}
\bar F(s\bar p-\phi^0)  \left.\bar F D^2\bar F\right|_{s\bar p-\phi^0}(s\bar p,q)\le c(q).
\end{equation*}
Rewriting the left-hand  side using homogeneity gives us
\begin{equation*}
s\bar F(\bar p-\phi^0/s)  \left.\bar F D^2\bar F\right|_{\bar p-\phi^0/s}(\phi^0,q)\le c(q).
\end{equation*}
As $s$ increases, the left-hand side is converging to a constant defined on the unit ball, $F(\bar p)  \left.\bar F D^2\bar F\right|_{\bar p}(\phi^0,q)$, multiplied by $s$.  Unless 
$F(\bar p)  \left.\bar F D^2\bar F\right|_{\bar p}(\phi^0,q)$
 is zero %
, the left-hand side will not remain bounded by the right-hand side as $s\rightarrow\infty$.   The estimate will not hold without further restrictions on $\bar F$.

\subsection*{Calculating with the homogeneous function $\bar F$}

We make some observations about properties arising directly from the homogeneity and convexity of $\bar F$, and introduce some notation.

Let $\lbrace \phi^0,\dots,\phi^n \rbrace$ be a basis for the cotangent space $V^*$ dual to $\lbrace e_0,\dots,e_n\rbrace$, the basis for the tangent space $V$.   Both $V^*$ and $V$ are copies of $\bigR^{n+1}$.

For $p=p_i\phi^i$ (all repeated indices are summed  from $1$ to $n$ unless indicated otherwise) we will write $$z:=p-\phi^0.$$

In general we will prefer to write all derivatives of $\bar F$ in a form that is homogeneous of   degree zero, that is, as $D\bar F$, $\bar FD^2\bar F$, or $\bar F^2D^3\bar F$.    This means that we can evaluate them on the unit ball, or scale as we wish --- for example, we can use $\left.\bar FD^2\bar F\right|_{p-\phi^0/t}$ instead of $\left.\bar FD^2\bar F\right|_{tp-\phi^0}$. 

Homogeneity also means that some derivatives in the radial direction disappear --- for all $\omega \in V^*$,%
{\allowdisplaybreaks{\begin{gather}
\left.D\bar F\right|_\omega(\omega)=\bar F(\omega) \label{purely radial derivatives} \\
\left.\bar F D^2\bar F\right|_\omega(\omega,\cdot)=\left.\bar F D^2\bar F\right|_\omega(\cdot,\omega)=0\label{no radial second derivatives}  \\
\left. D(\bar F D^2\bar F) \right|_\omega(\omega,\cdot,\cdot)=0. \label{third derivatives of F zero for some radial parts}
\end{gather}}}

The strict convexity of the unit ball of $\bar F$ means that for all $\omega,r\in V^*$ on the unit ball, with $r\not=\pm \omega$,  
$$\left.D^2\bar F\right|_\omega(r,r)>0.$$
As $\bar F$ is homogeneous, all the level sets of $\bar F$ are also strictly convex, so this holds for all non-zero $\omega$.

 We denote by  $\widehat{\phantom{i} }$  the removal of a component in the direction of $z$ from $\phi^k$, $k=0,\dots,n$,
$$ \widehat\phi^k:=\phi^k-c^kz,$$
where $c^k$ is such that $\widehat\phi^k$ is tangent to the unit ball of $\bar F$,
\begin{align*}
0&=\left.D\bar F\right|_{ z}(\widehat\phi^k) \\&=\left. D\bar F\right|_{z}(\phi^k-c^k z)\\
&=\left.D\bar F\right|_{ z}(\phi^k)-c^k \left. D\bar F\right|_{z}(z)\\
&=\left.D\bar F\right|_{ z}(\phi^k)-c^k \bar F (z).  
\end{align*}
We have used %
\eqref{purely radial derivatives} in the last   line.
It follows that \begin{equation*} c^k=\frac{\left.D\bar F\right|_{z}(\phi^k)}{ \bar F (z)}.\label{hat h} \end{equation*}

   In the next two lemmas, we show that the coefficients of the evolution operator satisfy a condition similar to the control on degeneracy that we required with condition \eqref{degeneracy for alpha} of  Chapter \ref{higher dimensional chapter}.    

\begin{lemma} \label{lemma defining A and P for anisotropic case}
Let $\bar F:V^*\rightarrow \bigR$ be a $C^2$, positive, homogeneous degree one function with a strictly convex unit ball $\bar F^{-1}(1)=\lbrace \omega: \bar F(\omega)=1\rbrace$.  

Let $\phi^0,\dots,\phi^n$ be a basis for $V^*$.
 
Then there exist positive constants $A$ and $P$ so that 
$$\left. \bar F D^2 \bar F\right|_{p-\phi^0}(p,p)\ge A $$
for all $p=\sum_{i=1}^n p_i\phi^i $ with  $\bar{ F}(p)\ge P.$\end{lemma}

\proof
Write 
\begin{equation*}
B(p):= \left. \bar F D^2 \bar F\right|_{p-\phi^0}(p,p)= \left. \bar F D^2 \bar F\right|_{p-\phi^0}(\widehat p,\widehat p),
\end{equation*}
where 
\begin{equation*}
\widehat p= p- \frac{\left.D\bar F\right|_{p-\phi^0}(p)}{F(p-\phi^0)}\left(p-\phi^0\right) 
\end{equation*}
is non-zero whenever $p\not=0$.  As $\widehat p$ is a non-zero tangent covector, the strict convexity of the unit ball means that $B(p)>0$ if $p\not=0$.  

Fix $p=p_i\phi^i$ on the unit ball of $\bar F$, $\bar F(p)=1$.   

Consider $B(sp)$.  We want to show that we can find some $P_p$ and some strictly positive $A_p$ so that for all $s\ge P_p$, $B(sp)\ge A_p$.    If  this is not possible, then we can find a sequence $s_k\rightarrow \infty$ with
$$\lim_{k\rightarrow\infty} B(s_kp)=0.$$
 
Since $\bar F$ is $C^2$, $\left.\bar F D^2 \bar F\right|_z$ is continuous in $z$, and we have that
\begin{align*}
\lim_{k\rightarrow\infty}  \left. \bar F D^2 \bar F\right|_{s_kp-\phi^0}(s_kp,s_kp)
&=\lim_{k\rightarrow\infty}  \left. \bar F D^2 \bar F\right|_{s_kp-\phi^0}(\phi^0,\phi^0) \\
&=\lim_{k\rightarrow\infty}  \left. \bar F D^2 \bar F\right|_{p-\phi^0/s_k}(\phi^0,\phi^0) \\
&= \left. \bar F D^2 \bar F\right|_{\lim_{k\rightarrow\infty} p-\phi^0/s_k}( \phi^0,\phi^0) .
\end{align*}
Clearly, 
$$ \lim_{k\rightarrow\infty} p-\phi^0/s_k= p,$$
so
\begin{align*} 
\lim_{k\rightarrow\infty}B(s_kp)%
&= \left. \bar F D^2 \bar F\right|_{p}\left(\phi^0, \phi^0\right)= \left. \bar F D^2 \bar F\right|_{p}\left(\widehat\phi^0, \widehat\phi^0\right)>0,
\end{align*}
by the strict convexity of the unit ball, as at $p$,
$$\widehat\phi^0= \phi^0- \frac{\left.D\bar F\right|_{p}(\phi^0)}{\bar F(p)}p$$ 
is a non-zero tangent covector.
The contradiction implies that we can indeed find such $P_p$ and $A_p$.  

We can find such $P_q$ and $A_q$ for every $q=q_i\phi^i$.  
  Let
$$A:=\inf_{q: \bar F(q)=1}A_q, \qquad\qquad P:=\sup_{q: \bar F(q)=1}P_q . $$
As we are optimizing over a compact space, %
$A>0$ and $P<\infty$. 
 
The result follows directly, for given any $p=p_i\phi^i$ with $\bar F(p)\ge P$,
\begin{equation*}
B(p)\,=\,B(\bar F(p) \bar p) 
\,\ge\, A_{\bar p}
\,\ge\, A,
\end{equation*}
where $\bar p=p/\bar F(p)$ is on the unit ball. 
\halmos

We use this to show that the anisotropic mean curvature flow satisfies \eqref{degeneracy for alpha}, 
 the condition controlling the degeneracy of the parabolic operator in Chapter \ref{higher dimensional chapter}.  

\begin{lemma} \label{lemma connecting amcf and alpha} For all non-zero $v=v_i\phi^i$ and $p=p_i\phi^i$, we can find positive constants $P$ and $A_0$ such that  
\begin{equation}
\frac{|p|^4}{(v\cdot p)^2} \left.\bar F D^2 \bar F\right|_{p-\phi^0}  (v,v) \ge A_0
\label{blah blah blah} \end{equation} 
whenever $\bar F(p)\ge P$.  Here, $(v\cdot p)^2=\sum (v_ip_i)^2$ and $|p|^2=p\cdot p$. 
\end{lemma}
\proof
Set $B(p,v)=\frac{|p|^4}{(v\cdot p)^2} \left.\bar F D^2 \bar F\right|_{p-\phi^0}  (v,v)$.
 
Since $B(\cdot,v)$  is invariant under $v\mapsto sv$, we need only to consider $v$ in the unit ball.

Suppose that $p$ is in the unit ball.  Since $\widehat v$ is a non-zero tangent covector at $p-\phi^0$, $\left.\bar F D^2 \bar F\right|_{p-\phi^0}  (v,v)>0$ by the strict convexity of the unit ball.  By compactness, 
\begin{equation*}
\inf_{p\in \bar F^{-1}(1)} \inf_{v\in \bar F^{-1}(1)}  \left.\bar F D^2 \bar F\right|_{p-\phi^0}  (v,v)\ge c_1 >0.
\end{equation*} 
Also, as neither $p$ nor $v$ are zero, 
\begin{equation*}
\frac{|p|^4}{(v\cdot p)^2} \ge c_2>0, 
\end{equation*}
and so $\inf_{p,v\in F^{-1}(1)} B(p,v)\ge c_1c_2>0$.

Suppose, in order to obtain a contradiction, that there is a pair $(v,p)$ in the unit ball for which there are no such constants $A_0$ and $P$.  That is,    
\begin{equation*}
\lim_{s\rightarrow \infty} B(sp,v) = \lim_{s\rightarrow \infty}\frac{s^2|p|^4}{(v\cdot p)^2}\left.\bar F D^2 \bar F\right|_{sp-\phi^0}  (v,v) =0.
\end{equation*} 

There are two possibilities here:   $v\not=p$ or $v=p$.  In the first case,  we must have
\begin{align*}
0&=\lim_{s\rightarrow \infty}  \left.\bar F D^2 \bar F\right|_{sp-\phi^0}  (v,v)
\\&= \lim_{s\rightarrow \infty}  \left.\bar F D^2 \bar F\right|_{p-\phi^0/s}  (v,v) 
\\&=\left.\bar F D^2 \bar F\right|_{p}  (v,v).
\end{align*} 
However, since $v\not=p$,  $\left.\bar F D^2 \bar F\right|_{p}  (v,v)>0$ which is a contradiction.

On the other hand, if $v=p$, then
\begin{align*}
\lim_{s\rightarrow \infty} B(sp,p) &= \lim_{s\rightarrow \infty} {s^2}\left.\bar F D^2 \bar F\right|_{sp-\phi^0}  (p,p)\\
&=\lim_{s\rightarrow \infty}\left.\bar F D^2 \bar F\right|_{sp-\phi^0}  (sp,sp)\\
&\ge A
\end{align*} 
by Lemma \ref{lemma defining A and P for anisotropic case}, which is again a contradiction.

Therefore for every pair of covectors $(v,p)$, there is a pair of positive constants $A_0$ and $P$ such that  $B(p,v)\ge A_0$ whenever $F(p)\ge P$.   To get bounds for all $(v,p)$ we take the infimum of the $A_0$ and the supremum of the $P$.  
\halmos

We will consider two different restrictions on $\bar F$.  The first is that third derivatives are small; the second is a symmetry in the distinguished direction $\phi^0$.

In order to define the first condition, consider the tensor 
\begin{equation}
{Q}(p,q,r):= \bar F^2(z) \left.D^3 \bar F \right|_z(p,q,r), \label{defn of Q}
\end{equation}
for $p,q,r$ covectors tangent to the unit ball of $\bar F$ at $z$, so that $\left.D\bar F\right|_z(p)=0$ (and similarly for $q$ and $r$).  
   
(This is the \emph{Cartan tensor} of Bao, Chern and Shen \cite{chern:finsler}, or the tensor $ Q$ of \cite{andrews:anisotropic} restricted to the tangent space of the unit ball.)

The \emph{smallness-of-third-derivatives condition} is that  $Q$ satisfies
\begin{equation}
Q(p,q,r)\le C_1\left[\bar F(z)^3\left.D^2\bar F\right|_z(p,p)\left.D^2\bar F\right|_z(q,q)\left.D^2\bar F\right|_z(r,r)\right]^{1/2} \label{condition on D^3 F},
\end{equation}
for all $p,q,r$ tangent to the unit ball of $\bar F$, where $C_1$ is a positive constant dependent on $n$.

The \emph{symmetry condition} is that 
\begin{equation}
\bar F(p+\phi^0)=\bar F(p-\phi^0) \text{ for all $p=\sum_{i=1}^n p_i\phi^i$.}\label{symmetry condition} 
\end{equation}

\begin{lemma} If $\bar F$ satisfies \eqref{symmetry condition}, then 
{\allowdisplaybreaks{\begin{gather}
\left.D\bar F \right|_p(\phi^0)=0  \label{symmetry condition implies some first derivatives zero} \\
\left. D^2 \bar F \right|_p(\phi^0,\phi^j)=0 \notag \\
\left. D^3 \bar F\right|_p(\phi^0,\phi^j,\phi^k)=0 \label{third consequence of symmetry}\\
\left. D^3 \bar F\right|_p(\phi^0,\phi^0,\phi^0)=0, \label{fourth consequence of symmetry}
\end{gather}}}
for all $p= \sum_{i=1}^n p_i\phi^i$ and all $j,k\not=0$.  
\end{lemma}
\proof This is a direct consequence of homogeneity.  \halmos

We can show that the symmetry condition \eqref{symmetry condition}  can be used in a similar way to the smallness-of-third-derivatives condition \eqref{condition on D^3 F}.

\begin{lemma}
 Suppose the symmetry condition \eqref{symmetry condition} holds.    Then a constant dependent only on $\bar F$ bounds
\begin{equation} \label{some equation in lemma 9.7} 
\left| \frac{\left.\bar F D(\bar FD^2\bar F)\right|_{p-\phi^0}(p,\widehat q,\widehat q)}
{G(p,p)^{1/2}G(q,q)} \right|,
\end{equation} 
for all  $p=\sum_{i=1}^np_i\phi^i$ and $q=\sum_{i=1}^nq_i\phi^i$, where $G=\left.\bar F D^2\bar F\right|_{p-\phi^0}$.     Furthermore,  for all $\epsilon>0$ we can find $S_\epsilon$  so that when 
$\bar F(p)\ge S_\epsilon$, %
\begin{equation*} 
\left| {\left.\bar F D(\bar FD^2\bar F)\right|_{p-\phi^0}(p,\widehat q,\widehat q)}\right| 
\le \epsilon {G(p,p)^{1/2}G(q,q)}. 
\end{equation*}

  \label{lemma showing that symmetry is as good as small third derivatives}
\end{lemma}

\proof 
Consider $p$ and $q$ %
on the unit ball and set 
\begin{equation*}  C:=\sup_{\substack{ p,q\in \bar F^{-1}(1) \\ p=p_i\phi^i, \, q=q_i\phi^i }}
\frac{\left.\bar F D(\bar FD^2\bar F)\right|_{p-\phi^0}(p,\widehat q,\widehat q)}
{G(p,p)^{1/2}G(q,q)} 
\end{equation*}

When we project  $p$ and $q$ onto the tangent plane at $p-\phi^0$ they give non-zero tangent covectors $\widehat p$ and $\widehat q$, %
\begin{gather*}
\widehat p= p- c^p(p-\phi^0), \qquad\qquad
\widehat q= q- c^q(p-\phi^0),
\end{gather*}  
where 
\begin{equation*}
c^p=\frac{\left.D\bar F\right|_{p-\phi^0}(p)}{\bar F(p-\phi^0)},\qquad\qquad
c^q=\frac{\left.D\bar F\right|_{p-\phi^0}(q)}{\bar F(p-\phi^0)}, 
\end{equation*}  
so the terms in the denominator of \eqref{some equation in lemma 9.7},  $G(p,p)$ and $G(q,q)$, %
are strictly positive,  and hence bounded below when we take the supremum over $p$ and $q$ in the unit ball.%

Also,  $\bar F D(\bar FD^2\bar F)$ is a homogeneous degree zero tensor, and so bounded above on the unit ball.   
 
It follows that $C$ is finite.

The constant $C$ is unchanged if we scale $q\mapsto sq$ %
so we only need to consider the behaviour of \eqref{some equation in lemma 9.7} as $\bar F(p)$ becomes large:  that is,
\begin{equation}
\lim_{s\rightarrow\infty} 
\frac{\left.\bar F D(\bar FD^2\bar F)\right|_{sp-\phi^0}(sp,\widehat q,\widehat q)}
{G(sp,sp)^{1/2}G(q,q)},\label{crucial quantity in third derivative estimate}
\end{equation}
where $G=\left.\bar F D^2\bar F \right|_{sp-\phi^0}$.  

Firstly,  consider the case that $q$ is parallel to $p$.   Let $p$ %
be on the unit ball, and without loss of generality, let $q=+p$.  

We note that the $G(q,q)=G(p,p)$ term in the denominator converges to zero,
\begin{align*}
\lim_{s\rightarrow\infty}  G(p,p)&=%
\lim_{s\rightarrow\infty} \left.\bar FD^2\bar F\right|_{p-\phi^0/s}(p,p)\\
&= \left.\bar FD^2\bar F\right|_{p}(p,p)=0,
\end{align*}  
so to deal with this we will multiply both top and bottom by $s^2$: %
\begin{align}
\lim_{s\rightarrow\infty} &
\frac{\left.\bar F D(\bar FD^2\bar F)\right|_{sp-\phi^0}(sp,\widehat p,\widehat p)}
{G(sp,sp)^{1/2}G(p,p)}
\notag \\&= 
\lim_{s\rightarrow\infty}
\frac{s^2\left.\bar F D(\bar FD^2\bar F)\right|_{sp-\phi^0}\left(\phi^0,p-c^p(sp-\phi^0),p-c^p(sp-\phi^0)\right)}
{s^2G(sp,sp)^{1/2}G(p,p)}
\notag \\&= \lim_{s\rightarrow\infty}
\frac{\left.\bar F D(\bar FD^2\bar F)\right|_{sp-\phi^0}\left(\phi^0,s(1-sc^p)p+sc^p\phi^0,s(1-sc^p)p+sc^p\phi^0\right)}
{G(sp,sp)^{3/2}} \label{line from the estimate on third derivatives}.
\end{align}

Now, the denominator is strictly positive, and by Lemma \ref{lemma defining A and P for anisotropic case}  bounded below~--- 
\begin{equation*}
\lim_{s\rightarrow\infty}  G(sp,sp)=\lim_{s\rightarrow\infty} \left.\bar FD^2\bar F\right|_{sp-\phi^0}(sp,sp)\ge A>0.
\end{equation*}
 
The limiting value of the coefficient of $\phi^0$ in \eqref{line from the estimate on third derivatives} is
{\allowdisplaybreaks[1]{\begin{align*}
\lim_{s\rightarrow\infty} sc^p&=\lim_{s\rightarrow\infty} s \frac{\left.D\bar F\right|_{p-\phi^0/s}(p)}{\bar F(sp-\phi^0)}\\*
&=\lim_{s\rightarrow\infty}  \frac{\left.D\bar F\right|_{p-\phi^0/s}(p)}{\bar F(p-\phi^0/s)}\\*
&= \frac{\left.D\bar F\right|_{p}(p)}{\bar F(p)}\\*
&=1,
\end{align*} }}
using \eqref{purely radial derivatives}.  

The limiting value of the coefficient of $p$ is   
\begin{align*}
\lim_{s\rightarrow\infty} s(1-sc^p)&=\lim_{s\rightarrow\infty} s\left(1-s\frac{\left.D\bar F\right|_{p-\phi^0/s}(p)}{\bar F(sp-\phi^0)}\right)\\*
&=\lim_{r\rightarrow0} \frac1r\left(1-\frac{\left.D\bar F\right|_{p-r\phi^0}(p)}{\bar F(p-r\phi^0)}\right)\\*
&=\lim_{r\rightarrow0} \frac1r \left(\frac{\left.D\bar F\right|_{p}(p)}{\bar F(p)}
-\frac{\left.D\bar F\right|_{p-r\phi^0}(p)}{\bar F(p-r\phi^0)}\right)
\end{align*}
where we use that ${\left.D\bar F\right|_{p}(p)}={\bar F(p)}$.   The above term is a derivative, so we have
\begin{align*}\lim_{s\rightarrow\infty} s(1-sc^p)
&= \dfrac{d  }{dr}  \left.\left(\frac{\left.D\bar F\right|_{p+r\phi^0}(p)}{\bar F(p+r\phi^0)}\right)\right|_{r=0} \\
&=\frac{\left.D^2\bar F\right|_{p}(\phi^0,p)}{\bar F(p)}-\frac{\left.D\bar F\right|_p(p) \left.D\bar F\right|_p(\phi^0)}{\bar F(p)^2} \\*
&=0,
\end{align*} 
where the first term of the second last line is zero due to \eqref{no radial second derivatives} while the second term is zero as consequence \eqref{symmetry condition implies some first derivatives zero}
of the symmetry condition. 

The equation \eqref{line from the estimate on third derivatives} is then 
\begin{align*} 
 \lim_{s\rightarrow\infty}
&\frac{\left.\bar F D(\bar FD^2\bar F)\right|_{sp-\phi^0}\left(\phi^0,s(1-sc^p)p+sc^p\phi^0,s(1-sc^p)p+sc^p\phi^0\right)}
{G(sp,sp)^{1/2}G(sp,sp)} \\
&\phantom{spacespacespacespacespace}= 
\frac{\left.\bar F D(\bar FD^2\bar F)\right|_{p}\left(\phi^0,\phi^0,\phi^0\right)}   
{\left.\bar FD^2\bar F\right|_p(\phi^0,\phi^0)^{3/2}}\\
&\phantom{spacespacespacespacespace}=0, 
\end{align*}
where symmetry has been used again in the form of \eqref{fourth consequence of symmetry}. 

Now consider the case that $q$ is \emph{not} parallel to $p$.   In this case, the denominator of \eqref{crucial quantity in third derivative estimate} is non-zero:
\begin{align*}
\lim_{s\rightarrow\infty} {G(sp,sp)^{1/2}} G(q,q)&=\lim_{s\rightarrow\infty} \left.\bar FD^2\bar F\right|_{p-\phi^0/s}(\phi^0,\phi^0) \left.\bar FD^2\bar F\right|_{p-\phi^0/s}(q,q)\\
&= \left.\bar FD^2\bar F\right|_{p}(\phi^0,\phi^0)\left.\bar FD^2\bar F\right|_{p}(q,q)\\
&>0, 
\end{align*}
since $\lim_{s\rightarrow\infty}\widehat q$ is a non-zero tangent vector at $p$,
\begin{equation*}
\lim_{s\rightarrow\infty}\widehat q= 
\lim_{s\rightarrow\infty}q-
\frac{\left.D\bar F\right|_{sp-\phi^0}(q)}{\bar F(sp-\phi^0)} (sp-\phi^0)
=q-\frac{\left.D\bar F\right|_{p}(q)}{\bar F(p)} p.
\end{equation*}
This cannot be zero as $q$ is not parallel to $p$.

Then 
\begin{align*}
\lim_{s\rightarrow\infty}  &
\frac{\left.\bar F D(\bar FD^2\bar F)\right|_{sp-\phi^0}(sp,\widehat q,\widehat q)}
{G(sp,sp)^{1/2}G(q,q)} \\
&=\lim_{s\rightarrow\infty} \left({G(sp,sp)^{1/2}G(q,q)}\right)^{-1}
\lim_{s\rightarrow\infty}
{\left.\bar F D(\bar FD^2\bar F)\right|_{p-\phi^0/s}(\phi^0,\widehat q,\widehat q)}\\
&= 0,
\end{align*}
where the second limit is zero by \eqref{third consequence of symmetry}, since $\lim_{s\rightarrow\infty}\widehat q$ has no component in the direction of $\phi^0$.

We have shown that for a fixed $p$ on the unit ball, the quantity \eqref{some equation in lemma 9.7} is bounded above (by $C$), and that as $p$ is scaled outwards this decreases to zero, so for fixed $p$ the quantity is bounded above.
    
By compactness of the unit ball, it is bounded for all $p$.  
\halmos

\subsection*{Estimates for periodic, anisotropic mean curvature flows}

Let $u:\bigR^n\times[0,T]\rightarrow \bigR$ be a $H_2$,  bounded 
$$|u(x,t)|\le M,$$
periodic
$$u(x+L,t)=u(x,t) \text{ for some lattice $L$,} $$
 solution
to the anisotropic curvature flow equation \eqref{anisotropic curvature flows}, where $\bar F$ is a positive convex function, homogeneous  degree one, with  a strictly convex unit ball.  

  \remark We have one estimate in the case that $\bar F$ satisfies the smallness-of-third-derivatives condition (Theorem \ref{First estimate for periodic anisotropic curvature flows}) and another in the case that $\bar F$ satisfies the symmetry condition (Theorem \ref{Second estimate for periodic anisotropic curve...}).  However,  in Theorem \ref{periodic higherdimensional theorem} and Corollary \ref{periodic higherdimensional corollary} we found an estimate for periodic anisotropic flows without the need to impose either of these conditions.   In that case, the fact that we were estimating on  the difference quotient, rather than the first derivative itself, avoided the need to take derivatives of the flow coefficients.   Strictly speaking, Theorems \ref{First estimate for periodic anisotropic curvature flows}  and \ref{Second estimate for periodic anisotropic curve...} are redundant, but they are a good introduction for the interior estimate of Theorem \ref{interior estimate for anisotropic curve...}.

\begin{theorem}%
\label{First estimate for periodic anisotropic curvature flows}
If the tensor $Q$ given by \eqref{defn of Q}   satisfies \eqref{condition on D^3 F} with 
\begin{equation} {C_1}^2<\frac4{\sqrt{n}}\label{condition on C_1},\end{equation}
 then 
 \begin{equation*}
F(Du)\le \max \left\lbrace t^{q/2}\exp\left(\frac{Aq(|u|-2M)^2}{4t}\right), P\right\rbrace
\end{equation*}
  for $0<t\le T'$, where %
 $T'$  depends on $M$, $A$ and $P$ (both given by Lemma  \ref{lemma defining A and P for anisotropic case}), and $1<q=(1-C_1^2\sqrt{n}/4)^{-1}$.  
\end{theorem}

\begin{theorem} \label{Second estimate for periodic anisotropic curve...}
If $\bar F$ satisfies the symmetry condition \eqref{symmetry condition},
then 
 \begin{equation*}
F(Du)\le \max\left\lbrace t\exp\left(\frac{A(|u|-2M)^2}{2t}\right), P, S_{(2/n)^{1/2}}\right\rbrace
\end{equation*}
  for $0<t\le T'$, where $T'$ depends on $M$, $A$ and $P$ (both given by Lemma \ref{lemma defining A and P for anisotropic case}), and $S$   is given by Lemma \ref{lemma showing that symmetry is as good as small third derivatives}.
\end{theorem}

\smallskip {{\noindent\textbf{Proof  of Theorem \ref{First estimate for periodic anisotropic curvature flows}:}}\quad}
 As in the previous sections, define the quantity
\begin{equation*}Z:= F(Du)-\varphi(u,t)\end{equation*}
where $P$ is given by Lemma \ref{lemma defining A and P for anisotropic case}, and where $\varphi$ is a positive function chosen so that $\varphi\rightarrow\infty$ as $t\rightarrow 0$.  Suppose that we are at the first point where $Z$ is no longer negative, and let us assume that at this point, $\varphi\ge P$.  This point will be a spatial maximum of $Z$,  due to the periodicity of $Z$.

Then at this point, $F(Du)=\varphi$ and the first derivative condition is
\begin{equation*}
0=D_kZ=\left.D\bar F\right.|_z(\phi^m)u_{mk}-\varphi_u u_k,
\end{equation*} 
where $z=Du-\phi^0$.    That is, for all vectors $v\in \text{span}\lbrace e_1,\dots,e_n\rbrace$, 
\begin{equation} \label{first derivative condition for bar F}  
D^2u\left(\left.D\bar F\right|_z(\phi^m)e_m,v\right)=\varphi_uDu\left(v\right).
\end{equation} 

Using \eqref{no radial second derivatives}, we can rewrite the evolution equation for $u$ in terms of the purely tangential directions $\widehat\phi^i$,
\begin{equation*}
u_t=\left.\bar F D^2\bar F \right|_z(\phi^i,\phi^j)D^2u(e_i,e_j)=\left.\bar F D^2\bar F \right|_z(\widehat\phi^i,\widehat\phi^j)D^2u(e_i,e_j).
\end{equation*}

We make use of this in finding an evolution equation for $F$   
{\allowdisplaybreaks[2]{\begin{align}
\pd {F} t&= \left.D\bar F\right|_z(\phi^k)u_{kt} \notag \\*
&= \left.D\bar F\right|_z(\phi^k)\left[\left.\bar F D^2\bar F \right|_z(\widehat\phi^i,\widehat\phi^j)u_{ij}\right]_k \notag \\*
&= 
\left.D\bar F\right|_z(\phi^k)
\Big[  
\left.D\left(\bar F D^2\bar F\right) \right|_z(D_k z,\widehat\phi^i,\widehat\phi^j)u_{ij} +\left.\bar F D^2\bar F \right|_z(D_k\widehat\phi^i,\widehat\phi^j)u_{ij}
 \notag \\*
&\phantom{spacespacespa}
+\left.\bar F D^2\bar F \right|_z(\widehat\phi^i,D_k\widehat\phi^j)u_{ij}   
+\left.\bar F D^2\bar F \right|_z(\widehat\phi^i,\widehat\phi^j)u_{ijk}
 \Big] 
\notag \\* &\phantom{spacespacespa}
+\left.\bar F D^2\bar F\right|_{z}(\phi^i,\phi^j)D_{ij}F 
\notag\\* &\phantom{spacespacespa}
-\left.\bar F D^2\bar F\right|_{z}(\phi^i,\phi^j)\left[ \left.D^2\bar F\right|_z(\phi^m,\phi^l)u_{mi}u_{lj}+\left.D\bar F\right|_z(\phi^m)u_{mij}\right]
\notag \\
&=
\left.D\bar F\right|_z(\phi^k)
\Big[  
\left.D\left(\bar F D^2\bar F\right) \right|_z(D_k z,\widehat\phi^i,\widehat\phi^j)u_{ij} +\left.\bar F D^2\bar F \right|_z(D_k\widehat\phi^i,\widehat\phi^j)u_{ij}
 \notag \\*
&\phantom{spacespacespa}
+\left.\bar F D^2\bar F \right|_z(\widehat\phi^i,D_k\widehat\phi^j)u_{ij}  \Big]
\notag \\* &\phantom{spacespa}
+\left.\bar F D^2\bar F\right|_{z}(\phi^i,\phi^j)D_{ij}F 
-\left.\bar F D^2\bar F\right|_{z}(\phi^i,\phi^j) \left.D^2\bar F\right|_z(\phi^m,\phi^l)u_{mi}u_{lj}
\notag, 
\end{align}}}
where in the third step we have added and subtracted second derivatives of $F$.

Derivatives of $z$ are $
D_kz=u_{mk}\phi^m \label{derivative of radial z}$.
We this to simplify those terms with derivatives of $\widehat\phi^i$---
 \begin{align*}
\left.D^2\bar F\right|_z(D_k\widehat\phi^i,\widehat\phi^j)&=\left.D^2\bar F\right|_z(D_k(-
c^{\phi^i}z),\widehat\phi^j)\\
&=\left.D^2\bar F\right|_z(-D_k(c^{\phi^i})z-c^{\phi^i}u_{mk}\phi^m,\widehat\phi^j),
\intertext{and remembering that $\left.D^2\bar F\right|_z(z,\cdot)=0$, this is}
&=-c^{\phi^i}\left.D^2\bar F\right|_z(u_{mk}\phi^m,\widehat\phi^j)\\
&=-\frac{\left.D\bar F\right|_{z}(\phi^i)}{ \bar F (z)}\left.D^2\bar F\right|_z(u_{mk}\phi^m,\widehat\phi^j).
\end{align*}

The evolution equation is now
\begin{align}\pd {F} t&=
\left.D\bar F\right|_z(\phi^k)  
\left.D\left(\bar F D^2\bar F\right) \right|_z(\phi^m,\widehat\phi^i,\widehat\phi^j)u_{mk}u_{ij} 
\notag 
\\ &\phantom{space}
-\left.D\bar F\right|_z(\phi^k) \left[ \left.D\bar F\right|_z(\phi^i)  \left. D^2\bar F \right|_z(\phi^m,\widehat\phi^j) +  \left.D\bar F\right|_z(\phi^j)  \left. D^2\bar F \right|_z(\phi^m,\widehat\phi^i) \right]u_{mk}u_{ij}
\notag 
\\ &\phantom{space}
+\left.\bar F D^2\bar F\right|_{z}(\phi^i,\phi^j)D_{ij}F 
-\left.\bar F D^2\bar F\right|_{z}(\phi^i,\phi^j) \left.D^2\bar F\right|_z(\phi^m,\phi^l)u_{mi}u_{lj}.
 \label{evolution eqn for F, first time}  
\end{align}

When we are at a critical point of $Z$,  we can use the first derivative condition \eqref{first derivative condition for bar F} to simplify further.  The first term of  \eqref{evolution eqn for F, first time} becomes
\begin{align*}
D^2u( \left.D\bar F\right|_z(\phi^k) e_k,e_m)
\left.D\left(\bar F D^2\bar F\right) \right|_z&(\phi^m,\widehat\phi^i,\widehat\phi^j)u_{ij}  \\&=
\varphi_uDu(e_m)
\left.D\left(\bar F D^2\bar F\right) \right|_z(\phi^m,\widehat\phi^i,\widehat\phi^j)u_{ij} \\
&= \varphi_u
\left.D\left(\bar F D^2\bar F\right) \right|_z(Du,\widehat\phi^i,\widehat\phi^j)u_{ij}, 
\end{align*}
while the second becomes 
\begin{align*}
&- \left.D\bar F\right|_z(\phi^k){\left.D\bar F\right|_{ z}(\phi^i)}\left.D^2\bar F\right|_z(\phi^m,\widehat\phi^j)u_{mk}u_{ij}\\
&\phantom{spacespacespace}=- \left. D^2\bar F \right|_z(\phi^m,\widehat\phi^j)D^2u\left( \left.D\bar F\right|_z(\phi^k)e_k,e_m\right)D^2u\left( \left.D\bar F\right|_z\left(\phi^i\right)e_i,e_j\right) \\
&\phantom{spacespacespace}=-\varphi_u^2\left. D^2\bar F \right|_z(Du(e_m)\phi^m,Du(e_j)\widehat\phi^j)\\
&\phantom{spacespacespace}=-\varphi_u^2\left. D^2\bar F \right|_z(Du,Du),
\end{align*}
as does the third.

The evolution equation for $F$, at the local maximum of $Z$, is now
\begin{align*}
\pd {F} t&= \frac{\varphi_u}\varphi\left.\bar F D\left(\bar F D^2\bar F\right) \right|_z(Du,\widehat\phi^i,\widehat\phi^j)u_{ij}-2\frac{\varphi_u^2}\varphi  \left. \bar F D^2\bar F \right|_z(Du,Du)  
\notag \\ &\phantom{space}
+\left.\bar F D^2\bar F\right|_{z}(\phi^i,\phi^j)D_{ij}F 
-\frac1\varphi\left.\bar F D^2\bar F\right|_{z}(\phi^i,\phi^j) \left.\bar F D^2\bar F\right|_z(\phi^m,\phi^l)u_{mi}u_{lj},
\end{align*}
where we have multiplied some terms through by $1=\bar F/\varphi$ (since we assume that $Z=0$ here) in order that derivatives of $\bar F$ appear as homogeneous degree zero terms.   

Derivatives of $\varphi$ are given by
{\allowdisplaybreaks{\begin{gather*}
D\varphi=\varphi_u Du \\
D_{ij}\varphi=\varphi_{uu}u_i u_j+\varphi_u u_{ij}\\
\dfrac{d \varphi}{dt}=\varphi_uu_t+\varphi_t
\end{gather*}   }}
for $i,j\not=0$, 
so  an evolution equation for $\varphi$ is
\begin{align*}\dfrac{d \varphi}{dt}
&=\varphi_uu_t+\varphi_t +\left.\bar F D^2\bar F\right|_z(\phi^i,\phi^j)\left( D_{ij}\varphi-\varphi_{uu}u_i u_j-\varphi_u u_{ij}\right) \\
&=\varphi_t +\left.\bar F D^2\bar F\right|_z(\phi^i,\phi^j)\left( D_{ij}\varphi-\varphi_{uu}u_i u_j\right),\end{align*}
and the entire evolution equation for $Z$, at a local maximum, is  
\begin{align}
\dfrac {dZ} {dt}
&=\left.\bar F D^2\bar F\right|_{z}(\phi^i,\phi^j)D_{ij}Z+  \frac{\varphi_u}{\varphi} \left.\bar F D\left(\bar F D^2\bar F\right) \right|_{z}(Du,\widehat\phi^i,\widehat\phi^j)u_{ij}
\notag\\&\phantom{spacespa} 
-2\frac{\varphi_u^2}{\varphi}\left. \bar F D^2\bar F \right|_{ z}(Du,Du)  
-\frac1\varphi\left.\bar F D^2\bar F\right|_{z}(\phi^i,\phi^j)\left.\bar F D^2\bar F\right|_{z}(\phi^m,\phi^l)u_{mi}u_{lj} \notag
\\&\phantom{spacespa}
-\varphi_t
+\left.\bar F D^2\bar F\right|_{z}(Du,Du)\varphi_{uu} .\notag
\end{align}

Notice that all the covectors $\phi^i$, $Du$  appear in places where they may be replaced by $\widehat\phi^i$, $\widehat{Du}$ respectively (using \eqref{no radial second derivatives}  and \eqref{third derivatives of F zero for some radial parts}). 
That is, we are working exclusively on the tangent space to the unit ball.     

Restricted to the tangent space, $D^2\bar F$ is positive definite so we can define a Riemannian metric on the tangent space $G^{\alpha \beta}:=\left.\bar F D^2\bar F\right|_{ z} (\widehat\phi^\alpha,\widehat\phi^\beta)$,  $\alpha,\beta\not=0$ .  We can choose the basis $\lbrace \phi^1,\dots,\phi^n\rbrace$ so that $G$ is the identity at our maximum point,  $G^{\alpha\beta}=\delta^{\alpha \beta}$.    The evolution equation for $Z$ is now
\begin{align}
\dfrac {dZ} {dt}&=G^{ij}D_{ij}Z+\frac{\varphi_u}{\varphi} \left.\bar F D\left(\bar F D^2\bar F\right) \right|_{z}(\widehat{Du},\widehat\phi^i,\widehat\phi^j)u_{ij}-2\frac{\varphi_u^2}{\varphi}G({Du},{Du})  
\notag \\
&\phantom{spacespacespacespace}-\frac1\varphi G^{ij}G^{ml}u_{mi}u_{lj} -\varphi_t
+G({Du},{Du})\varphi_{uu}.
\label{last line of evolution equation for Z}\end{align}

Recall the Cauchy-Schwarz inequality for positive matrices:  If $A$ is a positive semi-definite $n\times n$ matrix, then for $v,w\in \bigR^n$,
$$0\le \left( \left(2\epsilon\right)^{1/2}v-   \left(2\epsilon\right)^{-1/2}w\right)^T A \left( \left(2\epsilon\right)^{1/2}v-   \left(2\epsilon\right)^{-1/2}w\right)$$
and so 
$$ v^T A w\le \epsilon v^T A v +\frac1{4\epsilon}w^T A w.$$
If $A$ is  positive definite, then we can replace $w$ by $A^{-1}w$ to find that 
$$ v^T w \le \epsilon v^T A v +\frac1{4\epsilon}w^T A^{-1} w.$$

As we assume that $u$ is smooth, $G$ is positive definite and we can use the above inequality to estimate the second term of \eqref{last line of evolution equation for Z}  : \label{page for C-s}
\begin{align}
 \frac{\varphi_u}{\varphi} &\left.\bar F D\left(\bar F D^2\bar F\right) \right|_{ z}(\widehat{Du},\widehat\phi^i,\widehat\phi^j)u_{ij} \notag
\\&\phantom{space}= \frac{\varphi_u}{\varphi} \left[
\left.D\bar F \right|_{ z} (\widehat{Du})\left. \bar F D^2\bar F \right|_{z}(\widehat\phi^i,\widehat\phi^j)+ \left. \bar F^2D^3\bar F \right|_{z}(\widehat{Du},\widehat\phi^i,\widehat\phi^j)
\right]u_{ij} \notag \\
&\phantom{space}= \frac{\varphi_u}{\varphi}u_k {Q}^{kij}
u_{ij}  \notag \\
&\phantom{space}\le  
\epsilon \frac{{\varphi_u}^2}\varphi G(Du,Du)
+\frac1{4\epsilon\varphi} G_{\alpha\beta }Q^{\alpha i j}u_{ij}Q^{\beta kl} u_{kl},  \label{cauchy-schwarz for interior estimate}
\end{align}  
where the first term of the second line is zero because $\widehat{Du}$ is tangent to the unit ball, so $\left.D\bar F \right|_{ z} (\widehat{Du})=0$.  In the last line, we have used the notation for the inverse $G_{\alpha\beta}=(G^{-1})^{\alpha\beta}$.   We will choose $\epsilon<1$  later.

We can use \eqref{condition on D^3 F}, the smallness-of-third-derivatives condition, to estimate the second term in this inequality: 
{\allowdisplaybreaks[1]{\begin{align*}
\frac1{4\epsilon\varphi} G_{\alpha\beta} Q^{\alpha i j}u_{ij}Q^{\beta kl} u_{kl} 
&=\frac1{4\epsilon\varphi} Q\left(G_{\alpha\beta}\widehat\phi^\alpha,u_{ij}\widehat\phi^i,\widehat\phi^j\right)
 Q\left(\widehat\phi^\beta,u_{kl}\widehat\phi^k,\widehat\phi^l\right)\\*
&\le\frac{ {C_1}^2}{4\epsilon\varphi} 
\bigg( G(G_{\alpha \beta}\widehat\phi^\alpha,G_{\gamma \beta} \widehat\phi^\gamma)G(u_{ij}\widehat\phi^i,u_{mj}\widehat\phi^m)G(\widehat\phi^j,\widehat\phi^j) %
\\*&\phantom{spacespacespacespace}%
\times G(\widehat\phi^\beta, \widehat\phi^\beta)G(u_{kl}\widehat\phi^k,u_{pl}\widehat\phi^p)G(\widehat\phi^l,\widehat\phi^l)
\bigg)^{1/2} \\
&=\frac{ {C_1}^2}{4\epsilon\varphi} 
\left( G_{\alpha\beta}G^{\alpha\gamma}G_{\gamma\beta} u_{ij}u_{mj}G^{im}G^{jj}G^{\beta\beta}G^{kp}u_{kl}u_{pl}G^{ll}\right)^{1/2} \\*
&=
\frac{ {C_1}^2}{4\epsilon\varphi} 
\left( G_{\beta\beta}u_{ij}u_{ij}G^{jj}G^{\beta\beta}u_{kl}u_{kl}G^{ll}\right)^{1/2}\\*
&=\frac{ {C_1}^2}{4\epsilon\varphi} 
\sqrt{n}  \left( G^{ij}G^{kl}u_{ik}u_{jl}\right).
\end{align*}}}

Now we can estimate \eqref{last line of evolution equation for Z}  from above ---
\begin{align*}
\dfrac {dZ} {dt}&\le G^{ij}D_{ij}Z
+\frac{ {C_1}^2}{4\epsilon\varphi} 
\sqrt{n}  G^{ij}G^{kl}u_{ik}u_{jl}
+
\epsilon \frac{{\varphi_u}^2}\varphi G(Du,Du) \\&\phantom{spacespacespace} 
-2\frac{\varphi_u^2}{\varphi}G(Du,Du)  
 -\frac1{\varphi}G^{ij}G^{ml}u_{mi}u_{lj} -\varphi_t\notag
+\varphi_{uu}G(Du,Du)\\
&=G^{ij}D_{ij}Z
+\frac1\varphi\left(\frac{ {C_1}^2}{4\epsilon} 
\sqrt{n}    -1\right) G^{ij}G^{kl}u_{ik}u_{jl}\\*&\phantom{spacespacespace} 
+\frac{{\varphi_u}^2}\varphi \left(\epsilon-2\right)G(Du,Du) 
-\varphi_t\notag +\varphi_{uu}G(Du,Du).
\end{align*}
The second term is zero if we choose 
$$\epsilon\,=\,{C_1}^2\sqrt{n}/4\,<\,1;$$
the inequality here is a consequence of 
 \eqref{condition on C_1}.

As in the proof of Theorem \ref{interior mcf ch 9},  choose $\varphi=\Phi^{-q}$ for some $q>1$, with $\Phi$ given by \eqref{definition of p for the silly chapter}.  This satisfies the heat equation $\Phi_t=c\Phi''$.  We will choose $c=A$, where $A$ is given by Lemma  \eqref{lemma defining A and P for anisotropic case}.  Substituting $\Phi$ and its derivatives for $\varphi$ and its derivatives (see page \pageref{pagereference for this part}) we find that 
\begin{align*}
\frac{{\varphi_u}^2}\varphi \left(\epsilon\right.-&\left.2\right)G(Du,Du) 
-\varphi_t\notag +\varphi_{uu}G(Du,Du)
\\&\phantom{=} =
q^2\Phi^{-q-2}{\Phi'}^2(\epsilon-2)G(Du,Du)\\&\phantom{space}
+\left[q(q+1)\Phi^{-q-2}{\Phi'}^2-q\Phi^{-q-1}\Phi''\right]G(Du,Du)
        +Aq\Phi^{-q-1}\Phi'' \\
&\phantom{=}=q\Phi^{-q-2}{\Phi'}^2\left[q(\epsilon-1)+1\right]G(Du,Du)+ q\Phi^{-q-1}\Phi''\left[A-G(Du,Du)\right].
\end{align*}

Now, the first term is zero if we choose 
$$q=\frac1{1-\epsilon}=\frac1{1-{C_1}^2\sqrt{n}/4}.$$

As we assumed at the beginning that $\bar F(Du)\ge P$, Lemma \ref{lemma defining A and P for anisotropic case} ensures that 
$$G(Du,Du)\ge A.$$
Consequently, if we choose $T'>0$ small,  $\Phi''$ is positive and then the final term will be negative.   

On the other hand, if we consider the possibility that $\varphi<  P$ at this local maximum, we could replace $\varphi$ by $\sup \lbrace \varphi, P \rbrace$ in the definition of $Z$.    In that case, the first maximum of $Z$ occurs at a point where the barrier is flat, and so the first variation is  
$$0=D_kZ=\left.D\bar F\right|_z(\phi^k)u_{mk},$$
and the evolution equation for $Z$ at the local maximum is
\begin{align*}
\frac{dZ}{dt}&= \left.\bar F D^2\bar F\right|_z(\phi^i,\phi^j)D_{ij}Z-\left. D^2\bar F\right|_{z}(\phi^i,\phi^j) \left.\bar F D^2\bar F\right|_z(\phi^m,\phi^l)u_{mi}u_{lj} \\
&\le 0.
\end{align*}
 
Since $Z_t\le 0$ at the first point where $Z=0$,  $Z\le 0 $ for all $t<T'$ and the conclusion follows.
\halmos

\smallskip {{\noindent\textbf{Proof of Theorem \ref{Second estimate for periodic anisotropic curve...}:}}\quad}
We begin by assuming that at a local maximum point $F(Du)\ge \max\lbrace P,S_{(2/n)^{1/2}}\rbrace$,  and follow the  proof of Theorem \ref{First estimate for periodic anisotropic curvature flows} up to equation \eqref{last line of evolution equation for Z},
  the evolution equation for $Z$ at a local maximum:
\begin{align*}
\dfrac {dZ} {dt}
&=G^{ij}D_{ij}Z+\frac{\varphi_u}{\varphi} \left.\bar F D\left(\bar F D^2\bar F\right) \right|_{z}(\widehat{Du},\widehat\phi^i,\widehat\phi^j)u_{ij}-2\frac{\varphi_u^2}{\varphi}G({Du},{Du})  
\notag \\
&\phantom{spacespacespace}-\frac1\varphi G^{ij}G^{ml}u_{mi}u_{lj} -\varphi_t
+G({Du},{Du})\varphi_{uu}.
\end{align*}
This time we do not choose coordinates to make $G$ the identity.

We use Cauchy-Schwarz (with $\epsilon=1/2$) to estimate the second term --- 
\begin{align*}
\frac{\varphi_u}{\varphi} &\left.\bar F D\left(\bar F D^2\bar F\right) \right|_{z}(Du,\widehat\phi^i,\widehat\phi^j)u_{ij}\\
&\le \frac{{\varphi_u}^2}{2\varphi}G(Du,Du) + \frac1{2\varphi}G(Du,Du)^{-1}
\left[\left.\bar F D\left(\bar F D^2\bar F\right) \right|_{z}(Du,\widehat\phi^i,\widehat\phi^j)u_{ij}\right]^2.
\end{align*}

Choose the basis $\lbrace\phi^1,\dots,\phi^n\rbrace$ at this point so that $D^2u$ is diagonal.  Also, let $I$ be the $n\times n$ diagonal matrix with $+1$ or $-1$ as its diagonal entries, chosen so that  $I^{ii}u_{ii}=|u_{ii}|$.     As $I^2$ is the identity matrix,   $u_{ij}=I^{ik}|u_{kj}|$.   

In these coordinates,
{\allowdisplaybreaks[1]{\begin{align*}
\frac1{2\varphi}G(Du,Du)^{-1}&\left[\left.\bar F D\left(\bar F D^2\bar F\right) \right|_{z}(Du,\widehat\phi^i,\widehat\phi^j)u_{ij}\right]^2\\*
&=\frac1{2\varphi}\frac1{G(Du,Du)}%
\left[\sum_{i=1}^n\left.\bar F D\left(\bar F D^2\bar F\right) \right|_{z}(Du,\widehat\phi^i,\widehat\phi^i)I^{ii}|u_{ii}|\right]^2 
\\&\le \frac1{2\varphi}\frac1{G(Du,Du)}\left[\sum_{i=1}^n\left|\left.\bar F D\left(\bar F D^2\bar F\right) \right|_{z}(Du,\widehat\phi^i,\widehat\phi^i)\right||u_{ii}|\right]^2
\\&\le   \frac1{2\varphi}\frac1{G(Du,Du)}\left[\sum_{i=1}^n\left(\frac2n\right)^{1/2} G(Du,Du)^{1/2}G^{ii}|u_{ii}|\right]^2 \\*
&=\frac1{n\varphi}\left[\sum_{i=1}^n G^{ii}|u_{ii}|\right]^2,
 \end{align*} }}
where the term $(2/n)^{1/2}$ comes from  the estimation of $\left|\left.\bar F D\left(\bar F D^2\bar F\right) \right|_{z}(Du,\widehat\phi^i,\widehat\phi^i)\right|$ in  Lemma \ref{lemma showing that symmetry is as good as small third derivatives}, under the assumption that $\bar F(Du)\ge S_{(2/n)^{1/2}}$.  

This term is dominated by the fourth term of the evolution equation, since (in the same coordinates)
{\allowdisplaybreaks[1]{\begin{align*}
\frac1\varphi G^{ij}G^{ml}u_{mi}u_{lj} 
&=\frac1\varphi G^{ij} G^{ml} I^{m\alpha}|u_{\alpha i}|I^{l\beta} |u_{\beta j}| \\*
&=\frac1\varphi G^{ij}
\left[ I^{\alpha m}G^{ml} I^{l\beta}  \right]
|u_{\alpha i}| |u_{\beta j}|\\
&= \frac1\varphi G^{ij}
G^{\alpha \beta} |u_{\alpha i}| |u_{\beta j}| \\
&\ge \frac1\varphi \frac1n \left[\sum_{i,j}G^{ij}|u_{ji}|\right]^2 
\intertext{where we use the trace inequality $(\trace A)^2\le n\trace (A^2)$}
&= \frac1\varphi \frac1n \left[\sum_{i}G^{ii}|u_{ii}|\right]^2. 
\end{align*}}

What is left of the evolution equation is  
\begin{align*}
\dfrac {dZ} {dt}
&\le G^{ij}D_{ij}Z-\frac32 \frac{{\varphi_u}^2}{\varphi}G(Du,Du) 
-\varphi_t\notag +G({Du},{Du})\varphi_{uu} .
\end{align*}

This is negative at a local maximum if we make the same choice of barrier as before --- $\varphi=\Phi^{-q}$ for $q=2$ with $\Phi$ given by \eqref{definition of p for the silly chapter}.  

If our assumption that $F(Du)\ge \max\lbrace P,S_{(2/n)^{1/2}}\rbrace$ does not hold, then we can replace $\varphi$ by $\max\lbrace P,S_{(2/n)^{1/2}},\varphi\rbrace$.  At the local maximum, $Z_t\le 0$ and so the conclusion follows.
\halmos

\remark In the last theorem, we have chosen $q=2$ somewhat arbitrarily;  in fact $q$ needs only to be strictly greater than $1$, since we can set $q= (1-n\epsilon^2/4)^{-1}$, for  $\epsilon$ given by Lemma \ref{lemma showing that symmetry is as good as small third derivatives}.  However, a smaller $\epsilon$ may entail a larger $S_\epsilon$, so the optimal choice would depend on the exact form of $\bar F$.

\subsection*{Interior estimate for anisotropic mean curvature flow}

We begin by showing that when we have the symmetry condition, an estimate analogous to \eqref{estimate for nabla u nabla eta} in the isotropic case is possible.  
  
\begin{lemma} Suppose that $\bar F$ satisfies the symmetry condition \eqref{symmetry condition}.
Then  there exists a constant $C_2$ depending only on $\bar F$ such that 
\begin{equation*}\left.\bar F D^2\bar F\right|_{p-\phi^0}(p,q)\le C_2\frac {\bar F(q)}{\bar F(p-\phi^0)}\end{equation*}
for all $p=p_i\phi^i$ and  $q=q_i\phi^i$.     \label{cross-terms lemma}
\end{lemma}
\proof   This estimate is unchanged under $q\mapsto sq$, so  we need only to show that this holds for $q$ on the unit ball.  Let $q=q_i\phi^i$ be a fixed point on the unit ball.

By compactness, the estimate holds for all $p$ on the unit ball.  

Let $ p$ be a fixed point on the unit ball and consider%
\begin{align*}
\lim_{s\rightarrow \infty}  
\bar F(s p-\phi^0)&\left.\bar F D^2\bar F\right|_{s p-\phi^0}(s p,q)\\
&=\lim_{s\rightarrow 0}  \bar F(p-s\phi^0)\frac1s {\left.\bar F D^2\bar F\right|_{p-s\phi^0}(\phi^0,q)}\\
&=\lim_{s\rightarrow 0}  \bar F(p-s\phi^0)\frac1s \left[\left.\bar F D^2\bar F\right|_{p-s\phi^0}(\phi^0,q)-\left.\bar F D^2\bar F\right|_{p}(\phi^0,q)\right]
\intertext{where we have added zero in the form of $\left.\bar F D^2\bar F\right|_{p}(\phi^0,q)$}
&= \bar F( p)\left[-D\left(\left.\bar F D^2\bar F\right)\right|_{ p}(\phi^0,\phi^0,q)\right].
\end{align*}

So, either the supremum of $\bar F(s p-\phi^0)\left.\bar F D^2\bar F\right|_{s p-\phi^0}(s p,q)$ over $s\in[0,\infty)$   is the limit above, or else it is attained at some finite $s$.   In either case, 
$$C_p(q)=\sup_{s\ge 0} \bar F(sp-\phi^0) \left.\bar F D^2\bar F\right|_{sp-\phi^0}(sp,q)$$
is finite, and we can set 
$$C_2=\sup_{\lbrace q:\bar F(q)=1\rbrace} \, \sup_{\lbrace p: \bar F( p)=1\rbrace } C_{p}(q) $$
to complete the lemma.
\halmos

The next lemma shows that the trace of $\bar F D^2\bar F$ is bounded below. %

\begin{lemma} \label{lemma giving lower bound on trace G}
Let $\lbrace \phi^0,\phi^1,\dots,\phi^n\rbrace$ be a basis for $V^*$, where $n>1$.  Then there is a constant $k>0$ so that for all $p=\sum_{i=1}^n p_i\phi^i$,  
\begin{equation*}
\sum_{i=1}^n \left.G\right|_{p-\phi^0}(\phi^i,\phi^i)\ge k,
\end{equation*}
where $G=\bar F D^2\bar F$.  
\end{lemma}
\proof  
By compactness and strict convexity of the unit ball, 
\begin{equation*}
\sup_{\lbrace p:\bar F(p)\le 1\rbrace} \left.G\right|_{p-\phi^0}(\phi^i,\phi^i)=\sup_{\lbrace p:\bar F(p)\le 1\rbrace} \left.G\right|_{p-\phi^0}(\widehat\phi^i,\widehat\phi^i)>0,
\end{equation*}
since $\widehat \phi^i$ is a non-zero tangent covector.
 
Now consider 
\begin{align*}
\lim_{t\rightarrow\infty} \sum_{i=1}^n \left.G\right|_{tp-\phi^0}(\phi^i,\phi^i)=
 \sum_{i=1}^n \left.G\right|_{p}(\widehat\phi^i,\widehat\phi^i), \end{align*}
where, in the limit, $\widehat \phi^i=\phi^i- \frac{\left.D\bar F\right|p(\phi^i)}{\bar F(p)}p$.  

At most one of the $\phi^i$  may be parallel to $p$ --- suppose it is  $\phi^1$, in which case $\lim_{t\rightarrow\infty} \left.G\right|_{tp-\phi^0} (\phi^1,\phi^1)=0$.  For the remaining $(n-1)$ basis covectors, $\widehat\phi^i$ are non-zero tangent covectors (to the unit ball at $p$) and so  $\left.G\right.(\widehat\phi^i,\widehat\phi^i)$ is again bounded below for $i\not=1$.  

It follows that 
\begin{equation*}
\lim_{t\rightarrow\infty} \sum_{i=1}^n \left.G\right|_{tp-\phi^0}(\phi^i,\phi^i)\ge
 \sum_{i=2,n} \left.G\right|_{p}(\widehat\phi^i,\widehat\phi^i) >0.
\end{equation*}

If we take the infimum of all such lower bounds, over all $p$ in the unit ball, then the conclusion follows.
\halmos

Let $u:B_R(0)\times[0,T]\rightarrow \bigR$ be a $H_2$, bounded 
$$|u(x,t)|\le M$$
solution on the ball of radius $R$ to the anisotropic curvature flow equation \eqref{anisotropic curvature flows}, where $\bar F$ is a positive, convex  homogeneous degree one function, with a strictly convex unit ball.

\begin{theorem}[Interior estimate for anisotropic mean curvature flow] \label{interior estimate for anisotropic curve...}
If $\bar F$ satisfies both the smallness of third derivatives condition \eqref{condition on D^3 F} with some constant $$C_1^2<2/\sqrt n,$$ and the symmetry condition \eqref{symmetry condition},
then %
 \begin{equation*}
F(Du)\le t^{q/2}\exp\left(\frac{Aq(|u|-2M)^2}{4t}\right)\left(R^2-2kt-|x|^2\right)^{-r}
\end{equation*}
  for \phantom{.}$0<t\le T'$.  

Here,  $A$ is given by  Lemma \ref{lemma defining A and P for anisotropic case} and depends on $\bar F$; $k$ is given by Lemma  \ref{lemma giving lower bound on trace G} and depends on $\bar F$;  and $T'>0$, $q>1$, and $r>1$  depend on $M$, $A$, $P$ (which is also given by  Lemma \ref{lemma defining A and P for anisotropic case}) and $k$. 
\end{theorem}

\proof  We introduce the localising term $\eta$ into our definition of $Z$, now restricted to the shrinking ball:
\begin{equation*}
Z:=F(Du)-\frac\varphi\eta  
\end{equation*}
for $(x,t)\in B_{\sqrt{R^2-2kt}}\times[0,T]$, 
where $k$ is the constant given by Lemma \ref{lemma giving lower bound on trace G}, $\varphi=\varphi(u,t)$ is a smooth strictly positive function chosen so that $Z<0$ at the initial time, and $\eta$ is a smooth positive function chosen so that $\eta\rightarrow 0$ on the boundary of the shrinking ball.

Assume that at the first interior point where $Z=0$, $F(Du)\ge P$.

Then $F(Du)=\varphi/\eta$ and as this is a spatial maximum (since the choice of $\eta$ ensures that there are no boundary maxima)
we have a first derivative condition
\begin{equation}  \label{first derivative condition for interior anisotropic estimate}
0=D_kZ=\left.D\bar F\right.|_z(\phi^m)u_{mk}-D_k\left(\varphi/\eta\right).
\end{equation}

An evolution equation for  $\varphi/\eta$ is:
\begin{align*}
\dfrac{d}{dt} \left(\frac \varphi \eta\right)
&=\frac1\eta\left(\varphi_u u_t+\varphi_t\right) -\frac\varphi{\eta^2}\frac{d\eta}{dt}+ \left.\bar F D^2\bar F\right|_{z}(\phi^i,\phi^j)D_{ij}\left(\frac \varphi \eta\right)
\\&\phantom{spacespace}
-\left.\bar F D^2\bar F\right|_{z}(\phi^i,\phi^j)\bigg[\frac1\eta\left(\varphi_{uu}u_iu_j+\varphi_uu_{ij}\right)-\frac{\varphi_u}{\eta^2}\left(u_jD^i\eta+u_iD^j\eta\right)
\\&\phantom{spacespacespacespacespacespace}
+2\frac\varphi{\eta^3}D^i\eta D^j\eta-\frac{\varphi}{\eta^2}D_{ij}\eta \bigg]\\
&=G^{ij}D_{ij}\left(\frac\varphi\eta\right)+\frac1\eta\left[\varphi_t-G(Du,Du)\varphi_{uu} \right]
-\frac\varphi{\eta^2}\left(\frac{d}{dt}- G^{ij}D_{ij}\right)\eta 
\\&\phantom{spacespace}
+2\frac{\varphi_u}{\eta^2}G(Du,D\eta) -2\frac\varphi{\eta^3}G(D\eta,D\eta),
\end{align*}  
where in the last line we have used the notation $G^{ij}=\left.\bar F D^2\bar F\right|_{z}(\phi^i,\phi^j)$.

We can incorporate the first derivative condition \eqref{first derivative condition for interior anisotropic estimate} into \eqref{evolution eqn for F, first time},  the evolution equation for $ F$: {\allowdisplaybreaks{
\begin{align*}\dfrac {dF}{d t}&=\left.\bar F D^2\bar F\right|_{z}(\phi^i,\phi^j)D_{ij}F +
\left.D\bar F\right|_z(\phi^k)  
\left.D\left(\bar F D^2\bar F\right) \right|_z(\phi^m,\widehat\phi^i,\widehat\phi^j)u_{mk}u_{ij} 
\notag \\* &\phantom{===}
-\left.D\bar F\right|_z(\phi^k) \left[ \left.D\bar F\right|_z(\phi^i)  \left. D^2\bar F \right|_z(\phi^m,\widehat\phi^j) +  \left.D\bar F\right|_z(\phi^j)  \left. D^2\bar F \right|_z(\phi^m,\widehat\phi^i) \right]u_{mk}u_{ij} 
\\*&\phantom{spacespacespacespacespace}
-\left.\bar F D^2\bar F\right|_{z}(\phi^i,\phi^j) \left.D^2\bar F\right|_z(\phi^m,\phi^l)u_{mi}u_{lj} 
 \\ 
&= G^{ij}D_{ij}F  + 
D^m\left(\varphi/\eta\right)  
\left.D\left(\bar F D^2\bar F\right) \right|_z(\phi^m,\widehat\phi^i,\widehat\phi^j)u_{ij} 
\notag 
\\*&\phantom{spacespacespacespacespace}
-2\frac{\eta}{\varphi}G\left(D(\varphi/\eta),D(\varphi/\eta)\right)
-\frac{\eta}{\varphi}G^{ij}G^{ml}u_{mi}u_{lj}
 \\
&=G^{ij}D_{ij}F  + \frac{\eta}{\varphi} \left.\bar F
D\left(\bar F D^2\bar F\right) \right|_z\left(%
D(\varphi/\eta),\widehat\phi^i,\widehat\phi^j\right)u_{ij} \\*
&\phantom{===}-2\frac{\eta}{\varphi}\left[\frac{\varphi_u^2}{\eta^2}G(Du,Du)-2\frac{\varphi\varphi_u}{\eta^3}G(Du,D\eta)+\frac{\varphi^2}{\eta^4}G\left(D\eta,D\eta\right)\right] 
\notag \\*  &\phantom{spacespacespacespacespace}
-\frac{\eta}{\varphi}G^{ij}G^{ml}u_{mi}u_{lj}.
\end{align*} }}

Putting the last two steps together gives an evolution equation for $Z$ at a local maximum:  
\begin{align*}\dfrac {dZ}{d t} &= \dfrac {dF}{d t} -\dfrac{d }{dt}\left(\frac\varphi\eta\right) \\
&=G^{ij}D_{ij}Z +\frac{\eta}{\varphi} 
\left.\bar FD\left(\bar F D^2\bar F\right) \right|_z\left(D(\varphi/\eta),\widehat\phi^i,\widehat\phi^j\right)u_{ij} 
-\frac{\eta}{\varphi}G^{ij}G^{ml}u_{mi}u_{lj}
\\&\phantom{spacespacespace} 
-\frac1\eta\left[\varphi_t-G(Du,Du)\varphi_{uu} +2\frac{\varphi_u^2}{\varphi}G(Du,Du)\right] 
\\&\phantom{spacespacespace}+\frac\varphi{\eta^2}\left(\frac{d}{dt}-G^{ij}D_{ij}\right)\eta
+2\frac{\varphi_u}{\eta^2}G(Du,D\eta).
\end{align*}

The second term here may be split up into a part with $D\varphi$ and a part with $D\eta$:
\begin{align*}
\frac{\eta}{\varphi} &\left.\bar F
D\left(\bar F D^2\bar F\right) \right|_z \left(\frac{\varphi_u}\eta Du-\frac{\varphi}{\eta^2} D\eta,\widehat\phi^i,\widehat\phi^j\right)u_{ij} \\
&=\frac{\varphi_u}{\varphi} \left.\bar F
D\left(\bar F D^2\bar F\right) \right|_z 
        \left(Du,\widehat\phi^i,\widehat\phi^j\right)u_{ij}  
-\frac{1}{\eta} \left.\bar F D\left(\bar F D^2\bar F\right) \right|_z 
        \left(D\eta,\widehat\phi^i,\widehat\phi^j\right)u_{ij}. 
\end{align*}

These may be individually estimated using the Cauchy-Schwarz inequality and the smallness-of-third-derivatives condition, %
as described in the proof of Theorem \ref{First estimate for periodic anisotropic curvature flows} on page \pageref{page for C-s} ---
{\allowdisplaybreaks{\begin{gather*}
\begin{split}
\frac{\varphi_u}{\varphi} \left.\bar F
D\left(\bar F D^2\bar F\right) \right|_z 
        &\left(Du,\widehat\phi^i,\widehat\phi^j\right)u_{ij}  \\*
&\le \mu_1 \frac{\varphi_u^2}{\varphi\eta} G(Du,Du)  + \frac1{4\mu_1}\frac\eta{\varphi} C_1^2\sqrt{n} \left(G^{ij}G^{ml}u_{mi}u_{lj}\right),
\end{split}\\
\begin{split}
-\frac{1}{\eta} \left.\bar F D\left(\bar F D^2\bar F\right) \right|_z& 
\left(D\eta,\widehat\phi^i,\widehat\phi^j\right)u_{ij}\\*
&\le \mu_2\frac \varphi {\eta^3}G(D\eta,D\eta)+\frac1{4\mu_2}\frac\eta{\varphi}C_1^2\sqrt n\left(G^{ij}G^{kl}u_{ik}u_{jl}\right),
\end{split}\end{gather*} }}
for some $0<\mu_1,\mu_2<1$.

We choose the localising term to be $\eta:=\tilde \eta^r$ for some $r>1$,  $\tilde \eta=R^2-2kt-|x|^2$, and $k$ given in Lemma \ref{lemma giving lower bound on trace G}.  Then
\begin{gather*}
D_i\eta= r\tilde \eta^{r-1}D_i\tilde\eta\\
D_{ij}\eta=r\tilde \eta^{r-1}D_{ij}\tilde \eta + r(r-1)\tilde \eta^{r-2}D_i\tilde \eta D_j\tilde \eta,
\end{gather*}
and  the second-last term of the evolution equation is 
\begin{align*}
\frac\varphi{\eta^2}\left(\frac{d}{dt}-G^{ij}D_{ij}\right)\eta&=
\frac\varphi{\eta^2}r{\tilde\eta}^{r-1}\left[2k-2\trace G -(r-1){\tilde\eta}^{-1}G(D\tilde\eta, D\tilde\eta) \right]\\
&\le \frac\varphi{\eta^2}r{\tilde\eta}^{r-2}(1-r)G(D\tilde\eta, D\tilde\eta). \end{align*}

As $\bar F$ satisfies the symmetry condition \eqref{symmetry condition}, we may use  Lemma \ref{cross-terms lemma} to estimate the
 final term of the evolution equation:
\begin{align*}
2\frac{\varphi_u}{\eta^2}G(Du,D\eta)&=2\frac{\varphi_u}{\eta^2}\left.\bar FD^2\bar F\right|_{Du-\phi^0}(Du,D\eta)
\\&\le 2\frac{\varphi_u}{\eta^2} \frac{C_2\bar F(D\eta)}{\bar F(Du-\phi^0)}\\
&= 2C_2\bar F(D\eta)\frac{\varphi_u}{\varphi\eta}.
\end{align*}

The evolution equation can now be estimated from above --- 
\begin{align}\dfrac {dZ}{d t}
&\le G^{ij}D_{ij}Z   + \frac{\eta}\varphi\left(\frac1{4\mu_1} C_1^2\sqrt n+\frac1{4\mu_2}C_1^2\sqrt n -1 \right) G^{ij}G^{ml}u_{mi}u_{lj} 
\notag \\&\phantom{space}
-\frac1\eta\left[\varphi_t-G(Du,Du)\varphi_{uu} +(2-\mu_1)\frac{\varphi_u^2}{\varphi}G(Du,Du)
-2C_2\bar F(D\eta)\frac{\varphi_u}\varphi\right]
\notag \\&\phantom{spacespacespace}+\frac\varphi{\eta^2}r\tilde\eta^{r-2} \left(1-r+r\mu_2\right)G(D\tilde \eta,D\tilde \eta). \label{almost the last eqn}
\end{align}

Since ${C_1^2\sqrt n}/4<1/2$, we can choose $\mu_1<1$ and $\mu_2<1$ %
such that 
\begin{gather*}
\frac{C_1^2\sqrt{n}}{4}\left(\frac1\mu_1+\frac1\mu_2\right)\le 1. %
\end{gather*}

With such choices, the second term of the evolution inequality \eqref{almost the last eqn} will be negative.  We can also  set $r=(1-\mu_2)^{-1}>1$,  so the coefficient of $\tilde\eta^{-1}G(D\tilde \eta,D\tilde \eta)$ is zero.

As in the previous cases we can set $\varphi=\Phi^{-q}$ where $\Phi$ is given by \eqref{definition of p for the silly chapter} with $c=A$, where $A$ is given by Lemma \ref{lemma defining A and P for anisotropic case}.

The bracketted part  of the second line of the evolution equation is then
\begin{align}
&\varphi_t-G(Du,Du)\varphi_{uu} +(2-\mu_1)\frac{\varphi_u^2}{\varphi}G(Du,Du)
-2C_2\bar F(D\eta)\frac{\varphi_u}\varphi \notag\\
&\phantom{space}= -q\Phi^{-q-1}\left(\Phi_t-G(Du,Du)\Phi''\right)  \notag \\
&\phantom{spacespace}+ 
G(Du,Du)q{\Phi'}^2\Phi^{-q-2}
\left[-1+(1-\mu_1)q-2\frac{C_2\bar F(D\eta)}{G(Du,Du)}\frac{\Phi^{q+1}}{|\Phi'|}\right] \label{final eqn from newwork chap}
\end{align} 

If we choose $T'$ small enough that $\Phi''\ge0$, then the term $\Phi_t-G(Du,Du)\Phi''=\left(A-G(Du,Du)\right)\Phi''$ is negative.

As $r>1$, $\bar F(D\eta)= r\tilde\eta^{r-1}\bar F(D\tilde\eta)\le C_3rR^{2r-1}$, where $C_3>0$ depends only on $\bar F$.

With this choice of $\Phi$,
\begin{align*}
2\frac{C_2\bar F(D\eta)}{G(Du,Du)}\frac{\Phi^{q+1}}{|\Phi'|}&\le 
2\frac{C_2C_3rR^{2r-1}}{A}\frac{\Phi^{q+1}}{|\Phi'|} \\
&\le 4 \frac{C_2C_3rR^{2r-1} {\Phi^{q}t  }}{A^2M} \\
&\le 4 \frac{C_2C_3rR^{2r-1} {t  }}{A^2M} 
\end{align*}
if we choose $T'$ small enough that $\Phi\le 1$.  

In order to ensure that the last  part of \eqref{final eqn from newwork chap} is positive, we choose $q$ so that 
$$q\ge \frac1{1-\mu_1}\left(1+ 4 \frac{C_2C_3rR^{2r-1} {T' }}{A^2M} \right).$$

So, at such maxima, $Z_t\le0$.  

At local maxima where $F(Du)<P$, then in the definition of $Z$ we replace $\varphi/\eta$ by $\max\lbrace \varphi/\eta, P\rbrace$, in which case the barrier is flat at the local maxima, and we again find that 
$Z_t\le0$.

In either case, the maximum principle ensures that $Z$ is never greater than zero and the conclusion follows. 
\halmos

\appendix
\renewcommand{\chaptername}{Appendix}
\newchapter{Function spaces and regularity estimates for parabolic equations}{Function spaces and regularity estimates}

Here, we define relevant function spaces, and survey some regularity results used in the existence theorems of Chapters \ref{short-time chapter}, \ref{neumann chapter} and \ref{dirichlet chapter}.  This treatment follows the books of Krylov \cite{krylov:parabolic} and Lieberman \cite{li:parabolic}. 

\section{Function spaces}

On the space of continuous functions $u:\Omega\rightarrow \mathbb{R}$,  we have the supremum norm
$$\vert u\vert _{0;\Omega}:=\sup_{x\in\Omega} |u(x)|.$$

Define the H\"older semi-norm with exponent $\alpha\in(0,1]$ by 
$$[u]_{\alpha;\Omega}:=\sup_{\substack{ x,y\in \Omega \\ x\not= y}}\frac{|u(x)-u(y)|}{|x-y|^\alpha}.$$
For an integer $k\ge 0$, we define the H\"older $(k+\alpha)$-norm by 
$$\vert u\vert_{k+\alpha;\Omega}:=\sum_{|\beta|\le k}\vert D^\beta u\vert_{0;\Omega} + \sum_{|\beta| =k} [D^\beta u]_{\alpha;\Omega},$$
where $\beta$ is a multi-index --- an $n$-tuple of non-negative integers with $|\beta|=\sum \beta_i$, and where $D^\beta u:= \dfrac{\partial^{|\beta|} u}{\partial ({x^1})^{\beta_1}\dots \partial( {x^n})^{\beta_n}}$.
   
The Banach space associated with this norm is $C^{k+\alpha}(\Omega)=\{u:\vert u\vert_{k+\alpha}<\infty\}$.

\subsection*{ The parabolic H\"older spaces $H_\alpha$}

With parabolic equations, it is useful to weight the space and time variables differently- that is, two space derivatives to one time derivative.    Following Lieberman, we will denote parabolic H\"older spaces by $H$ rather than $C$.  
\nomenclature[H4alpha]{$H_\alpha $   }{a parabolic H\"older space\refpage }
For points $z=(x,t)$ in a domain $Q\subset \mathbb R^n\times \mathbb R$,  define 
the parabolic H\"older semi-norm by %
$$[u]_{\alpha, \alpha /2 ; Q}:=\sup_{\substack{z_1,z_2\in Q \\ z_1\not= z_2}}\frac{|u(z_1)-u(z_2)|}{(|x_1-x_2|+|t_1-t_2|^{\frac12})^\alpha},$$
where $\alpha\in(0,1]$%
.  The parabolic norm is
$$\vert u\vert_{\alpha,\alpha/2;Q}:= |u|_{0;Q} + [u]_{\alpha,\frac\alpha2;Q}.$$

Higher spatial derivatives and derivatives in time are bounded by $k+\alpha$ norms, where $k\ge 0$ is an integer:
\begin{equation*}
|u|_{k+\alpha,\alpha/2;Q}=\sum_{|\beta|+2j\le k} \sup |D^\beta_x D^j_t u| + \sum_{|\beta|+2j=k}[D^\beta_xD^j_t u]_{\alpha,\alpha /2} .
\end{equation*}

The Banach space associated with the $|\cdot|_{k+\alpha, \alpha/ 2 }$ norm is $$H_{k+\alpha}(Q)=\{u:|u|_{k+\alpha, \alpha /2 } <\infty\}.$$

When the region $Q$ is a cylinder, in the sense that $Q=\Omega\times[0,T]$, the parabolic boundary is given by 
\begin{equation*}
{\cal P}\left(\Omega\times[0,T]\right)=\Omega\times\lbrace 0\rbrace \cup \partial\Omega\times[0,T].
\end{equation*}
On the boundary, we can define parabolic norms $H_{k+\alpha}({\cal P})$ exactly as above.

\section{Regularity estimates} 
In the following, $\boldsymbol{P}$ is a quasilinear parabolic operator 
$$\boldsymbol{P}u=-u_t+a^{ij}(Du,x,t)D_{ij}u+b(Du,u,x,t)$$ 
with positive constants $\lambda_K$ and $\Lambda_K$ such that 
\begin{gather*}
 a^{ij}(x,p)\xi_i\xi_j\geq \lambda_K|\xi|^2 \text{ for }\xi\in\bigR^n \\
 |a^{ij}(x,p)|\leq \Lambda_K,
\label{parabolicity for a} 
\end{gather*}
whenever $|p|\le K$.   We work on a domain $\Omega\times[0,T]$ for some smoothly bounded $\Omega\subseteq \bigR^n$.

Here $Q_r$ is the intersection of the region and a cylinder:
$$Q_R(x_0,t_0):=\{(x,t)\in\Omega\times[0,T]:|x-x_0|<R, t_0-R^2<t<t_0\}.$$ 

We begin with an oscillation estimate for the gradient of a solution for a Neumann problem near a flat boundary:
\begin{theorem}%
\label{Lieb}  
Suppose that $a^{ij}=a^{ij}(Du,x)$, $b=b(x)$, and that inside $Q_R$, the boundary of $\Omega$ is  $x^n=0$.   Let $u\in C^{2,1}\cap H_1(Q_R)$  be a solution of  $\boldsymbol{P}u=0$ when $x^n>0$ and $D_nu=0$, when $x^n=0$.  

Suppose there are positive constants $b_0$ and $\lambda_0$ 
such that
\begin{equation*}
 |{a^{ij}}_p|\leq \lambda_0,  \quad |b(x)|\leq b_0 \end{equation*}
for all $(x,p)\in Q_R\times {\bigR}^n$ with $|p|\le K$.

If $|Du|\le K$, then there are positive constants $\theta$ and $\sigma$ determined only by $K$, $n$, $\lambda$, $\Lambda$ and $\lambda_0$ such that 
$$\osc_{\Omega \cap Q_R} a^{ij}(\cdot,p)\le\sigma \text{ for all }|p|\le K$$ 
implies
$$\osc_{\Omega \cap Q_r} Du \le C(K,n,\lambda, \Lambda,\lambda_0)\left(\frac rR\right)^\theta [\osc_{\Omega \cap Q_R} Du + b_0 R].$$
\end{theorem}

\smallskip

Similarly, we can find an estimate near the boundary for problems with Dirichlet boundary conditions:

\begin{theorem}
\label{side-estimate}
Suppose that  $a^{ij}$ and $b$ are uniformly continuous,  that $a^{ij}$ is differentiable with respect to $(p,q,x)$,  %
and where, if $|q|+|p|\le K$, we can find a positive constant %
$\mu_K$ such that
\begin{gather}
K[|{a^{ij}}_x|+|{a^{ij}}_q||p|] + |b|\le \mu_K.   \label{equation for mu_k}
\end{gather}

If $u\in C^{2,1}\cap H_1(\Omega\times[0,T])$ satisfies $\boldsymbol{P}u=0$ on $\Omega\times[0,T]$ and $u=0$ on $\partial\Omega\times(0,T)$, 
then there are positive constants $C$ and $\alpha$ depending on $n, \lambda_K, \Lambda_K$ and ${\mu_K R}/K$ such that for any $(x_0,t_0)\in\partial\Omega\times(0,T)$
\begin{equation*}
\osc_{Q_r(x_0,t_0)} Du \le 
        C\left[ \left(\frac rR\right)^\alpha K+%
                + \mu_K r\right]
\end{equation*} 
for $0<r<R<t_0^{1/2}$.
\end{theorem}

On the interior of the domain, one can also find a H\"older bound for the gradient:

\begin{theorem}
\label{interior holder gradient estimate} 
Suppose that $a^{ij}$ and $a$ are continuous;  $a^{ij}$ is differentiable with respect to $(p,q,x,t)$; and where  if $|q|+|p|\le K$,  we can find a positive constant $\mu_K$ satisfying \eqref{equation for mu_k}.  

If $u\in C^{2,1}$ satisfies $\boldsymbol{P}u=0$  and $|Du|+|u|\le K$ in $\Omega\times[0,T]$,  
then there is a positive constant $\alpha$ determined by $n,\lambda_K, \Lambda_K, \sup|{a^{ij}}_p|K,$ and $\mu_KR/K$ such that for interior sets $\Omega'\times[t_1,t_2]\subset\subset \Omega\times[0,T]$ we have 
$$[Du]_{\alpha;\Omega'\times[t_1,t_2]}\le C\left(n,K,\lambda_K,\Lambda_K,\mu_K,{\diam}(\Omega\times[0,T])\right)d^{-\alpha}, $$
where $d=\dist(\Omega'\times[t_1,t_2],{\cal P} \Omega\times[0,T])$.  

If $Q_r$ is a cylinder in the interior of the domain, we also have 
$$\osc_{Q_r(x_0,t_0)}Du\le C(n,K,\lambda_K,\Lambda_K)\left(\frac r R\right)^\alpha
\left(\osc_{Q_R(x_0,t_0)}Du + \mu_KR\right)$$
as long as $0<r\le R \le d(x_0,t_0)$, where $d(x_0,t_0)$ is the distance from $(x_0,t_0)$ to the parabolic boundary ${\cal P}(\Omega\times[0,t_0])$.

\end{theorem}

The following interior estimate is a $H_{2+\alpha}$ bound for $u$ when the coefficients of $\boldsymbol{P}$ are smooth:
\begin{theorem}
\label{interior1}%
Suppose that $a^{ij}=a^{ij}(x,t)$, and $b(Du,u,x,t)=b^i(x,t)D_i u + c(x,t)u$, and that there is a constant $K$ such that
$$|a,b,c|_{\alpha,\alpha/2}\le K,$$ 
for $\alpha\in(0,1)$.  

Then for any $R>0$ there is a constant $C$ dependent on $R$, $\lambda$, $\alpha$, $K$ and $n$ such that if $u\in H_{2+\alpha}(Q_{3R})$, %
then 
\begin{equation*}
|u|_{2+\alpha,1+\alpha/2;Q_R}\le C \left(|\boldsymbol{P}u|_{\alpha,\alpha/2;Q_{2R}} + |u|_{0;Q_{2R}}\right).
\end{equation*}
\end{theorem}

In a similar vein, there are higher regularity estimates on the interior of a domain:
\begin{theorem}
\label{regularity1}
If we have $|D^\alpha a|_{\alpha ,\alpha/2; Q_{2R}}$, $|D^\alpha b|_{\alpha ,\alpha/2; Q_{2R}}$, $|D^\alpha c|_{\alpha ,\alpha/2; Q_{2R}} \le K$ for any $ |\alpha|\le k$, and if $u\in H_{ 2+\alpha}(Q_{2R})$ and $D^\alpha (\boldsymbol{P}u)\in H_{ \alpha }(Q_{2R})$, then  $D^\alpha u\in H_{\alpha}(Q_{R})$ for $|\alpha|\le k$, and there is a constant $C=C(R,K,\alpha,k,d,\lambda)$ such that
$$\sum_{|\beta|\le k}|D^\beta u|_{2+\alpha,1+\alpha/2;Q_{R}}
\le  \sum_{|\beta|\le k}|D^\beta|_{2+\alpha  ,1+\alpha/2;Q_{2R}}
+ C|u|_{0;Q_{2R}}.$$
\end{theorem}

\markright{\color{draftgrey}{Notation}} 
\renewcommand{\nomname}{Notation} %

\addcontentsline{toc}{chapter}{Notation}

\printglossary[1.8cm]

\nomenclature[Q]{$Q$}{the {Cartan tensor} restricted to the tangent space of the unit ball of $\bar F$\refpage}
\nomenclature[a]{$a$}{coefficient of second derivative in a parabolic operator; also appears as $a^{ij}$} 

\nomenclature[b]{$b$}{ lower order terms in a parabolic operator  }

\nomenclature[Qr]{$Q_R(x_0,t_0)$}{{the parabolic cylinder of radius $R$ centred at $(x_0,t_0)$, $Q_R(x_0,t_0):=\lbrace(x,t):\vert x-x_0\vert <R, t_0-R^2<t<t_0\rbrace$}}

\nomenclature[E]{$\cal E$}{ the Bernstein $\cal E$-function, ${\cal E} (p,q,x,t)=a^{ij}(p,q,x,t)p^ip^j $}    
\nomenclature[Erf]{{${\mathcal{E}} {\rm{rf}}$}}{{the error function $\erf{x}:=\frac 2 {\sqrt{\pi}}\int_0^x e^{-y^2}dy $}}

\nomenclature[F]{$\mathbf{F}$}{embedding of a manifold}

\nomenclature[H3nu]{$\mathbf{H}$}{mean curvature vector  $\mathbf{H}=H\nu$}

\nomenclature[B]{$B_R(x_0)$ }{the open ball of radius $R$ centred at $x_0$, {$B_R(x_0)=\lbrace x:\vert   x-x_0\vert <R\rbrace$ } }

\nomenclature[P]{$\cal P$}{ the parabolic boundary of a set:  ${\cal P}\left(\Omega\times[0,T]\right)=\Omega\times\lbrace 0\rbrace \cup \partial\Omega\times[0,T]$. } 

\nomenclature[P1]{$\boldsymbol{P}$}{a parabolic operator, page \pageref{page where parabolic operator is defined}} 
\nomenclature[sn]{$\bigS^n $}{  the set of symmetric $n\times n$ matrices }

\nomenclature[la]{$\lambda_{\cal K}$}{ lower parabolicity constant on a bounded set $\cal K$ }

\nomenclature[Lb]{$\Lambda_{\cal K}$ }{ upper parabolicity constant on a bounded set $\cal K$ }

\nomenclature[lc]{$\lambda(K)$}{ lower parabolicity constant on a set with $\vert p\vert< K$, page \pageref{page ref for K} }
\nomenclature[Ld]{$\Lambda(K)$}{ upper parabolicity constant on a set with $\vert p\vert< K$, page \pageref{page ref for K} }

\nomenclature[mij]{$m^{ij}$}{the second-order coefficients of the equation for mean curvature flow of graphs, $m^{ij}(p)=\delta^{ij}-\delta^{ik}\delta^{jl}p_kp_l(1+\vert p\vert^2)^{-1}$}

\nomenclature[nu]{$\nu$       }{ outwards unit normal --- to the surface, page \pageref{nu}; to the boundary of the region, page \pageref{page where nu is second time around}  }

\nomenclature[phi]{$\Phi$       }{  the fundamental solution to the heat equation, $\Phi(y,t)=t^{-1/2}\exp\left(-cy^2/t\right)$ } 

\nomenclature[psi]{$\psi    $   }{  a special  barrier defined in Section \ref{section where psi is defined}, page \pageref{section where psi is defined}}

\nomenclature[Psi]{$\Psi$ }{ a boundary-straightening change of coordinates }

\nomenclature[Omega]{$\Omega^\epsilon$}{an interior set:  given $\Omega$,  $\Omega^\epsilon=\lbrace x\in\Omega: \dist(x,\partial\Omega)>\epsilon\rbrace$}

\markright{\color{draftgrey}{Bibliography}} 

\bibliographystyle{plain}
\bibliography{../latex-sources/2000}

\begin{thebibliography}{10}

\bibitem{almgren:1993}
Fred Almgren, Jean~E. Taylor, and Lihe Wang.
\newblock Curvature-driven flows: a variational approach.
\newblock {\em SIAM J. Control Optim.}, 31(2):387--438, 1993.

\bibitem{aag:rotation}
Steven Altschuler, Sigurd~B. Angenent, and Yoshikazu Giga.
\newblock Mean curvature flow through singularities for surfaces of rotation.
\newblock {\em J. Geom. Anal.}, 5(3):293--358, 1995.

\bibitem{andrews:convexcurves}
Ben Andrews.
\newblock Evolving convex curves.
\newblock {\em Calc. Var. Partial Differential Equations}, 7(4):315--371, 1998.

\bibitem{andrews:anisotropic}
Ben Andrews.
\newblock Volume-preserving anisotropic mean curvature flow.
\newblock {\em Indiana Univ. Math. J.}, 50(2):783--827, 2001.

\bibitem{ang:zeroset}
Sigurd Angenent.
\newblock The zero set of a solution of a parabolic equation.
\newblock {\em J. Reine Angew. Math.}, 390:79--96, 1988.

\bibitem{angenent:nodal}
Sigurd Angenent.
\newblock Nodal properties of solutions of parabolic equations.
\newblock {\em Rocky Mountain J. Math.}, 21(2):585--592, 1991.
\newblock Current directions in nonlinear partial differential equations
  (Provo, UT, 1987).

\bibitem{angenent:formation}
Sigurd Angenent.
\newblock On the formation of singularities in the curve shortening flow.
\newblock {\em J. Differential Geom.}, 33(3):601--633, 1991.

\bibitem{angenent:curves-II}
Sigurd Angenent.
\newblock Parabolic equations for curves on surfaces. {II}. {I}ntersections,
  blow-up and generalized solutions.
\newblock {\em Ann. of Math. (2)}, 133(1):171--215, 1991.

\bibitem{angenent-gurtin:aniso}
Sigurd~B. Angenent and Morton~E. Gurtin.
\newblock Anisotropic motion of a phase interface. {W}ell-posedness of the
  initial value problem and qualitative properties of the interface.
\newblock {\em J. Reine Angew. Math.}, 446:1--47, 1994.

\bibitem{chern:finsler}
David Bao, Shiing-Shen Chern, and Zhongmin Shen.
\newblock {\em An introduction to {R}iemann-{F}insler geometry}, volume 200 of
  {\em Graduate Texts in Mathematics}.
\newblock Springer-Verlag, New York, 2000.

\bibitem{bellettini:1996}
G.~Bellettini and M.~Paolini.
\newblock Anisotropic motion by mean curvature in the context of {F}insler
  geometry.
\newblock {\em Hokkaido Math. J.}, 25(3):537--566, 1996.

\bibitem{colding:2003}
Tobias~H. Colding and William~P. Minicozzi~II.
\newblock Sharp estimates for mean curvature flow of graphs.
\newblock {\em \tt {arXiv:math.AP/0305099}}, 2003.

\bibitem{eh:mean}
Klaus Ecker and Gerhard Huisken.
\newblock Mean curvature evolution of entire graphs.
\newblock {\em Ann. of Math. (2)}, 130(3):453--471, 1989.

\bibitem{eh:interior}
Klaus Ecker and Gerhard Huisken.
\newblock Interior estimates for hypersurfaces moving by mean curvature.
\newblock {\em Invent. Math.}, 105(3):547--569, 1991.

\bibitem{finn}
Robert Finn.
\newblock On equations of minimal surface type.
\newblock {\em Ann. of Math. (2)}, 60:397--416, 1954.

\bibitem{trudinger}
David Gilbarg and Neil~S. Trudinger.
\newblock {\em Elliptic partial differential equations of second order}.
\newblock Classics in Mathematics. Springer-Verlag, Berlin, 2001.
\newblock Reprint of the 1998 edition.

\bibitem{grayson:plane-curves}
Matthew~A. Grayson.
\newblock The heat equation shrinks embedded plane curves to round points.
\newblock {\em J. Differential Geom.}, 26(2):285--314, 1987.

\bibitem{huisken:spheres}
Gerhard Huisken.
\newblock Flow by mean curvature of convex surfaces into spheres.
\newblock {\em J. Differential Geom.}, 20(1):237--266, 1984.

\bibitem{hu:boundary}
Gerhard Huisken.
\newblock Nonparametric mean curvature evolution with boundary conditions.
\newblock {\em J. Differential Equations}, 77(2):369--378, 1989.

\bibitem{hu:local}
Gerhard Huisken.
\newblock Local and global behaviour of hypersurfaces moving by mean curvature.
\newblock In {\em Differential geometry: partial differential equations on
  manifolds}, volume~54 of {\em Proc. Sympos. Pure Math.}, pages 175--191.
  Amer. Math. Soc., Providence, RI, 1993.

\bibitem{gh:dist}
Gerhard Huisken.
\newblock A distance comparison principle for evolving curves.
\newblock {\em Asian J. Math.}, 2(1):127--133, 1998.

\bibitem{kruzhkov:nonlinear}
S.~N. Kru{\v{z}}kov.
\newblock Nonlinear parabolic equations in two independent variables.
\newblock {\em Trudy Moscov. Mat. Obshch}, (16):329--346, 1967.
\newblock English translation in {\it{Trans. Moscow Math. Soc.}} {\bf{16}}
  (1967), 355-373.

\bibitem{krylov:parabolic}
N.V. Krylov.
\newblock {\em Lectures on Elliptic and Parabolic Equations in H\"older
  Spaces}.
\newblock American Mathematical Society, 1996.

\bibitem{li:first-bvp}
Gary~M. Lieberman.
\newblock The first initial-boundary value problem for quasilinear second order
  parabolic equations.
\newblock {\em Ann. Scuola Norm. Sup. Pisa Cl. Sci. (4)}, 13(3):347--387, 1986.

\bibitem{li:parabolic}
Gary~M. Lieberman.
\newblock {\em Second order parabolic differential equations}.
\newblock World Scientific Publishing Co. Inc., River Edge, NJ, 1996.

\bibitem{yoshi:grad-estimates}
Yuko Nagase and Yoshihiro Tonegawa.
\newblock Interior gradient estimate for 1-d anisotropic curvature flow.
\newblock 2003.
\newblock
  \newline\underline{http://www.math.hokudai.ac.jp/~tonegawa/research.htm}.

\bibitem{simon}
Leon Simon.
\newblock Equations of mean curvature type in {$2$} independent variables.
\newblock {\em Pacific J. Math.}, 69(1):245--268, 1977.

\end{thebibliography}

\markright{   }

\end{document}